\documentclass{siamart1116}

\usepackage{lipsum}
\usepackage{amsfonts}
\usepackage{graphicx}
\usepackage{epstopdf}
\usepackage{algorithmic}
\ifpdf
  \DeclareGraphicsExtensions{.eps,.pdf,.png,.jpg}
\else
  \DeclareGraphicsExtensions{.eps}
\fi

\usepackage{bbm}
\usepackage{bm}
\usepackage{color}
\usepackage{xcolor}
\usepackage{mathdots}
\usepackage{graphicx}
\usepackage{epstopdf}
\usepackage{psfrag}

\usepackage{amssymb}
\usepackage{amsmath}
\usepackage{dsfont}
\usepackage{rotating}
\usepackage{mathrsfs,amsfonts,latexsym}
\usepackage{amssymb}

\newcommand{\diff}[2]{\frac{\partial {#1} }{\partial {#2} } }
\newcommand{\nt}{\tau}	
\newcommand{\p}{\mathbf{p}}
\newcommand{\pan}{p}

\graphicspath{{pix/}}
\usepackage{array}
\newcolumntype{C}[1]{>{\centering\arraybackslash}p{#1}}
\usepackage[font=footnotesize]{caption}
\usepackage{algorithm}
\usepackage{algorithmic}
\usepackage{multirow}

\newtheorem{Theorem}{Theorem}
\theoremstyle{definition}\newtheorem{Definition}{Definition}\newtheorem{Remark}{Remark}

\usepackage[labelformat=simple]{subcaption}

\newcommand{\bx}[1]{{\bf #1}} 

\newcommand{\Mass}[1]{\mathbf{M}_{\textsl{#1}} }
\newcommand{\iMass}[1]{\mathbf{M}^{\mathbf{-1}}_{\textsl{#1}} }
\newcommand{\barr}[2]{\mathcal{#1}_{\textsl{#2}}^\top  }
\newcommand{\tild}[2]{\mathcal{#1}_{\textsl{#2}}  }

\newcommand{\Hoperator}[2]{\mathbb{#1}^{\textsl{#2}}  }
\newcommand{\dof}[1]{\widehat{\mathbf{#1}}}

\renewcommand{\vec}{\mathbf}

\numberwithin{theorem}{section}

\newcommand{\TheTitle}{Staggered discontinuous Galerkin methods for the incompressible Navier-Stokes equations: 
spectral analysis and computational results} 
\newcommand{\TheAuthors}{M. Dumbser,  F. Fambri, I. Furci, M. Mazza, M. Tavelli, and S. Serra-Capizzano }

\headers{Staggered DG methods: spectral analysis and computational results}{M. Dumbser et al.}

\title{{\TheTitle}\thanks{Submitted to the editors DATE.}}

\author{
Michael Dumbser \footnotemark[3]
\and 
  Francesco Fambri \footnotemark[3] 
  	\and 
  Isabella Furci\thanks{Department of  Science and High Technology, University of Insubria, via Valleggio 11, 
	                      I-22100 Como, Italy 
    (\email{ifurci@uninsubria.it}).}
  \and
  Mariarosa Mazza\footnotemark[4] 
	\and
	Stefano Serra-Capizzano \footnotemark[6]	
	\and
	 Maurizio Tavelli \thanks{Laboratory of Applied Mathematics, Department of Civil, Environmental and Mechanical Engineering,
          					   University of Trento, Via Mesiano 77, I-38123 Trento, Italy (\email{m.tavelli@unitn.it}).}
}

\usepackage{amsopn}

\ifpdf
\hypersetup{
  pdftitle={\TheTitle},
  pdfauthor={\TheAuthors}
}
\fi

\begin{document}

\maketitle

\begin{abstract}
The goal of this paper is to create a fruitful bridge between the numerical methods for approximating partial differential equations (PDEs) in fluid dynamics and the (iterative) numerical methods for dealing with the resulting large linear systems. Among the main objectives are the design of new efficient iterative solvers and a rigorous analysis of their convergence speed.
The link we have in mind is either the structure or the hidden structure that the involved coefficient matrices inherit, both from the continuous PDE and from the approximation scheme: in turn, the resulting structure is used for deducing spectral information, crucial for the conditioning and convergence analysis, and for the design of more efficient solvers.

As specific problem we consider the incompressible Navier-Stokes equations, as numerical technique we consider a novel family of high order accurate Discontinuous Galerkin methods on \textit{staggered} meshes, and as tools we use the theory of Toeplitz matrices generated by a function (in the most general block, multi-level form) and the more recent theory of Generalized Locally Toeplitz matrix-sequences. We arrive at a quite complete picture of the spectral features of the underlying matrices and this information is employed for giving a forecast of the convergence history of the conjugate gradient method, together with a discussion on new more advanced techniques (involving preconditioning, multigrid, multi-iterative solvers). Several numerical tests are provided and critically illustrated in order to show the validity and the potential of our analysis.
\end{abstract}



\section{Introduction}
Computational fluid dynamics (CFD) represents a vast sector of ongoing research in engineering and applied mathematics, which has also a wide applicability to real world problems, such as aerodynamics of airplanes and cars, geophysical flows in oceans, lakes and rivers, Tsunami wave propagation, blood flow in the human cardiovascular system, weather forecasting and many others. The governing equations for incompressible fluids are given by the incompressible Navier-Stokes equations that consist in a divergence-free condition for the velocity
\begin{eqnarray}
	\nabla \cdot \vec{v}=0,
\label{eq:intro1}
\end{eqnarray}
and a momentum equation that involves convective, pressure and viscosity effects:
\begin{eqnarray}
	\diff{\vec{v}}{t}+\nabla \cdot \mathbf{F} + \nabla p= \nabla \cdot \left( \nu \nabla \vec{v} \right).
\label{eq:intro2}
\end{eqnarray}
Here, $\vec{v}$ is the velocity field; $p$ is the pressure; $\nu$ is the kinematic viscosity coefficient and $\mathbf{F}=\vec{v}\otimes \vec{v}$ is the tensor containing the nonlinear convective term.
The dynamics induced  by equations \eqref{eq:intro1}-\eqref{eq:intro2} can be rather complex and has been observed in various  experiments, see \cite{Armaly1983,Sakamoto1990,Williamson1988}.
 
In the last decades a lot of effort was made to numerically solve the incompressible Navier-Stokes equations using finite difference schemes (see \cite{markerandcell,patankarspalding,patankar,vanKan}), continuous finite elements (see \cite{TaylorHood,SUPG,SUPG2,Fortin,Verfuerth,Rannacher1,Rannacher3}) and more recently high order Discontinuous Galerkin (DG) methods, see e.g. \cite{Bassi2006,Bassi2007,Shahbazi2007,Ferrer2011,Nguyen2011,Rhebergen2012,Rhebergen2013,Crivellini2013,KleinKummerOberlack2013}.

The main difficulty in the numerical solution of the incompressible Navier-Stokes equations \eqref{eq:intro1}-\eqref{eq:intro2} lies in the elliptic pressure Poisson equation and the associated linear equation system that needs to be solved. On the discrete level the pressure system is obtained by substitution of the discrete momentum equation \eqref{eq:intro2} into the discrete form of the divergence-free condition \eqref{eq:intro1}.

Since the solution of the incompressible Navier-Stokes equations requires necessarily the solution of large systems of algebraic equations, it is indeed very important 
to have a scheme that uses a stencil that is as small as possible, in order to improve the sparsity pattern of the resulting system matrix. It is also desirable to use methods that lead to reasonably well conditioned systems that can be solved with iterative solvers, like the conjugate gradient (CG) method \cite{cgmethod} or the GMRES algorithm \cite{GMRES}. 

Very recently, a new class of arbitrary high order accurate semi-implicit DG schemes for the solution of the incompressible Navier-Stokes equations on structured and unstructured edge-based \textit{staggered} grids was proposed in \cite{Fambri2016, 2SINS, 2STINS, 3STINS}, following a philosophy that had been first introduced in finite difference schemes, see \cite{markerandcell,patankarspalding,patankar,vanKan,HirtNichols,CasulliCheng1992,Casulli1999,CasulliWalters2000,Casulli2009,CasulliStelling2011,CasulliVOF}. 
All those approaches have in common that the pressure is defined on a main grid, while the velocity field is defined on an appropriate edge-based staggered grid.
The nonlinear convective terms are discretized explicitly by using a standard DG scheme based on the local Lax-Friedrichs (Rusanov) flux \cite{Rusanov:1961a}.
Then, the discrete momentum equation is inserted into the discrete continuity equation in order to obtain the discrete form of the pressure Poisson equation.

The advantage in using staggered grids 
 is that they allow to improve significantly the sparsity pattern of the final linear system that has to be solved for the pressure. 
For the structured case the resulting main linear system is a sparse block penta-diagonal and hepta-diagonal one in two and three space dimensions, respectively.
Furthermore, several desirable properties, such as the symmetry and the positive definiteness can be achieved see e.g. \cite{Fambri2016, 3STINS}.

The regular shape of the structured case allows to further describe the structure of the main linear system for the pressure in the framework of multi-level block Toeplitz matrices:
in this setting we can deliver spectral and computational
properties, including specific preconditioners and specific multigrid methods for the preconditioning matrices.

The rest of this paper is organized as follows. Section \ref{sec:overv} is devoted to a brief overview of the numerical methods used in this paper for the solution of the incompressible  Navier-Stokes equations.
Section \ref{MainSec:SpectralAnalysis} studies the linear systems stemming from the considered approximations in a setting of structured linear algebra: by using known properties of multi-level block Toeplitz (and circulant) matrices, we are able 
to provide a detailed structural and spectral analysis of the involved matrices, including conditioning, extremal eigenvalues, and spectral distribution results. In Section \ref{sec:alg-num} the spectral features are used for proposing
specific (preconditioned) Krylov methods with a study of the complexity and of the convergence speed: several numerical experiments are reported and critically discussed. Finally, Section \ref{sec:final} deals with conclusions, open problems, and future lines of research.

\section{Overview}\label{sec:overv}

\begin{figure}[htb]
\begin{subfigure}[c]{.47\textwidth}
\includegraphics[width=1.0\textwidth]{Staggering-eps-converted-to.pdf}
\subcaption{Two dimensional case }
\end{subfigure}
\begin{subfigure}[c]{.47\textwidth}
\includegraphics[width=1.0\textwidth]{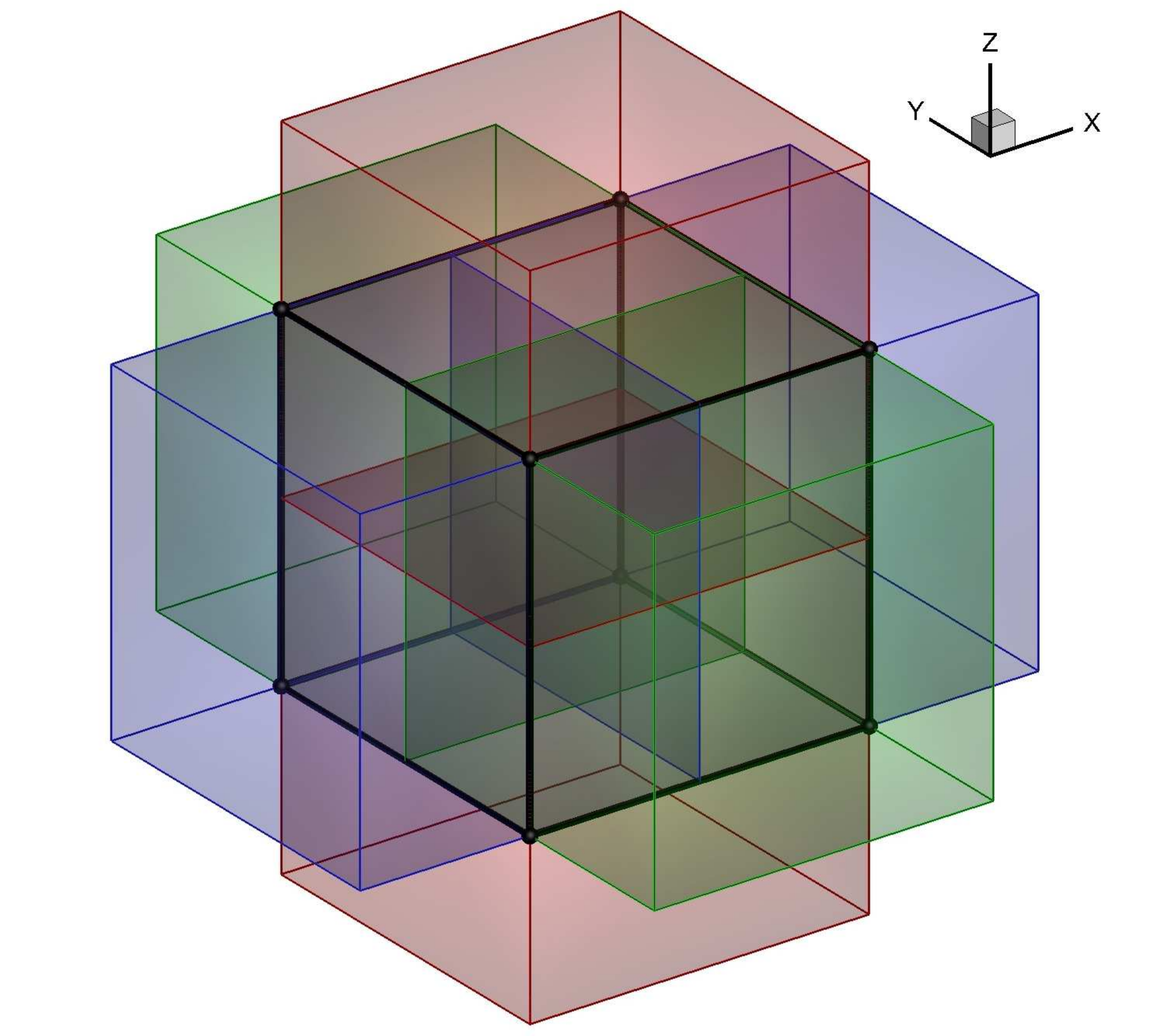}
\subcaption{Three dimensional case }
\end{subfigure}
\caption{Mesh-staggering for the two dimensional case (left) and for the three-dimensional case (right).}\label{fig:staggering}
\end{figure}

In the framework of high order semi-implicit \textit{staggered} discontinuous Galerkin schemes for the incompressible Navier-Stokes equations, the numerical solution for the velocity $\vec{v}=(u,v,w)$ and the pressure $p$ is represented by piecewise polynomials on overlapping staggered grids.  
The numerical solution can be written as a linear combination of polynomial basis functions, i.e. $p_h(\vec{x},t)= \sum_l \mathbf{\phi}_l(\vec{x}) \, \hat{p}_l(t)$ and $\vec{v}_h(\vec{x},t)=\sum_l \mathbf{\psi}_l(\vec{x}) \, \hat{\vec{v}}_l(t)$. Here, $\phi_l$ represents the vector of piecewise polynomial basis functions computed in $\vec{x}$ on the \textit{main grid}, while $\mathbf{\psi}_l$ are the basis functions on the edge-based staggered dual grid; the $\hat{\vec{v}}_l$ and $\hat{p}_l$ are the vectors of the so called \emph{degrees of freedom} associated with the discrete solution $\vec{v}_h$ and $p_h$, respectively. The chosen staggered grid is an \textit{edge based} staggering, corresponding to the one used in \cite{DumbserCasulli}. The staggering of the flow quantities is briefly depicted in Figure 
$\ref{fig:staggering}$, where also the main indexing used for the numerical solution is reported, together with fractional indices referring to staggered grids.

We can derive a weak formulation of the momentum and continuity equation in the form 
\begin{align}
\Mass{xyz} \cdot\left( \dof{u}^{\nt+\delta\tau}_{i+\frac{1}{2},j,k}- \dof{Fu}^{\nt}_{i+\frac{1}{2},j,k}  \right) + \frac{\delta \tau}{\Delta x}\Mass{yz}\cdot \left( \tild{R}{x} \cdot\dof{p}^{\nt+\delta\tau}_{i+1,j,k} - \tild{L}{x} \cdot\dof{p}^{\nt+\delta\tau}_{i,j,k} \right)&=0,\label{eq:uSDG1} \\
\Mass{xyz}\cdot \left( \dof{v}^{\nt+\delta\tau}_{i,j+\frac{1}{2},k}- \dof{Fv}^{\nt}_{i,j+\frac{1}{2},k}  \right) + \frac{\delta \tau}{\Delta y}\Mass{zx} \cdot\left( \tild{R}{y} \cdot\dof{p}^{\nt+\delta\tau}_{i,j+1,k} - \tild{L}{y}\cdot \dof{p}^{\nt+\delta\tau}_{i,j,k} \right)&=0,\label{eq:vSDG1} \\
\Mass{xyz} \cdot\left( \dof{w}^{\nt+\delta\tau}_{i,j,k+\frac{1}{2}}- \dof{Fw}^{\nt}_{i,j,k+\frac{1}{2}}  \right) + \frac{\delta \tau}{\Delta z} \Mass{xy}\cdot\left( \tild{R}{z}\cdot \dof{p}^{\nt+\delta\tau}_{i,j,k+1} - \tild{L}{z} \cdot\dof{p}^{\nt+\delta\tau}_{i,j,k} \right)&=0,\label{eq:wSDG1} 
\end{align}
\begin{eqnarray}
 \frac{\Mass{yz} \left(\barr{L}{x} \cdot\dof{u}^{\nt+\delta\tau}_{i+\frac{1}{2},j,k} \!-\! \barr{R}{x} \cdot\dof{u}^{\nt+\delta\tau}_{i-\frac{1}{2},j,k}\right)}{\Delta x} + 
 \frac{\Mass{zx} \left(\barr{L}{y} \cdot\dof{v}^{\nt+\delta\tau}_{i,j+\frac{1}{2},k} \!-\! \barr{R}{y}\cdot \dof{v}^{\nt+\delta\tau}_{i,j-\frac{1}{2},k}\right)}{\Delta y} +  && \nonumber \\ 
 \frac{\Mass{xy} \left(\barr{L}{z} \cdot\dof{w}^{\nt+\delta\tau}_{i,j,k+\frac{1}{2}} \!-\! \barr{R}{z} \cdot\dof{w}^{\nt+\delta\tau}_{i,j,k-\frac{1}{2}}\right)}{\Delta z} = 0, && \label{eq:pSDG1}
\end{eqnarray}
where $\Mass{}$ is the so called \emph{mass matrix}, $\tild{R}{}$ and $\tild{L}{}$ are some real valued matrices related to the discrete form of the gradient operator, see \cite{Fambri2016} for the complete definitions, $\Delta x$, $\Delta y$, $\Delta z$ and $\delta \tau$ are the space and time step size.  
Formal substitution of the implicit velocities $[\dof{u}^{\nt+\delta\tau}_\cdot,\dof{v}^{\nt+\delta\tau}_\cdot,\dof{w}^{\nt+\delta\tau}_\cdot]$ given in equations \eqref{eq:uSDG1}-\eqref{eq:wSDG1} into \eqref{eq:pSDG1} leads to a linear system for the new pressure $\dof{p}^{\nt+\delta\tau}_\cdot$ that reads
\begin{align}
 \frac{ \delta \tau}{\Delta x^2}\left(\Mass{yz}\cdot \Hoperator{R}{x} \right)\cdot\dof{p}^{\nt+\delta\tau}_{i+1,j,k} + \frac{ \delta \tau}{\Delta y^2}\left(\Mass{zx} \cdot \Hoperator{R}{y} \right)\cdot\dof{p}^{\nt+\delta\tau}_{i,j+1,k}+ \frac{ \delta \tau}{\Delta z^2}\left(\Mass{xy}\cdot  \Hoperator{R}{z} \right)\cdot\dof{p}^{\nt+\delta\tau}_{i,j,k+1} + \nonumber \\
+ \left( \frac{ \delta \tau}{\Delta x^2}\Mass{yz}\cdot  \Hoperator{C}{x} + \frac{\delta \tau}{\Delta y^2}\Mass{zx} \cdot \Hoperator{C}{y}+ \frac{ \delta \tau}{\Delta z^2}\Mass{xy} \cdot \Hoperator{C}{z}\right)\cdot \dof{p}^{\nt+\delta\tau}_{i,j,k} + \nonumber \\
 + \frac{ \delta \tau}{\Delta x^2}\left(\Mass{yz}\cdot  \Hoperator{L}{x} \right)\cdot\dof{p}^{\nt+\delta\tau}_{i-1,j,k} + \frac{\delta \tau}{\Delta y^2}\left(\Mass{zx} \cdot \Hoperator{L}{y} \right)\cdot\dof{p}^{\nt+\delta\tau}_{i,j-1,k}+ \frac{ \delta \tau}{\Delta z^2}\left(\Mass{xy} \cdot \Hoperator{L}{z} \right)\cdot\dof{p}^{\nt+\delta\tau}_{i,j,k-1}  \nonumber \\ 
 = \dof{b}^\nt_{i,j,k},\label{eq:pSyst} \\
\text{for}\;\;\;i=2,...,n_1-1;\;\;\;j=2,...,n_2-1;\;\;\,k=2,...,n_3-1  \nonumber 
\end{align}
where
\begin{align}
\Hoperator{R}{} = - \left(\barr{L}{}\cdot\iMass{}\cdot\tild{R}{} \right),\;\;\; 
\Hoperator{L}{} = - \left(\barr{R}{}\cdot\iMass{}\cdot\tild{L}{}\right),\;\;\;   \nonumber \\ 
\Hoperator{C}{} =  \left(\barr{L}{}\cdot\iMass{}\cdot\tild{L}{}\right) + \left(\barr{R}{}\cdot\iMass{}\cdot\tild{R}{} \right).
\end{align}
and $n_1,n_2,n_3$ are the total number of elements in the $x$, $y$ and $z$ direction, respectively.
System ($\ref{eq:pSyst}$) is then written in compact form as $K_N \cdot \mathbf{p}^{\nt+\delta\tau} = b^\nt$. Here $\mathbf{p}^{\nt+\delta\tau}$ collects all the unknown pressure degrees of freedom at the new time step $\nt+\delta\tau$ and $b^\nt$ contains all the terms known at the time step $\nt$, see again \cite{Fambri2016} for more details. In particular, in \cite{Fambri2016} it has been shown that the resulting linear system is \textit{symmetric}. Furthermore, it is clear from system $\eqref{eq:pSyst}$ that the stencil involves only the direct neighbors, and hence it is a symmetric $7$ block-diagonal system for the $3D$ case and a $5$ block-diagonal system for the $2D$ case.

Once the new pressure $\mathbf{p}^{\nt+\delta\tau}$ is known, we can readily compute the new velocity field $[\dof{u}^{\nt+\delta\tau}_\cdot,\dof{v}^{\nt+\delta\tau}_\cdot,\dof{w}^{\nt+\delta\tau}_\cdot]$ from equations \eqref{eq:uSDG1}-\eqref{eq:wSDG1}.


\section{Spectral analysis}
\label{MainSec:SpectralAnalysis}

This section is devoted to the structural and spectral analysis of the linear systems arising from the staggered DG approximation of incompressible two-dimensional incompressible Navier-Stokes equations, with special attention to the following items:
\begin{itemize}
\item structural properties, in connection with multi-level block Toeplitz (and circulant) matrices,
\item distribution spectral analysis in the Weyl sense, 
\item conditioning and asymptotic behaviour of the extremal eigenvalues.
\end{itemize}

In particular, the first item is used for the second two, which in turn are of interest in the analysis of the intrinsic difficulty of the problem and in the design and convergence analysis of (preconditioned) Krylov methods \cite{AxL,BeKu}.

\subsection{Problem setting}
Our aim is to efficiently solve large linear systems arising from the staggered DG approximation of incompressible two-dimensional Navier-Stokes equations taking advantage of the structure of the coefficient matrix and especially of its spectral features. More precisely, when discretizing the problem of interest for a sequence of discretization parameters $h_N$ we obtain a sequence of linear systems, in which the $N$-th component is of the form
\begin{equation}\label{system}
K_Nx=b, \quad K_N\in\mathbb{R}^{N\times N}, \quad x,b\in\mathbb{R}^{N},
\end{equation}
whose coefficient matrix size $N$ grows to infinity as the approximation error tends to zero. In order to analyze standard methods and for designing new efficient solvers for the considered linear systems, it is of crucial importance to have a spectral analysis of the matrix-sequence $\{K_N\}_N$. As we will show in the next sections, the coefficient matrix $K_N$ is, up to low-rank perturbations, a $2$-level block Toeplitz matrix: however, when considering variable coefficients or for the study of the preconditioning, standard Toeplitz structures are not sufficient. For this reason, we need to introduce the notion of multi-level block-Toeplitz sequences associated with a matrix-valued symbol and of Generalized Locally Toeplitz (GLT) algebra.

\subsection{Background and definitions}
Throughout this paper, we use the following notation. Let $\mathcal{M}_s$ be the linear space of the complex $s\times s$ matrices and let $f:G\to\mathcal M_s$, with $G\subseteq\mathbb R^\ell$, $\ell\ge 1$, measurable set. We say that $f$ belongs to $L^p(G)$ (resp. is measurable) if all its components $f_{ij}:G\to\mathbb C,\ i,j=1,\ldots,s,$ belong to $L^p(G)$ (resp. are measurable) for $1\le p\le\infty$. Moreover, we denote by $\mathcal{I}_k$ the $k$-dimensional
cube $(-\pi,\pi)^k$ and define $L^p(k,s)$ as the linear space of $k$-variate functions $f:\mathcal{I}_k\to\mathcal M_s$, $f\in
L^p(\mathcal{I}_k)$. Let ${\bf n}:=(n_1,\ldots,n_k)$ be a multi-index
in $\mathbb N^k$ and set $\hat n:=\prod_{i=1}^kn_i$.

\begin{Definition}\label{def-multilevel}
Let the Fourier coefficients of a given function $f\in L^1(k,s)$ be defined as
\begin{align}\label{fhat}
  \hat f_{\bf j}:=\frac1{(2\pi)^k}\int_{\mathcal{I}_k}f(\theta){\rm e}^{-{\iota}\left\langle {\bf j},\theta\right\rangle}\ d\theta\in\mathcal M_s,
  \qquad {\bf j}=(j_1,\ldots,j_k)\in\mathbb Z^k,
\end{align}
where $\left\langle { \bf j},\theta\right\rangle=\sum_{t=1}^kj_t\theta_t$ and the integrals in \eqref{fhat} are computed componentwise. \\
Then, the ${\bf n}$th
Toeplitz matrix associated with $f$ is the matrix of order $s\hat n$ given by

\begin{equation}
T_{\bf n}(f) =\sum_{\bf j=-(\bf n-\bf e)}^{\bf n-\bf e} J_{n_1}^{j_1} \otimes \cdots\otimes J_{n_k}^{j_k}\otimes \hat{f}_{\bf j}
\end{equation}
and the ${\bf n}$th
circulant matrix associated with $f$ is the matrix of order $s\hat n$ given by

\begin{equation}
C_{\bf n}(f) =\sum_{\bf j=-(\bf n-\bf e)}^{\bf n-\bf e} Z_{n_1}^{j_1} \otimes \cdots\otimes Z_{n_k}^{j_k}\otimes \hat{f}_{\bf j}
\end{equation} 
where $\bx{e}=(1,\ldots,1)\in\mathbb{N}^k, \,\bx{j}=(j_1,\ldots,j_k)\in\mathbb{N}^k$ and $ Z^{j_{\xi}}_{n_{\xi}}$ (resp. $ J^{j_{\xi}}_{n_\xi}$) is the $n_{\xi} \times n_{\xi}$ matrix whose $(i,h)$-th entry equals 1 if $(i-h)$mod$n_{\xi}=j_{\xi}$ (resp. if $(i-h)=j_{\xi}$) and $0$ otherwise.\\
The sets {\color{blue}}$\{T_{\bf n}(f)\}_{{\bf n}\in\mathbb N^k}$ and $\{C_{\bf n}(f)\}_{{\bf n}\in\mathbb N^k}$ are 
called the \emph{family of $k$-level Toeplitz matrices and $k$-level circulant matrices respectively, generated by $f$}, that in turn is referred to as the \emph{generating function or the symbol of either
$\{T_{\bf n}(f)\}_{{\bf n}\in\mathbb N^k}$ or $\{C_{\bf n}(f)\}_{{\bf n}\in\mathbb N^k}$}.
\end{Definition}

In order to deal with low-rank perturbations and to show that they do not
affect the symbol of a Toeplitz sequence, we need first to introduce the definition of spectral
distribution in the sense of the eigenvalues and of the singular values for a
generic matrix-sequence $\{A_{\bf n}\}_{{\bf n}\in\mathbb{N}^v}$, $v\ge 1$, and then the notion of 
GLT algebra. In short, the latter is an algebra
containing sequences of matrices including the Toeplitz sequences with Lebesgue
integrable symbols and virtually any sequence of matrices coming from
`reasonable' approximations by local discretization methods (finite differences,
finite elements, isogeometric analysis, etc.) of partial differential equations.

\begin{Definition}\label{def-distribution}
Let $f:G\to\mathcal{M}_s$ be a measurable function, defined on a measurable set $G\subset\mathbb R^\ell$ with $\ell\ge 1$,
$0<m_\ell(G)<\infty$. Let $\mathcal C_0(\mathbb K)$ be the set of continuous functions with compact support over $\mathbb
K\in \{\mathbb C, \mathbb R_0^+\}$ and let $\{A_{\bf n}\}_{{\bf n}\in\mathbb{N}^v}$, $v\ge 1$, be a sequence of matrices with
eigenvalues $\lambda_j(A_{\bf n})$, $j=1,\ldots,N$ and singular
values $\sigma_j(A_{\bf n})$, $j=1,\ldots,N$, where $N\equiv N({\bf n})$ is the size of $A_{\bf n}$ and has to be a monotonic function with respect to every single variable $n_i$, $i=1,\ldots,v$.

\begin{itemize}
    \item $\{A_{\bf n}\}_{{\bf n}\in\mathbb{N}^v}$ is {\em distributed as the pair
    	  $(f,G)$ in the sense of the eigenvalues,} in symbols
    	  $$\{A_{\bf n}\}_{{\bf n}\in\mathbb N^v}\sim_\lambda(f,G),$$ if the following limit relation holds for all $F\in\mathcal C_0(\mathbb C)$:
		\begin{align}\label{distribution:sv-eig}
		  \lim_{{\bf n}\to\infty}\frac{1}{N}\sum_{j=1}^{N}F(\lambda_j(A_{\bf n}))=
		  \frac1{m_\ell(G)}\int_G \frac{{\rm tr}(F(f(\theta)))}{s} d\theta.
		\end{align}
    \item $\{A_{\bf n}\}_{{\bf n}\in\mathbb N^v}$ is {\em distributed as the pair
    	  $(f,G)$ in the sense of the singular values,} in symbols
    	  $$\{A_{\bf n}\}_{{\bf n}\in\mathbb N^v}\sim_\sigma(f,G),$$ if the following
    	  limit relation holds for all $F\in\mathcal C_0(\mathbb R_0^+)$:
		\begin{align}\label{distribution:sv-eig-bis}
		  \lim_{{\bf n}\to\infty}\frac{1}{N}\sum_{j=1}^{N}F(\sigma_j(A_{\bf n}))=
		  \frac1{m_\ell(G)}\int_G \frac{{\rm tr}(F(|f(\theta)|))}{s} d\theta.
		\end{align}
\end{itemize}
In this setting the expression ${{\bf n}\to\infty}$ means that every component of the vector ${\bf n}$ tends to infinity, that is,
$\min_{i=1,\dots,v} n_i\to\infty$, while $|f(\theta)|=\left(f(\theta)f^*(\theta)\right)^{1/2}$ with $^*$ meaning transpose conjugate. 
\end{Definition}

\begin{Remark}\label{rem_distr}
Denote by $\lambda_1(f), \ldots,\lambda_s(f)$ and by $\sigma_1(f),
\ldots,\sigma_s(f)$ the eigenvalues and the singular values of a
$s\times s$ matrix-valued function $f$, respectively. If $f$ is
smooth enough, an informal interpretation of the limit relation
\eqref{distribution:sv-eig} (resp. \eqref{distribution:sv-eig-bis})
is that when the matrix-size of $A_{\bf n}$ is sufficiently large,
then $N/s$ eigenvalues (resp. singular values) of $A_{\bf n}$ can
be approximated by a sampling of $\lambda_1(f)$ (resp. $\sigma_1(f)$)
on a uniform equispaced grid of the domain $G$, and so on until the
last $N/s$ eigenvalues which can be approximated by an equispaced sampling
of $\lambda_s(f)$ (resp. $\sigma_s(f)$) in the domain.

For example, 
take $G$ any domain as in Definition \ref{def-distribution} and let $F=\chi_{[a,b]}(\cdot)$ for a fixed real interval $[a,b]$ such that 
\begin{equation}\label{hp:misura zero}
 m_\ell \left\{\theta \in G  \, : \,\lambda_r(f(\theta))=a\right\}= m_\ell \left\{\theta \in G  \, : \,\lambda_r(f(\theta))=b\right\}=0
 \end{equation}
for every $r=1,\ldots,s$. 
 Note that $F=\chi_{[a,b]}(\cdot)$ is a discontinuous function, but, under the assumptions in (\ref{hp:misura zero}), the limit relation \eqref{distribution:sv-eig} still holds true. The argument of the proof relies in choosing two families of continuous approximations $\{F^{-}_\delta\}_\delta$, $\{F^{+}_\delta\}_\delta$ of $\chi_{[a,b]}$ such that $F^{+}_\delta<\chi_{[a,b]} <F^{+}_\delta$ (see \cite{serra1998} for more details).
If we define
\begin{equation*}
m_r={\rm essinf}_{G} \lambda_r(f(\theta)),\qquad M_r={\rm esssup}_{G}\lambda_r(f(\theta)), \qquad r=1,\dots ,s,
\end{equation*}
when $F=\chi_{[m_r,M_r]}(\cdot)$, then equation \eqref{distribution:sv-eig} becomes
 \begin{equation}
 \lim_{n \to \infty}\frac{1}{N} \sum _{i=1}^{N} \chi_{[m_r,M_r]}\left(\lambda_i(A_{\textbf{n}})\right)=\frac{1}{s m_\ell(G)} \int_{G} {\rm tr}\left(\chi_{[m_r,M_r]}(f(\theta))\right) \, d\theta,
 \label{conta_autovalori}
 \end{equation}
and hence
\begin{equation}\label{conta}
\begin{split}
\lim_{N \to \infty}\frac{1}{N} &\# \left\{ i \, : \, \lambda_i(A_{\bf n}) \in [m_r,M_r]\right\}=\\
 &\frac{1}{s m_\ell(G)} \sum _{j=1}^{s} m_\ell \left\{\theta \in G  \, : \,\lambda_j(f(\theta)) \in [m_r,M_r]\right\}.
\end{split}
 \end{equation}
Moreover, if
\begin{equation*}
{\rm esssup}_{G}\lambda_r (f(\theta) )\le{\rm essinf}_{G} \lambda_{r+1} (f(\theta)), \quad r=1,\ldots, s-1,
\end{equation*}
then equation \eqref{conta} in turn becomes
\begin{equation*}
  \begin{split}
 \lim_{n \to \infty}\frac{1}{N} &\# \left\{ i \, : \, \lambda_i(A_{\bf n}) \in [m_r,M_r]\right\}=\\
 &\frac{1}{s m_\ell(G)}  m_\ell \left\{\theta \in G  \, : \,\lambda_r(f(\theta)) \in [m_r,M_r]\right\}=\frac{1}{s}
 \end{split}
 \end{equation*}
which means that
\begin{equation*}
\# \left\{ i \, : \, \lambda_i(A_{\bf n}) \in [m_r,M_r]\right\}=\frac{N}{s}+o(N).
\end{equation*}

In the Toeplitz setting, when $f$ is a $k$-variate polynomial, the quantity $o(N)$ becomes proportional to
$N^{1-\frac{1}{k}}$, with constant proportional to $s$ and to the degree of the polynomial, and with $N=\hat n$. 

\end{Remark}

\subsection{Spectral analysis of Hermitian (block)  Toeplitz sequences: distribution results}

Concerning the spectral distribution of Toeplitz sequences, if $f$ is a real-valued function, then the following theorem holds:
 Szeg\"o stated this result for $f$ essentially bounded, while the extension to Lebesgue integrable generating function is due to
Tyrtyshnikov and Zamarashkin \cite{TyZ}.

\begin{Theorem}[\cite{GSz}]\label{szego}
Let $f\in L^1(k,1)$ be a real-valued function with $k\ge 1$. Then, $$\left\{T_{\bf n}(f)\right\}_{{\bf n}\in\mathbb N^k}\sim_\lambda(f,\mathcal{I}_k).$$
\end{Theorem}

In the case where $f$ is a Hermitian matrix-valued function, the previous theorem can be extended as follows:

\begin{Theorem}[\cite{Tillinota}]\label{szego-herm}
Let $f\in L^1(k,s)$ be a Hermitian matrix-valued function with $k\ge 1,s\ge 2$. Then, $$\{T_{\bf n}(f)\}_{{\bf n}\in\mathbb N^k}\sim_\lambda~(f,\mathcal{I}_k).$$
\end{Theorem}

\begin{Remark}\label{rem_simmetria}
If $\{T_{\bf n}(f)\}_{{\bf n}\in\mathbb N^k}$ is such that each $T_{\bf n}(f)$ is symmetric with symmetric and real blocks, then the symbol has the additional property that
$f(\pm\theta_1,\ldots,\pm\theta_k)\equiv f(\theta_1,\ldots,\theta_k)$, $\forall(\theta_1,\ldots,\theta_k)\in \mathcal{I}_k^+=
[0,\pi]^k$ and therefore Theorem \ref{szego-herm} can be rephrased as
$$\{T_{\bf n}(f)\}_{{\bf n}\in\mathbb N^k}\sim_\lambda~(f,\mathcal{I}_k^+).$$

\end{Remark}

\subsection{Spectral analysis of Hermitian (block) Toeplitz sequences: extremal eigenvalues}

Concerning the localization and the extremal behaviour of the spectra of Toeplitz sequences, if $f$ is a real-valued function, then we have the following result.

\begin{Theorem}[\cite{serra98, Grudsky}]\label{loc-extr}
Let $f\in L^1(k,1)$ be a real-valued function with $k\ge 1$. Let $m$ be the essential infimum of $f$ and $M$ be the essential supremum of $f$. 
\begin{enumerate}
\item If $m=M$ then $f$ is the constant $m$ a.e. and $T_{\bf n}(f)$ coincides with $m$ times the identity of size $\hat n$.
\item If $m<M$ then all the eigenvalues of $T_{\bf n}(f)$ belong to the open set $(m,M)$ for every ${\bf n}\in\mathbb N^k$.
\item If $m=0$ and $\tilde \theta$ is the unique zero of $f$ such that there exist positive constants $c,C,\alpha$ for which
\[ 
c\|\theta -\tilde \theta\|^\alpha \le f(\theta) \le C \|\theta -\tilde \theta\|^\alpha,
\]
then the minimal eigenvalue of $T_{\bf n}(f)$ goes to zero as $(\hat n)^{-\alpha/k}$.
\end{enumerate}
\end{Theorem}

In the case where $f$ is a Hermitian matrix-valued function, the previous theorem can be extended as follows:

\begin{Theorem}[\cite{serra1998, serra99}]\label{loc-extr-s}
Let $f\in L^1(k,s)$ be a Hermitian matrix-valued function with $k\ge 1,s\ge 2$. Let $m_1$ be the essential infimum of 
the minimal eigenvalue of $f$, $M_1$ be the essential supremum of the minimal eigenvalue of $f$, $m_s$ be the essential infimum of 
the maximal eigenvalue of $f$, and $M_s$ be the essential supremum of the maximal eigenvalue of $f$. 
\begin{enumerate}
\item If $m_1=M_s$ then $f$ is the constant $m I_s$ a.e. and $T_{\bf n}(f)$ coincides with $m$ times the identity of size $s\hat n$.
\item If $m_1<M_1$ then all the eigenvalues of $T_{\bf n}(f)$ belong to the open set $(m_1,M_s]$ for every ${\bf n}\in\mathbb N^k$. If $m_s<M_s$ then all the eigenvalues of $T_{\bf n}(f)$ belong to the open set $[m_1,M_s)$ for every ${\bf n}\in\mathbb N^k$.
\item If $m_1=0$ and $\tilde \theta$ is the unique zero of $\lambda_{\min}(f)$ such that there exist positive constants $c,C,\alpha$ for which
\[ 
c\|\theta -\tilde \theta\|^\alpha \le \lambda_{\min}(f(\theta)) \le C \|\theta -\tilde \theta\|^\alpha,
\]
then the minimal eigenvalue of $T_{\bf n}(f)$ goes to zero as $(\hat n)^{-\alpha/k}$.
\end{enumerate}
\end{Theorem}

\subsection{Spectral analysis and computational features of (block) circulant matrices}
In this subsection we report key features of the (block) circulant matrices, also in connection with the generating function.
\begin{Theorem}[\cite{davis}]\label{circ}
Let $f\in L^1(k,1)$ be a complex-valued function with $k\ge 1$. Then, the following (Schur) decomposition of 
$C_{\bf n}(f)$ is valid:
\begin{equation}
\label{schur}
C_{\bf n}(f)=F_{\bf n}  D_{\bf n}(f) F^*_{\bf n}, 
\end{equation}
where 
\begin{equation}\label{eig-circ}
  D_{\bf n}(f) ={\rm diag}_{ {\bf 0}\le {\bf r}\le {\bf n}-{\bf e}}\left(S_{\bf n}(f)\left(\theta_{\bf r}^{({\bf n})}\right)\right), \quad \theta_{\bf r}^{({\bf n})}=2\pi\frac{ {\bf r}}{{\bf n}}, \quad F_{\bf n}=\frac{1}{\sqrt{\hat n}} \left({\rm e}^{-{\iota}\left\langle {\bf j},\theta_{\bf r}^{({\bf n})}\right\rangle}\right)_{\bx{j},\bx{r}=\bx{0}}^{\bx{n-e}}
\end{equation}
with $\left\langle {\bf j},\theta_{\bf r}^{({\bf n})}\right\rangle=\sum_{t=1}^k 2\pi \frac{j_tr_t}{n_t} $. Here $S_{\bf n}(f)(\cdot)$ is the $\bf n$-th Fourier sum of $f$ given by
\begin{equation}\label{fourier-sum}
S_{\bf n}(f)(\theta) = \sum_{j_1=1-n_1}^{n_1-1} \cdots  \sum_{j_k=1-n_k}^{n_k-1} \hat f_{\bf j} 
{\rm e}^{{\iota}\left\langle {\bf j},\theta\right\rangle}, \ \ \ \  \left\langle {\bf j},\theta\right\rangle=\sum_{t=1}^k  j_t \theta_t.
\end{equation}
Here $F_{\bf n}$ is the $k$-level Fourier matrix, $F_{\bf n}=F_{n_1}\otimes \cdots \otimes F_{n_k}$, and its columns are the eigenvectors of $C_{\bf n}(f)$ with eigenvalues given by the evaluations of the $\bf n$-th Fourier sum $S_{\bf n}(f)(\cdot)$ at the grid points
\[
\theta_{\bf r}^{({\bf n})}=2\pi\frac{ {\bf r}}{{\bf n}}.
\] 
\end{Theorem}

In the case where $f$ is a Hermitian matrix-valued function, the previous theorem can be extended as follows:

\begin{Theorem}[\cite{garoniserrasesana15}]\label{circ-s}
Let $f\in L^1(k,s)$ be a matrix-valued function with $k\ge 1,s\ge 2$. Then, the following (block-Schur) decomposition of 
$C_{\bf n}(f)$ is valid:
\begin{equation}
\label{schur-s}
C_{\bf n}(f)=(F_{\bf n}\otimes I_s) D_{\bf n}(f) (F_{\bf n}\otimes I_s)^*, 
\end{equation}
where 
\begin{equation}\label{eig-circ-s}
  D_{\bf n}(f) ={\rm diag}_{ {\bf 0}\le {\bf r}\le {\bf n}-{\bf e}}\left(S_{\bf n}(f)\left(\theta_{\bf r}^{({\bf n})}\right)\right), \quad \theta_{\bf r}^{({\bf n})}=2\pi\frac{ {\bf r}}{{\bf n}}, \quad F_{\bf n}=\frac{1}{\sqrt{\hat n}} \left({\rm e}^{-{\iota}\left\langle {\bf j},\theta_{\bf r}^{({\bf n})}\right\rangle}\right)_{\bx{j},\bx{r}=\bx{0}}^{\bx{n-e}}
\end{equation}
with $\left\langle {\bf j},\theta_{\bf r}^{({\bf n})}\right\rangle=\sum_{t=1}^k 2\pi \frac{j_tr_t}{n_t} $ and $I_s$ the $s \times s$ identity matrix. Here $S_{\bf n}(f)(\cdot)$ is the $\bf n$-th Fourier sum of $f$ given by
\begin{equation}\label{fourier-sum-s}
S_{\bf n}(f)(\theta) = \sum_{j_1=1-n_1}^{n_1-1} \cdots  \sum_{j_k=1-n_k}^{n_k-1} \hat f_{\bf j} 
{\rm e}^{{\iota}\left\langle {\bf j},\theta\right\rangle}, \ \ \ \  \left\langle {\bf j},\theta\right\rangle=\sum_{t=1}^k  j_t \theta_t.
\end{equation}
Here the eigenvalues of $C_{\bf n}(f)$ are given by the evaluations of $\lambda_t(S_{\bf n}(f)(\cdot))$, $t=1,\ldots,s$, at the grid points
\[
\theta_{\bf r}^{({\bf n})}=2\pi\frac{ {\bf r}}{{\bf n}}.
\] 
\end{Theorem}

\begin{Remark}\label{fourier-vs-funzione}
If $f$ is a trigonometric polynomial of fixed degree (with respect to $\bf n$), then it is worth noticing that
 $S_{\bf n}(f)(\cdot) = f(\cdot)$ for $\bf n$ large enough: more precisely, every $n_j$ should be larger than the double of the degree with respect to the $j$-th variable.
Therefore, in such a setting, the eigenvalues of $C_{\bf n}(f)$ are either the evaluations of $f$ at the grid points if $s=1$ or
the evaluations of $\lambda_t(f(\cdot))$, $t=1,\ldots,s$, at the very same grid points.
\end{Remark}

\begin{Remark}\label{distrib-circ}
Thought the eigenvalues of any $C_{\bf n}(f)$ are explicitly known, results like Theorem \ref{szego} and Theorem \ref{szego-herm}
do not hold for sequences $\left\{C_{\bf n}(f)\right\}_{{\bf n}\in\mathbb N^k}$ in full generality: this is due to the fact that the Fourier sum of $f$ converges to $f$ under quite restrictive assumptions (see \cite{zygmund}). 
In fact if $f$ is continuous $2\pi$-periodic and its modulus of continuity evaluated at $\delta$ goes to zero faster than $1/|\log(\delta)|$, that is
\[
\lim_{\delta\to 0^+} \log(\delta)\omega_f(\delta)=0,
\]
then $\left\{C_{\bf n}(f)\right\}_{{\bf n}\in\mathbb N^k}\sim_\lambda(f,\mathcal{I}_k)$, simply because 
$S_{\bf n}(f)(\cdot)$ uniformly converges to $f$ (see \cite{estatico-serra} for more relations between circulant sequences and 
spectral distribution results)

\end{Remark}
We end this subsection by recalling the computational properties of (block) circulants. Every matrix/vector operation with circulants has cost $O(\hat{ n}\log \hat {n})$ with moderate multiplicative constants: in particular, this is true for the matrix-vector product, for the solution of a linear system, for the computation of the blocks
 $S_{\bf n}(f)\left(\theta_{\bf r}^{({\bf n})}\right)$ and consequently of
the eigenvalues (see e.g. \cite{fft}). 

\subsection{GLT sequences: operative features}

Without going into details of the GLT algebra (see the pioneering work \cite{Tilli-lt} by Tilli for describing the
spectrum of one-dimensional differential operators and the generalization contained in \cite{glt,glt-vs-fourier} for
multi-variate differential operators), here we list some properties of the GLT sequences in their block form (see \cite{glt-vs-fourier}), used when proving that a sequence of Toeplitz matrices, up to low-rank corrections, is a GLT sequence and that its symbol is not affected by the low-rank perturbation.
\begin{description}
	\item[GLT1] Each GLT sequence has a singular value symbol $f(x,\theta)$ for $(x,\theta)\in [0,1]^k\times [-\pi,\pi]^k$ according to the second item in Definition \ref{def-distribution} with $\ell=2k$.
	If the sequence is Hermitian, then
	the distribution also holds in the eigenvalue sense.

	\item[GLT2] The set of GLT sequences form a $*$-algebra, i.e., it is closed under linear combinations, products,
	inversion (whenever the symbol vanishes, at most, in a set of zero Lebesgue measure), conjugation. Hence, the sequence
	obtained via algebraic operations on a finite set of given GLT sequences is still a GLT sequence and its symbol is
	obtained by performing the same algebraic manipulations on the corresponding symbols of the input GLT sequences.

	\item[GLT3] Every Toeplitz sequence generated by an $L^1(k,s)$ function $f=f(\theta)$ is a GLT sequences and its symbol
	is $f$, with the specifications reported in item {\bf GLT1}. We note that the function $f$ does not depend on the
	spacial variables $x\in [0,1]^k$.

	\item[GLT4] Every sequence which is distributed as the constant zero in the singular value sense is a GLT sequence with
	symbol~$0$.
\end{description}

\subsection{Analysis of the spectral symbol}
\label{sec:spectralsymbol}
Using Definition \ref{def-multilevel}, we can now explicitly express the symbol of the matrix $K_N$ in \eqref{system}. Let ${\bf n}=(n_1,n_2)$ be a $2$-index and let $\hat{n}=n_1n_2$. If $p$ is the degree of the basis functions used for the staggered DG, we obtain the following Hermitian matrix
\begin{equation}\label{K_N}
K_N=T_{\bf n}(f)+E_{\bf n}, \quad N=(p+1)^2\hat{n},
\end{equation}
where
\begin{equation*}
T_{\bf{n}}(f) =\left[\hat{f}_{\bx{i}-\bx{j}}\right]_{\bx{i},\bx{j}=\bx{e}}^{\bf n}\in\mathcal{M}_N
\end{equation*}
and $f:\mathcal{I}_2\rightarrow \mathcal{M}_{s}$, $s=(p+1)^2$, while $E_{\bf n}$ is a low-rank perturbation whose rank grows at most proportionally to $\sqrt{{\hat n}}$ and with constant depending on the bandwidths of $K_N$. The nonzero coefficients of $T_{\bf n}(f)=[\hat f_{\bx{i}-\bx{j}}]_{\bx{i},\bx{j}=\bx{e}}^{\bx{n}}$ correspond to the indexes $\bx{i}=(i_1,i_2),\bx{j}=(j_1,j_2)$ such that
\[|i_1-j_1|+|i_2-j_2|\le1.\]
Therefore, in the two-dimensional case ($k=2$) the symbol $f$ is given by
\begin{equation}
f(\theta_1 , \theta_2 ) =  \hat f_{(0,0)}+  \hat f_{(-1,0)} e^{-\mathbf i \theta_1}+  \hat f_{(0,-1)} e^{-\mathbf i \theta_2}+  \hat f_{(1,0)} e^{\mathbf i \theta_1}+  \hat f_{(0,1)} e^{\mathbf i \theta_2},
\label{eqn.symbol.f} 
\end{equation}
where $\hat f_{(0,0)}, \hat f_{(-1,0)},\hat f_{(0,-1)},\hat f_{(1,0)},\hat f_{(0,1)}\in\mathbb{R}^{(p+1)^2\times(p+1)^2}$, that is $f$ is a linear trigonometric polynomial in the variables $\theta_1$ and $\theta_2$. For detailed expressions of these matrices in the particular case $k=2$ and $p=3$, see \ref{app.symbol.k2p2}.  
Furthermore, the coefficients of $T_{\bf n}(f)$ verify the following relations
\begin{equation*}
\hat f_{(0,0)}^T=\hat f_{(0,0)}, \qquad \hat f_{(-1,0)}^T= \hat f_{(1,0)}, \qquad \hat f_{(0,-1)}^T= \hat f_{(0,1)}.
\end{equation*}
As a consequence,
\begin{equation*}
 f^*(\theta_1 , \theta_2 )=  f(\theta_1 , \theta_2 ),
\end{equation*}
that is $f$ is a Hermitian matrix-valued function which implies that $T_{\bf n}(f)$ is a Hermitian matrix. Using Theorem \ref{szego-herm}, we can conclude that
\begin{equation}\label{distr_Tn}
\{T_{\bf n}(f)\}_{{\bf n}\in\mathbb N^2}\sim_\lambda(f,\mathcal{I}_2).
\end{equation}
From \textbf{GLT3}, we know that $\{T_{\bf n}(f)\}_{{\bf n}\in\mathbb N^2}$ is a GLT sequence with symbol $f$. Moreover, let us observe that, $\{E_{\bf n}\}_{{\bf n}\in\mathbb N^2}\sim_\sigma0$ and so, by the \textbf{GLT4}, the sequence $\{E_{\bf n}\}_{{\bf n}\in\mathbb N^2}$ is a GLT sequence with symbol identically zero. Therefore, by \textbf{GLT2} and by relation \eqref{distr_Tn}, the sequence $\{T_{\bf n}(f)+E_{\bf n}\}_{{\bf n}\in\mathbb N^2}$ is a GLT sequence with symbol $f$, and
\begin{equation}\label{eig_distr_dum}
\{K_N\}_N\sim_\lambda(f,\mathcal{I}_2).
\end{equation}
Furthermore, since each $K_N$ is symmetric and its blocks are symmetric and real, from Remark \ref{rem_simmetria} with k=2, we have
\begin{equation}\label{eig_distr_dum_piu}
\{K_N\}_N\sim_\lambda(f,\mathcal{I}^+_2).
\end{equation}
Let
\begin{equation*}
\lambda_1 (K_N) \le \lambda_2 (K_N)\le \dots \le \lambda_N (K_N).
\end{equation*}
be the eigenvalues of $K_N$.
Recalling Remark \ref{rem_distr}, from equation \eqref{eig_distr_dum_piu}, we know that for $N$ sufficiently large, $N/(p+1)^2$ eigenvalues of $K_N$, up to outliers, can be approximated by a sampling of $\lambda_1(f)$ on a uniform equispaced grid of the domain $\mathcal{I}^+_2$, and so on until the last $N/(p+1)^2$ eigenvalues which can be approximated by an equispaced sampling of $\lambda_{(p+1)^2}(f)$ in the domain. 
%
In the following section we give numerical evidence of this result.
\subsection{Numerical tests}
\label{sub:Numericaltests}
Let us fix $\bx{n}=(n_1,n_2)$, with $n_1,n_2=n$, and let $p=2$. Within these choices, the matrix-size of $K_N$ defined as in \eqref{K_N} is $N=9n^2$. This section is devoted to the comparison of the eigenvalues of $K_N$ with a sampling of the eigenvalue functions $\lambda_1(f),\ldots,\lambda_9(f)$. Actually, we do not analytically compute the eigenvalue functions, but, according
to Theorem \ref{circ-s} and Remark \ref{fourier-vs-funzione}, we are able to provide an 'exact' evaluation of them on an equispaced grid on $\mathcal{I}^+_2$ (see Subsection \ref{sub:approx}) and this is sufficient for our aims.

\subsubsection{Evaluation of the eigenvalue functions of the symbol}\label{sub:approx}
Let us define the following equispaced grid on $\mathcal{I}^+_2$
\begin{equation*}
G_n=\left\{(\theta_1^{(j)},\theta_2^{(k)})=\left(\frac{j\pi}{n},\frac{k\pi}{n}\right), \qquad j,k=0, \dots, n-1\right\}
\end{equation*}
and let us consider the following $n^2$ Hermitian matrices of size $9\times9$ 
\begin{equation}
A_{j,k}:=f(\theta_1^{(j)},\theta_2^{(k)}), \quad j,k=0, \dots, n-1.
\end{equation}
Ordering in ascending way the eigenvalues of $A_{j,k}$
\begin{equation*}
\lambda_1 (A_{j,k}) \le \lambda_2  (A_{j,k}) \le \dots \le \lambda_9  (A_{j,k}), \quad j,k=0,\ldots,n-1,
\end{equation*}
for a fixed $l=1,\ldots,9$, an evaluation of $\lambda_l(f)$ at $(\theta_1^{(j)},\theta_2^{(k)})$ is given by $\lambda_l  (A_{j,k})$, $j,k=0,\ldots,n-1$. From now onwards, fixed $l$, we will denote by $P^{(n)}_l$ the vector of all eigenvalues $\lambda_l  (A_{j,k})$, $j,k=0,\ldots,n-1$, that is
\begin{equation*}
P^{(n)}_l:=\left[\lambda_l (A_{0,0}), \lambda_l  (A_{0,1}), \dots ,\lambda_l  (A_{n-1,n-1})\right],
\end{equation*}
and by $P^{(n)}$ the vector of all eigenvalues $\lambda_l  (A_{j,k})$, $j,k=0,\ldots,n-1$ varying $l$
\begin{equation*}
P^{(n)}:=\left[\lambda_1 (A_{0,0}),\ldots, \lambda_1  (A_{n-1,n-1}), \dots ,\lambda_9 (A_{0,0}),\ldots, \lambda_9  (A_{n-1,n-1})\right].
\end{equation*}


\begin{figure}[htb]
\begin{subfigure}[c]{.47\textwidth}
\includegraphics[width=\textwidth]{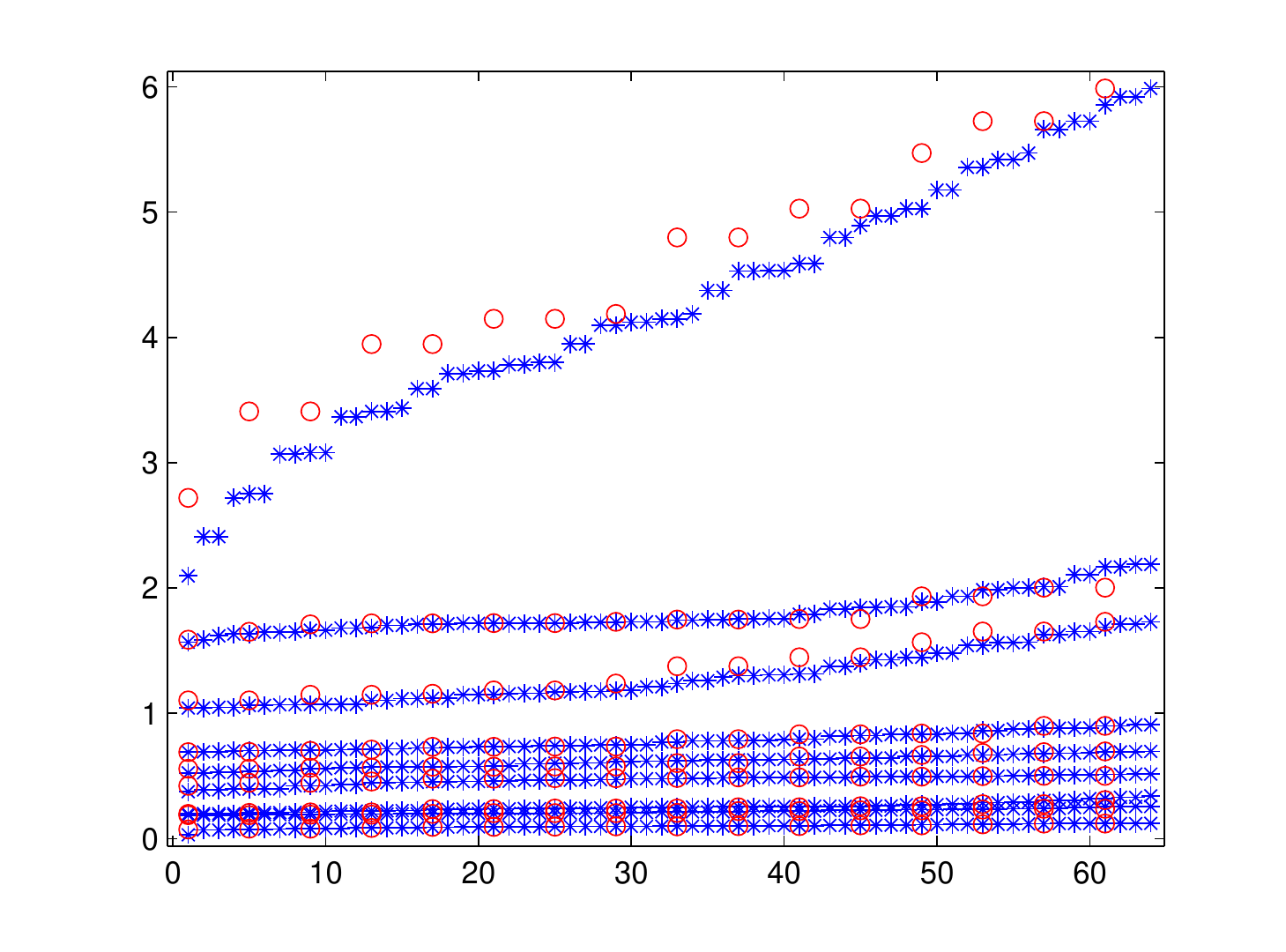}
\subcaption{$n=4$ }\label{n=4_n=8}
\end{subfigure}
\begin{subfigure}[c]{.47\textwidth}
\includegraphics[width=\textwidth]{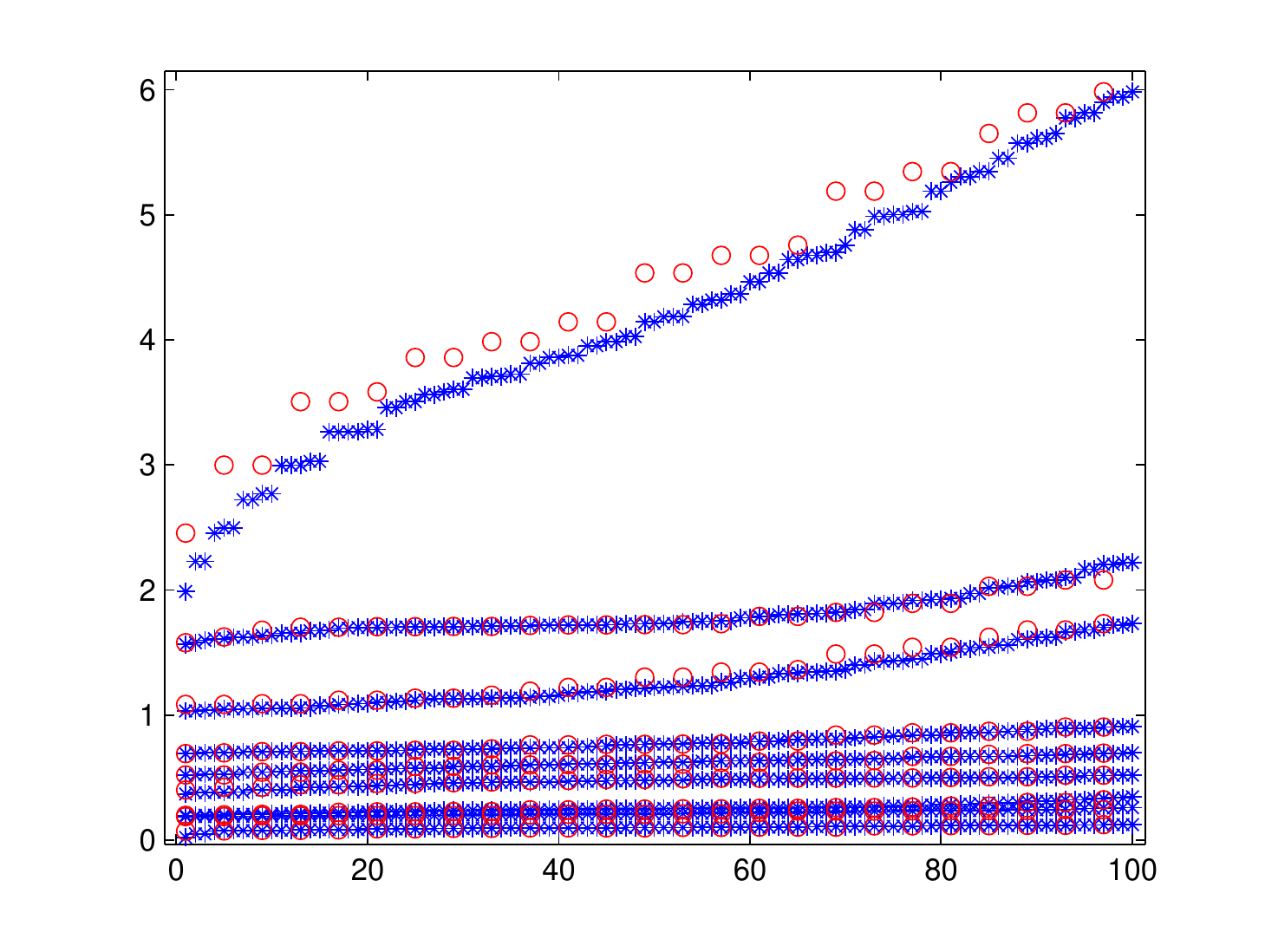}
\subcaption{$n=5$ }\label{n=5_n=10}
\end{subfigure}
\caption{Comparison between $P^{(n)}_l$(${\color{red}\circ}$) and $P^{(2n)}_l$(${\color{blue}\ast}$), $l=1,\ldots,9$ with $n=4,5$ }\label{n=4_n=8andn=5_n=10}
\end{figure}


Refining the grid $G_n$ by increasing $n$, we can provide the evaluation of the eigenvalue functions of $f$ in a larger number of
grid points: convincing numerical evidences of the latter claim are reported in Figure \ref{n=4_n=8andn=5_n=10}. More specifically, 
in Figures \ref{n=4_n=8}, \ref{n=5_n=10} we compare the approximation of $\lambda_l(f)$ on $G_n$, $n=4,5$ contained in $P^{(n)}_l$ (ordered in ascending way) with the approximation of the same eigenvalue function on a grid that is twice as fine $G_{2n}$, $n=4,5$ contained in $P^{(2n)}_l$ (ordered in ascending way as well) for every $l=1,\dots,9$.

\begin{figure}[htb]
\begin{subfigure}[c]{.30\textwidth}
\includegraphics[width=\textwidth]{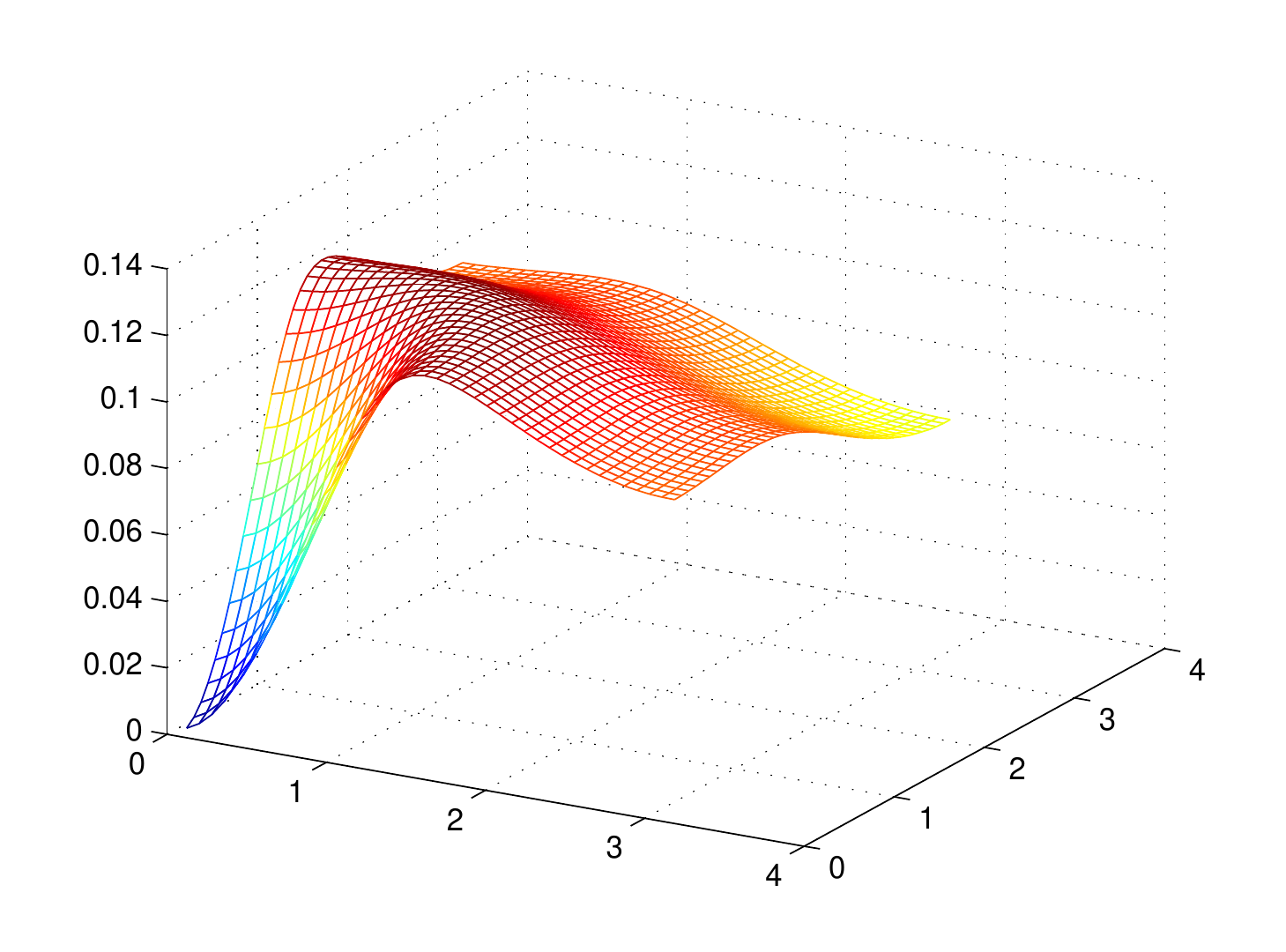}
\subcaption{$\lambda_1(f)$ }\label{valutazioni_lambda1_n=40}
\end{subfigure}
\begin{subfigure}[c]{.30\textwidth}
\includegraphics[width=\textwidth]{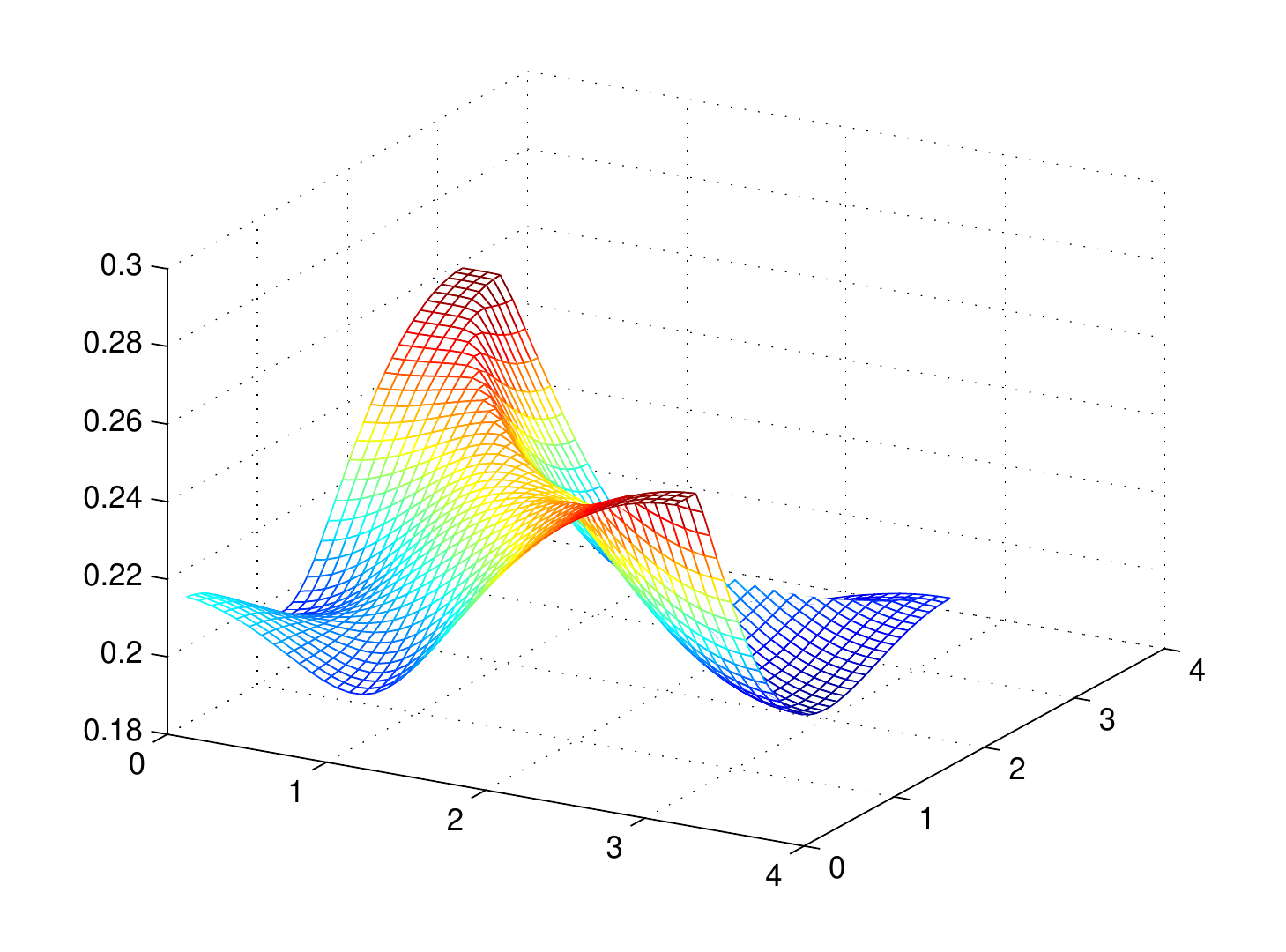}
\subcaption{$\lambda_2(f)$ }\label{valutazioni_lambda2_n=40}
\end{subfigure}
\begin{subfigure}[c]{.30\textwidth}
\includegraphics[width=\textwidth]{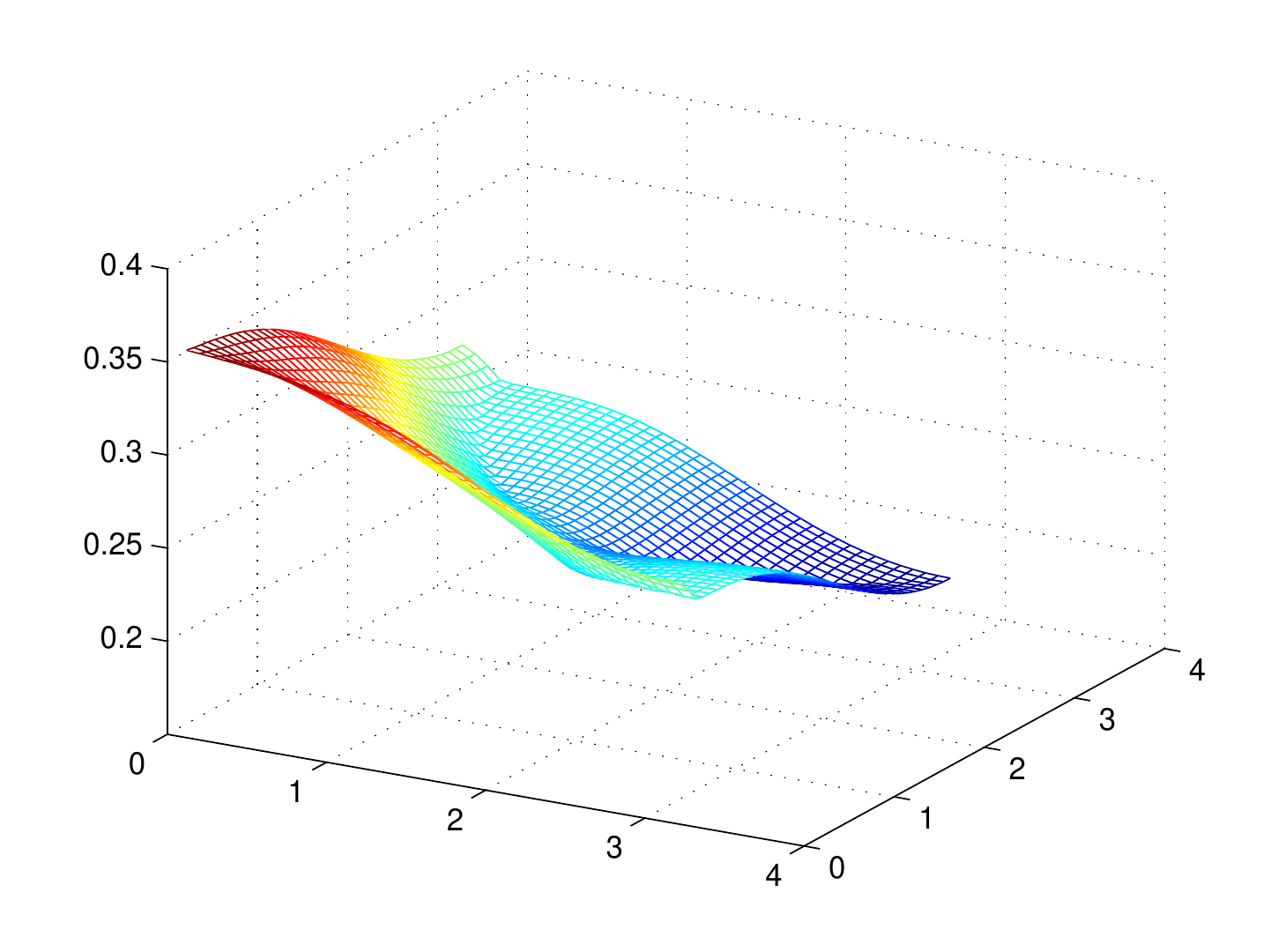}
\subcaption{$\lambda_3(f)$ }\label{valutazioni_lambda3_n=40}
\end{subfigure}
\\
\begin{subfigure}[c]{.30\textwidth}
\includegraphics[width=\textwidth]{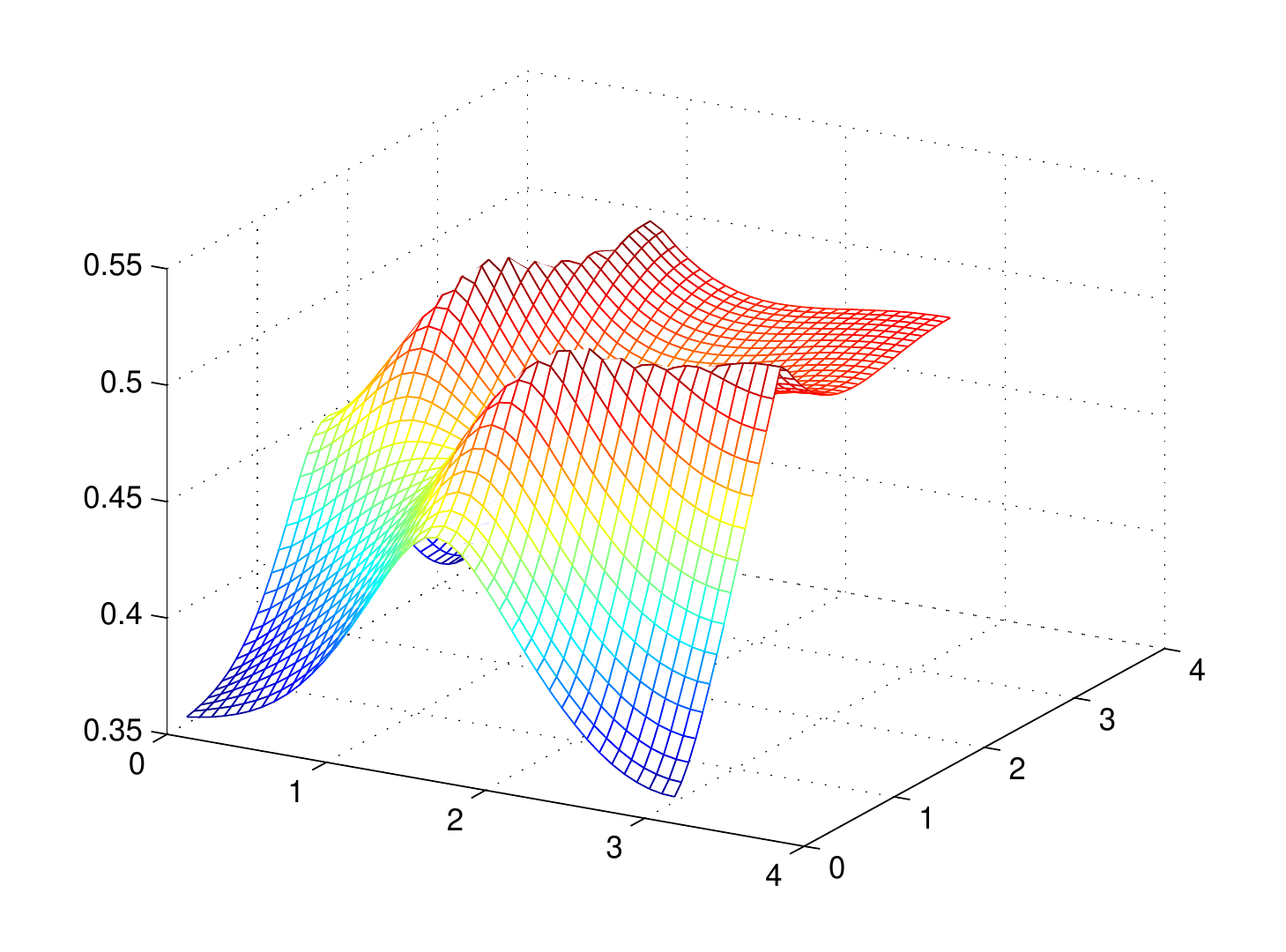}
\subcaption{$\lambda_4(f)$}\label{valutazioni_lambda4n=40}
\end{subfigure}
\begin{subfigure}[c]{.30\textwidth}
\includegraphics[width=\textwidth]{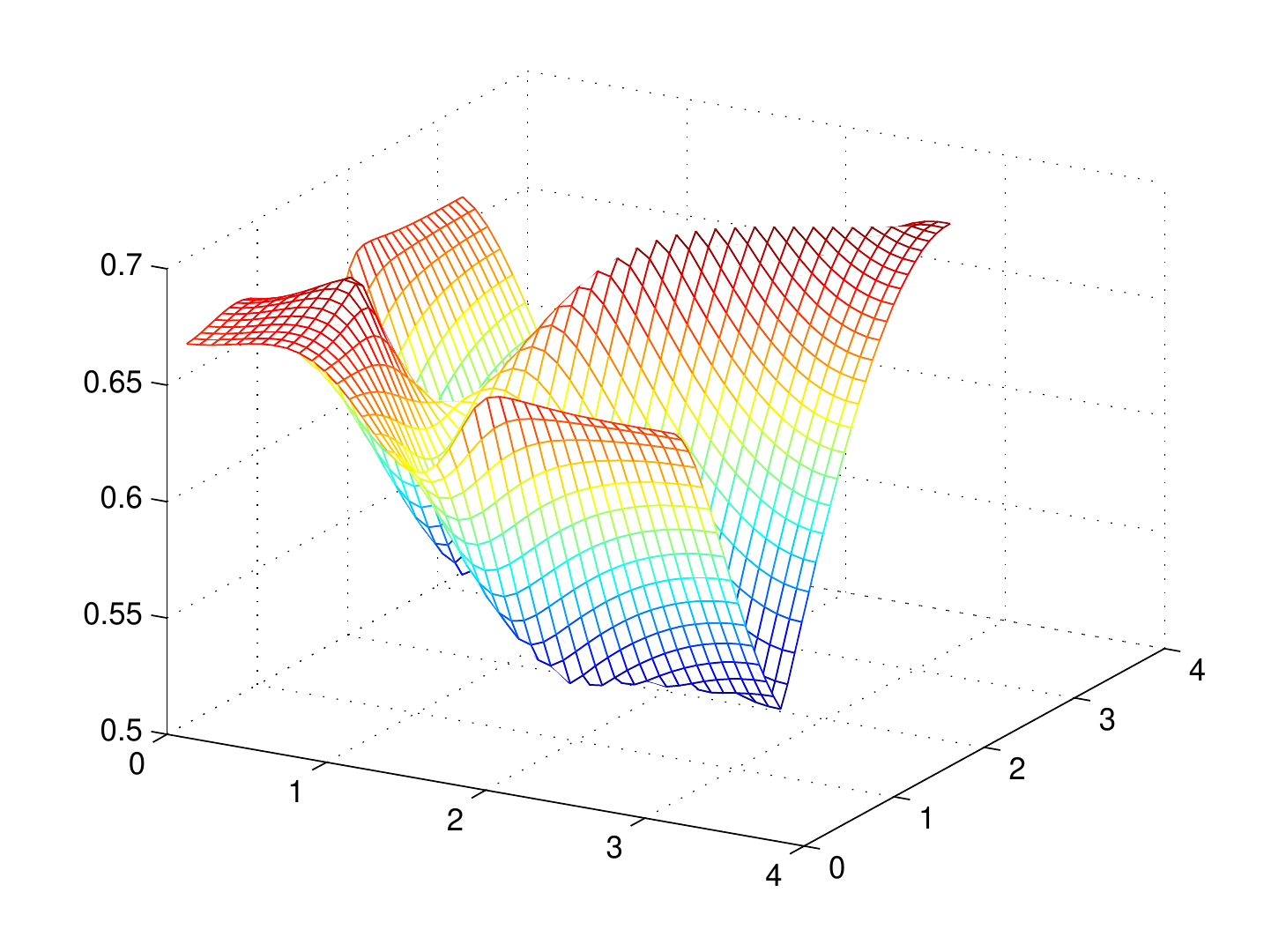}
\subcaption{$\lambda_5(f)$ }\label{valutazioni_lambda5_n=40}
\end{subfigure}
\begin{subfigure}[c]{.30\textwidth}
\includegraphics[width=\textwidth]{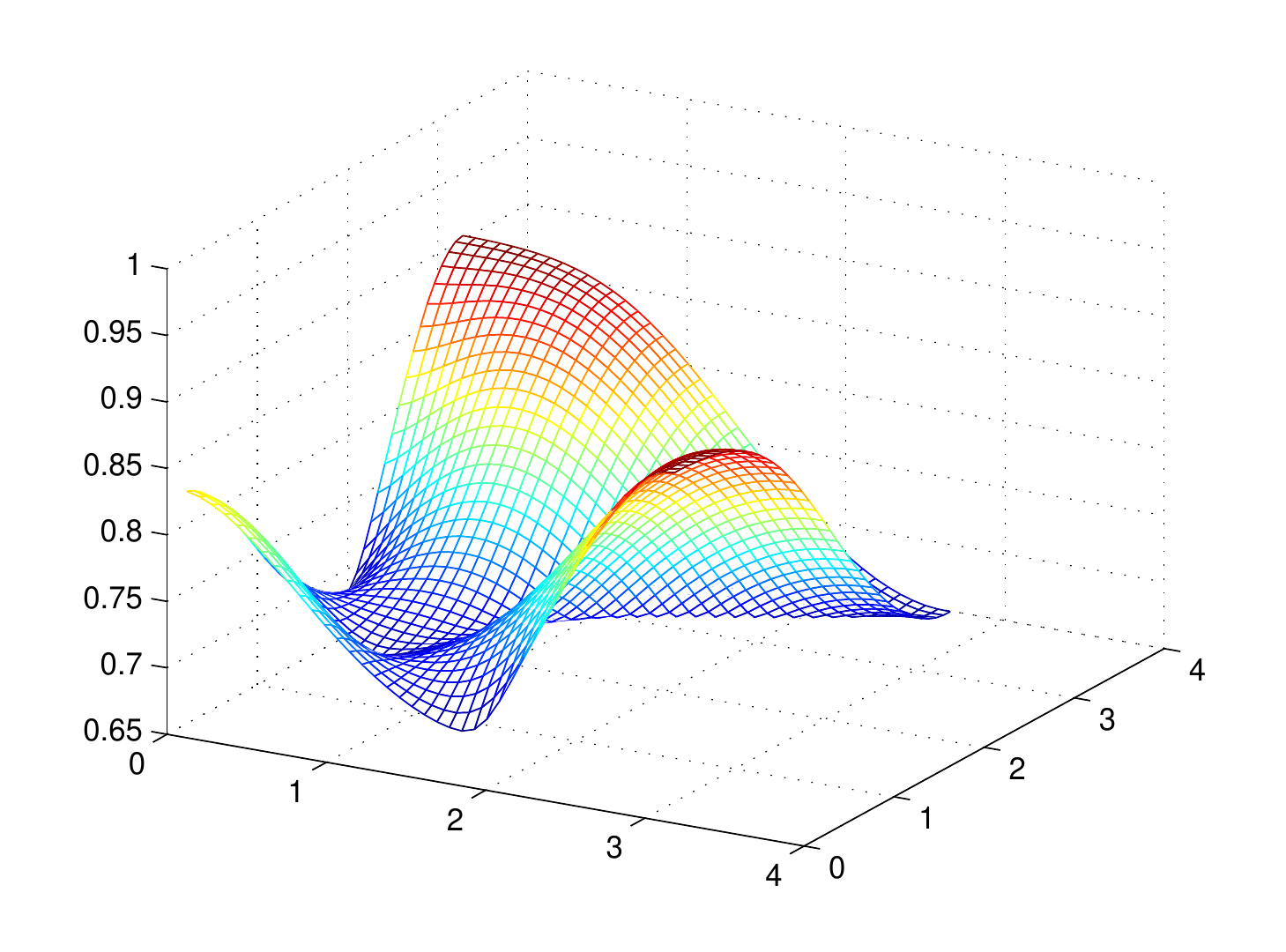}
\subcaption{$\lambda_6(f)$ }\label{valutazioni_lambda6_n=40}
\end{subfigure}
\\
\begin{subfigure}[c]{.30\textwidth}
\includegraphics[width=\textwidth]{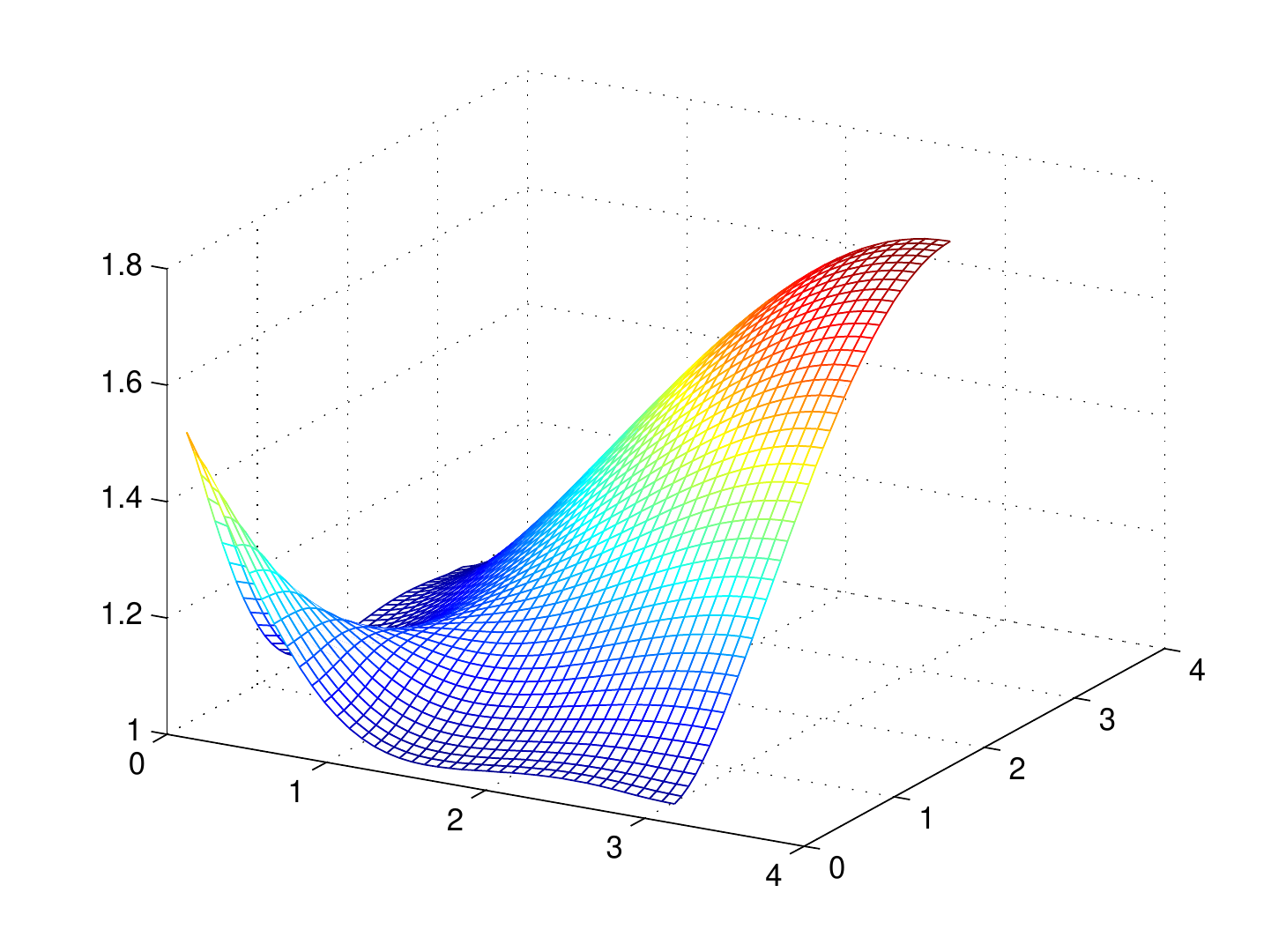}
\subcaption{$\lambda_7(f)$ }\label{valutazioni_lambda7_n=40}
\end{subfigure}
\begin{subfigure}[c]{.30\textwidth}
\includegraphics[width=\textwidth]{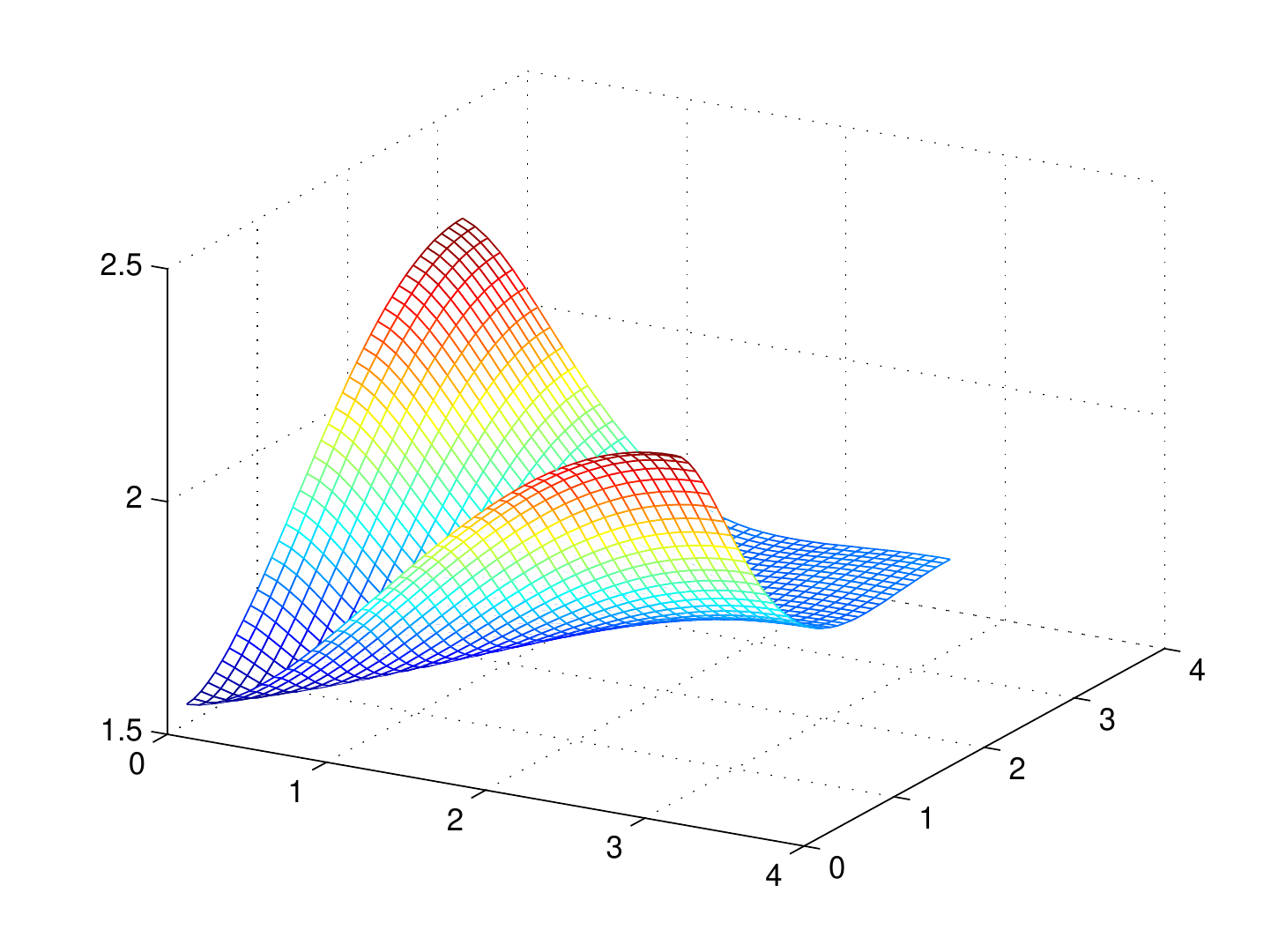}
\subcaption{$\lambda_8(f)$ }\label{valutazioni_lambda8_n=40}
\end{subfigure}
\begin{subfigure}[c]{.30\textwidth}
\includegraphics[width=\textwidth]{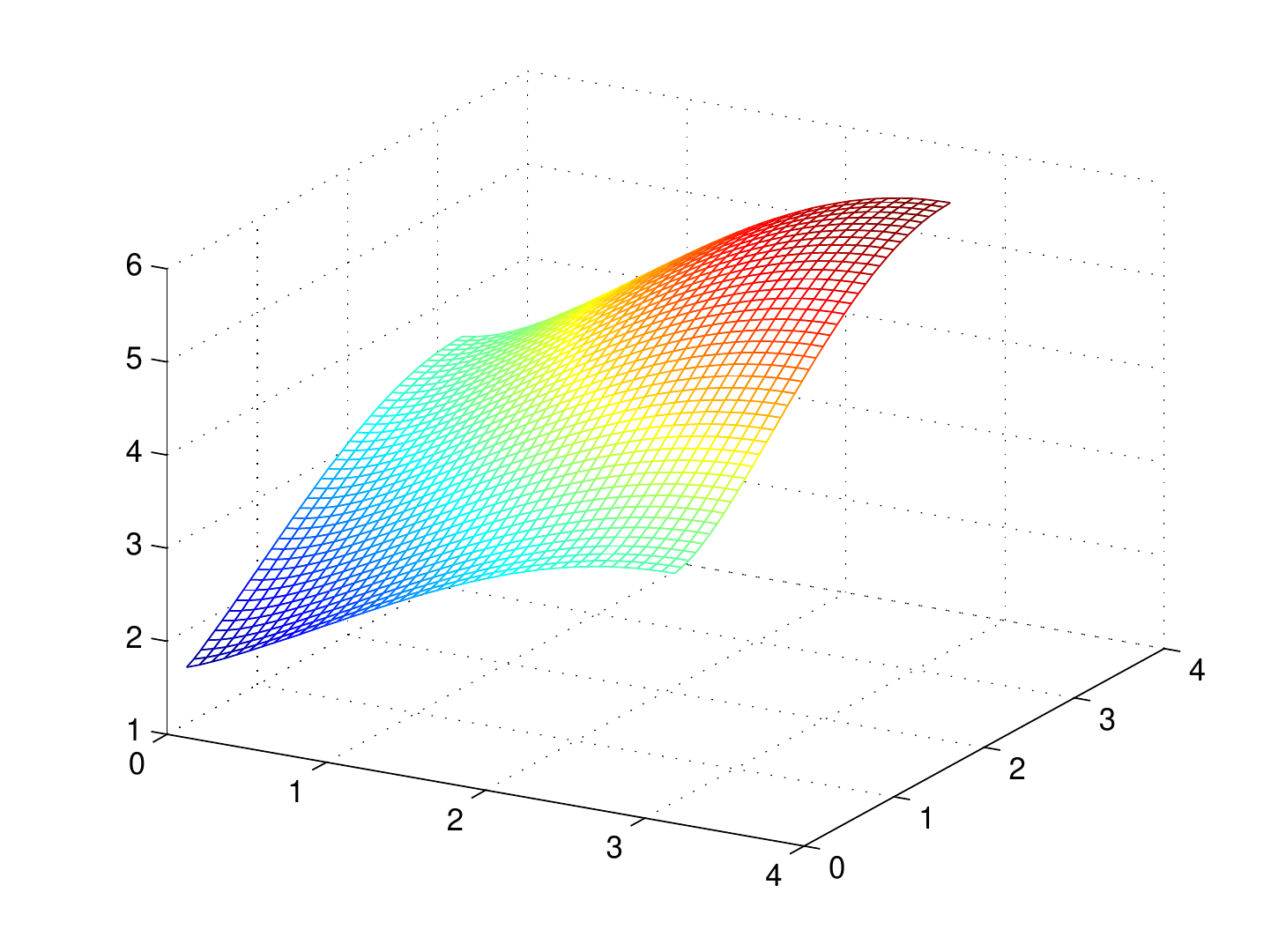}
\subcaption{$\lambda_9(f)$ }\label{valutazioni_lambda9_n=40}
\end{subfigure}
\caption{Approximation of the eigenvalues functions $\lambda_l(f)$, $l=1,\ldots,9$ as a mesh on $G_n$, when $n=40$}\label{eig_2d}
\end{figure}

Therefore, for $n$ sufficiently large, a feasible approximation of $\lambda_l(f)$, $l=1,\ldots,9$, can be obtained by displaying $P^{(n)}_l$ as a mesh on $G_n$ (see Figure \ref{eig_2d}, for $n=40$).

\subsubsection{Spectral distribution of $\{K_N\}_N$}
\label{sub:spectral_distribution_KN}

In this subsection we provide numerical evidences of the distribution result \eqref{eig_distr_dum_piu}, making use of the strategy for computing an approximation of $\lambda_l(f)$ on an equispaced grid showed in Subsection \ref{sub:approx}.

\begin{figure}[htb]
\centering
\includegraphics[scale=0.5,keepaspectratio]{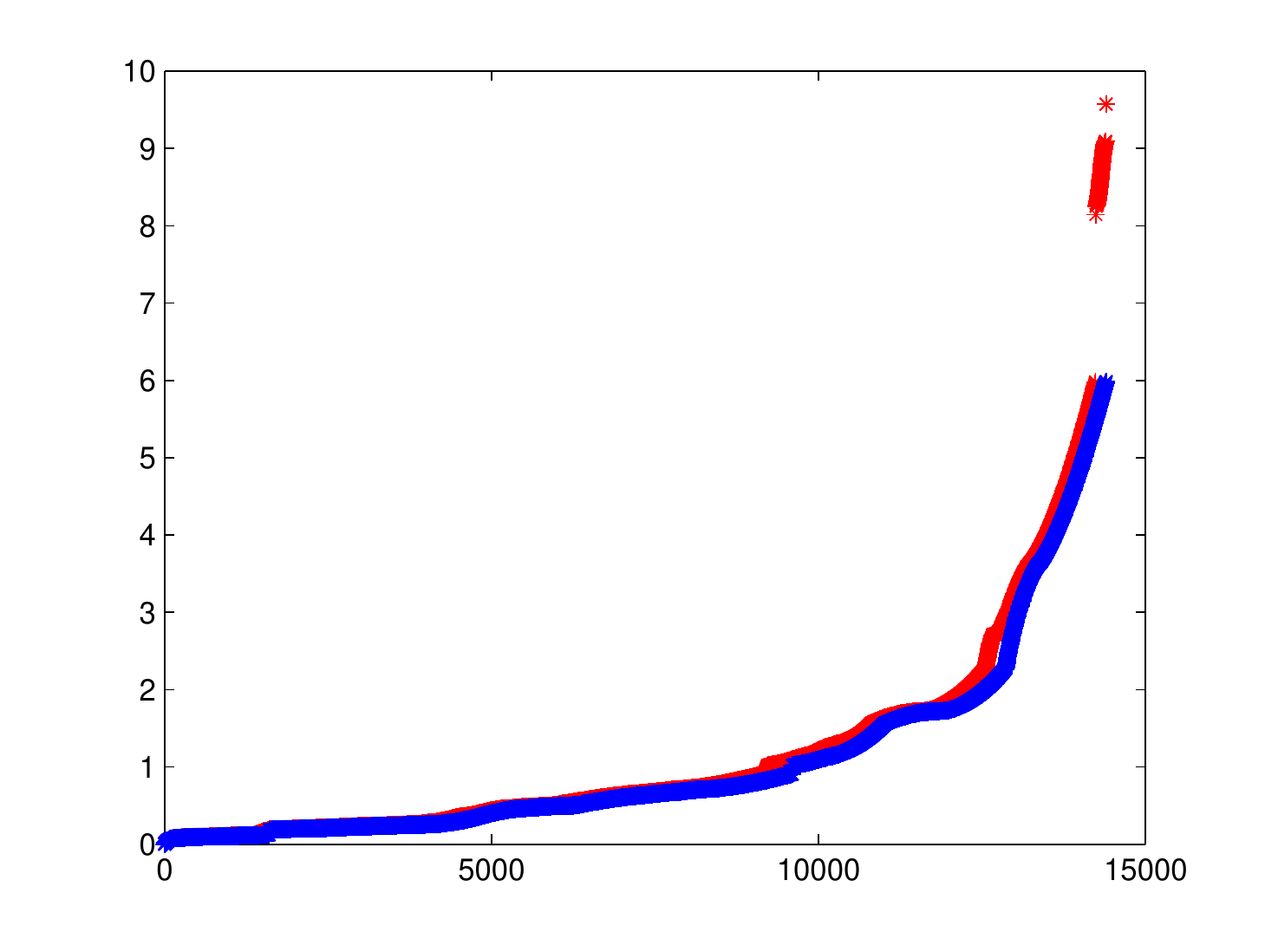}
\caption{Comparison of the eigenvalues of $K_{N}$ (\textcolor{red}{$\ast$}) with the approximation of $\lambda_l(f)$ $l=1,\ldots,9$ on $G_n$ given by a proper ordering of $P^{(n)}$ (\textcolor{blue}{$\ast$}), for $n=40$.}\label{confronto_matrice_simbolo_dum}
\end{figure}

As a first evidence, we compare the eigenvalues of $K_N$ with the evaluation of $\lambda_l(f)$ $l=1,\ldots,9$ at $G_n$ given by a proper ordering of $P^{(n)}$. As shown in Figure \ref{confronto_matrice_simbolo_dum} in which we fixed $n=40$, the eigenvalues of $K_N$ mimic, up to outliers, the sampling of the eigenvalue functions. This agrees with relation \eqref{eig_distr_dum_piu}.

Aside from such a global comparison, if
\begin{equation*}
{\rm esssup}_{\mathcal{I}^+_2}\lambda_l (f(\theta) )\le{\rm essinf}_{\mathcal{I}^+_2} \lambda_{l+1} (f(\theta)),
\end{equation*}
for some $l=1,\ldots,8$, exploiting Remark \ref{rem_distr}, we can provide a more accurate analysis of the spectrum of $K_N$ determining how many blocks it is made up of and how many eigenvalues contains each block.
With this aim, let us observe that, for a sufficiently large $n$, if we order in ascending way $P^{(n)}_l$, its extremes satisfy the following relation
\begin{equation*}
(P^{(n)}_l)_1\approx m_l,\qquad (P^{(n)}_l)_{n^2}\approx M_l, \quad l=1,\ldots,9.
\end{equation*}
A satisfactory approximation of $[m_l,M_l]$ can be numerically computed by setting $n=500$; as a result we obtain the following approximations

\begin{align*}
[m_1,M_1]&\approx[ 0.000000000,0.123775621],\\
[m_2,M_2]&\approx[0.186715287, 0.260786617],\\
[m_3,M_3]&\approx[0.197732806,0.355965321],\\
[m_4,M_4]&\approx[0.355965321, 0.524158720],\\
[m_5,M_5]&\approx[0.520903995,0.696882517],\\
[m_6,M_6]&\approx[0.677870643, 0.910001758],\\
[m_7,M_7]&\approx[1.015599697,1.731431133],\\
[m_8,M_8]&\approx[1.560701345, 2.284336270],\\
[m_9,M_9]&\approx[1.651355307, 5.985129348].\\
\end{align*}

However, looking at 
\begin{equation*}
f(0,0)= \begin{bmatrix}
  
   \frac  { 19}{45  }    &   \frac  { 1}{60    }& \frac{      -7}{40   }& \frac{        1}{60     }& \frac{     -4}{45    }& \frac{      -7}{180  }& \frac{       -7}{40     }& \frac{     -7}{180        }& \frac{ 11}{180} \\  
   &         &    &  &  &    &  & &\\ 
     \frac {  1}{60      }& \frac{    46}{45   }& \frac{        1}{60  }& \frac{        -4}{45  }& \frac{        -4}{15  }& \frac{        -4}{45  }& \frac{        -7}{180    }& \frac{     -8}{15  }& \frac{        -7}{180 } \\ 
     &         &    &  &  &    &  & &\\ 
   \frac  {  -7}{40    }& \frac{       1}{60      }& \frac{    19}{45   }& \frac{       -7}{180  }& \frac{       -4}{45   }& \frac{        1}{60   }& \frac{       11}{180   }& \frac{      -7}{180   }& \frac{      -7}{40  } \\ 
   &         &    &  &  &    &  & &\\ 
    \frac  {  1}{60     }& \frac{     -4}{45     }& \frac{     -7}{180   }& \frac{      46}{45     }& \frac{     -4}{15       }& \frac{   -8}{15  }& \frac{         1}{60     }& \frac{     -4}{45 }& \frac{         -7}{180}\\  
    &         &    &  &  &    &  & &\\  
   \frac   { -4}{45   }& \frac{       -4}{15  }& \frac{        -4}{45   }& \frac{       -4}{15    }& \frac{      64}{45    }& \frac{      -4}{15      }& \frac{    -4}{45   }& \frac{       -4}{15 }& \frac{         -4}{45 } \\  
   &         &    &  &  &    &  & &\\ 
   \frac  {  -7}{180   }& \frac{      -4}{45     }& \frac{      1}{60    }& \frac{      -8}{15   }& \frac{       -4}{15    }& \frac{      46}{45     }& \frac{     -7}{180   }& \frac{      -4}{45 }& \frac{          1}{60 } \\  
   &         &    &  &  &    &  & &\\ 
    \frac {  -7}{40    }& \frac{      -7}{180    }& \frac{     11}{180    }& \frac{      1}{60   }& \frac{       -4}{45    }& \frac{      -7}{180      }& \frac{   19}{45     }& \frac{      1}{60}& \frac{          -7}{40   } \\
    &         &    &  &  &    &  & &\\ 
    \frac  { -7}{180   }& \frac{      -8}{15     }& \frac{     -7}{180    }& \frac{     -4}{45   }& \frac{       -4}{15      }& \frac{    -4}{45   }& \frac{        1}{60   }& \frac{       46}{45 }& \frac{          1}{60  } \\ 
    &         &    &  &  &    &  & &\\ 
 \frac {     11}{180      }& \frac{   -7}{180    }& \frac{     -7}{40     }& \frac{     -7}{180   }& \frac{      -4}{45   }& \frac{        1}{60     }& \frac{     -7}{40   }& \frac{        1}{60 }& \frac{         19}{45 }\\
\end{bmatrix}
\end{equation*}

we observe that the matrix has row sum equal to zero for every row.\\
 This means that $f(0,0)e=0$ where $e\in \mathbb{R}^9$ is the vector of all ones. Therefore $f(0,0)$ is analytically singular and $m_1=0$, since the symbol is theoretically nonnegative definite due to the Galerkin approach.
Now, recalling the second item of Theorem \ref{loc-extr-s} and observing that $f(\pi,\pi)$ is positive definite, we deduce that $\lambda_1(f(\theta_1,\theta_2))$ has positive maximum and therefore the interval $[m_1,M_1]$ can be replaced by $(0,M_1]$.

From now onwards, we assume $(0,M_1] $, $\left( m_l,M_l\right)$, $l=2,\ldots,9$, to be equal to its estimate. Let us observe that the following relations hold

\begin{equation}\label{M_3}
\begin{split}
M_1&<m_2,\\
 M_3&=m_4,\\
M_6&<m_7.
\end{split}
\end{equation}

In other words, according to relations \eqref{eig_distr_dum_piu}, \eqref{M_3}, and Remark \ref{rem_distr}, 
we expect the eigenvalues of $K_N$ to verify
\begin{equation} \label{rel_aut}
\begin{split}
\# \left\{ i \, : \, \lambda_i(K_N) \in [m_1,M_1]\right\}&=\frac{9n^2}{9}+o(9n^2),\\
\# \left\{ i \, : \, \lambda_i(K_N) \in [m_2,M_3]\right\}&=2\frac{9n^2}{9}+o(9n^2),\\
\# \left\{ i \, : \, \lambda_i(K_N) \in [m_4,M_6]\right\}&=3\frac{9n^2}{9}+o(9n^2),\\
\# \left\{ i \, : \, \lambda_i(K_N) \in [m_7,M_9]\right\}&=3\frac{9n^2}{9}+o(9n^2),
\end{split}
\end{equation}
and then to identify $4$ blocks
\begin{align*}
{\rm Bl}_1&=\left[\lambda_1 (K_N), \dots,\lambda_{n^2}(K_N)\right],\\
{\rm Bl}_2&=\left[\lambda_{n^2+1} (K_N), \dots,\lambda_{3n^2}(K_N)\right],\\
{\rm Bl}_3&=\left[\lambda_{3n^2+1} (K_N), \dots,\lambda_{6n^2}(K_N)\right],\\
{\rm Bl}_4&=\left[\lambda_{6n^2+1} (K_N), \dots,\lambda_{9n^2}(K_N)\right].
\end{align*}
Correspondingly, we can split the vector $P^{(n)}$ containing the sampling of the eigenvalue functions on $G_n$ as follows
\begin{align*}
{\rm Eval}_1&= [(P^{(n)})_{1}, \ldots, (P^{(n)})_{n^2}],\\
{\rm Eval}_2&= [(P^{(n)})_{n^2+1},\ldots, (P^{(n)})_{3n^2}],\\
{\rm Eval}_3&= [(P^{(n)})_{3n^2+1}, \ldots, (P^{(n)})_{6n^2}], \\
{\rm Eval}_4&= [(P^{(n)})_{6n^2+1}, \ldots, (P^{(n)})_{9n^2}].\\
\end{align*}

Note that because of \eqref{rel_aut}, a number of outliers infinitesimal in the dimension $N$ is allowed.
For instance, when $n=40$ ($N=14400$), we find
$$\frac{9n^2}{9}=1600, \qquad 2\frac{9n^2}{9}=3200, \qquad  3\frac{9n^2}{9}=4800,$$
and
\begin{align}\label{rel_aut40}
\begin{split}
\# \left\{ i \, : \, \lambda_i(K_N) \in [m_1,M_1]\right\}=1444,\\
\# \left\{ i \, : \, \lambda_i(K_N) \in [m_2,M_3]\right\}=2911,\\
\# \left\{ i \, : \, \lambda_i(K_N) \in [m_4,M_6]\right\}=4670,\\
\# \left\{ i \, : \, \lambda_i(K_N) \in [m_7,M_9]\right\}=5016.
\end{split}
\end{align}

\begin{figure}[htb]
\centering
\includegraphics[scale=0.5,keepaspectratio]{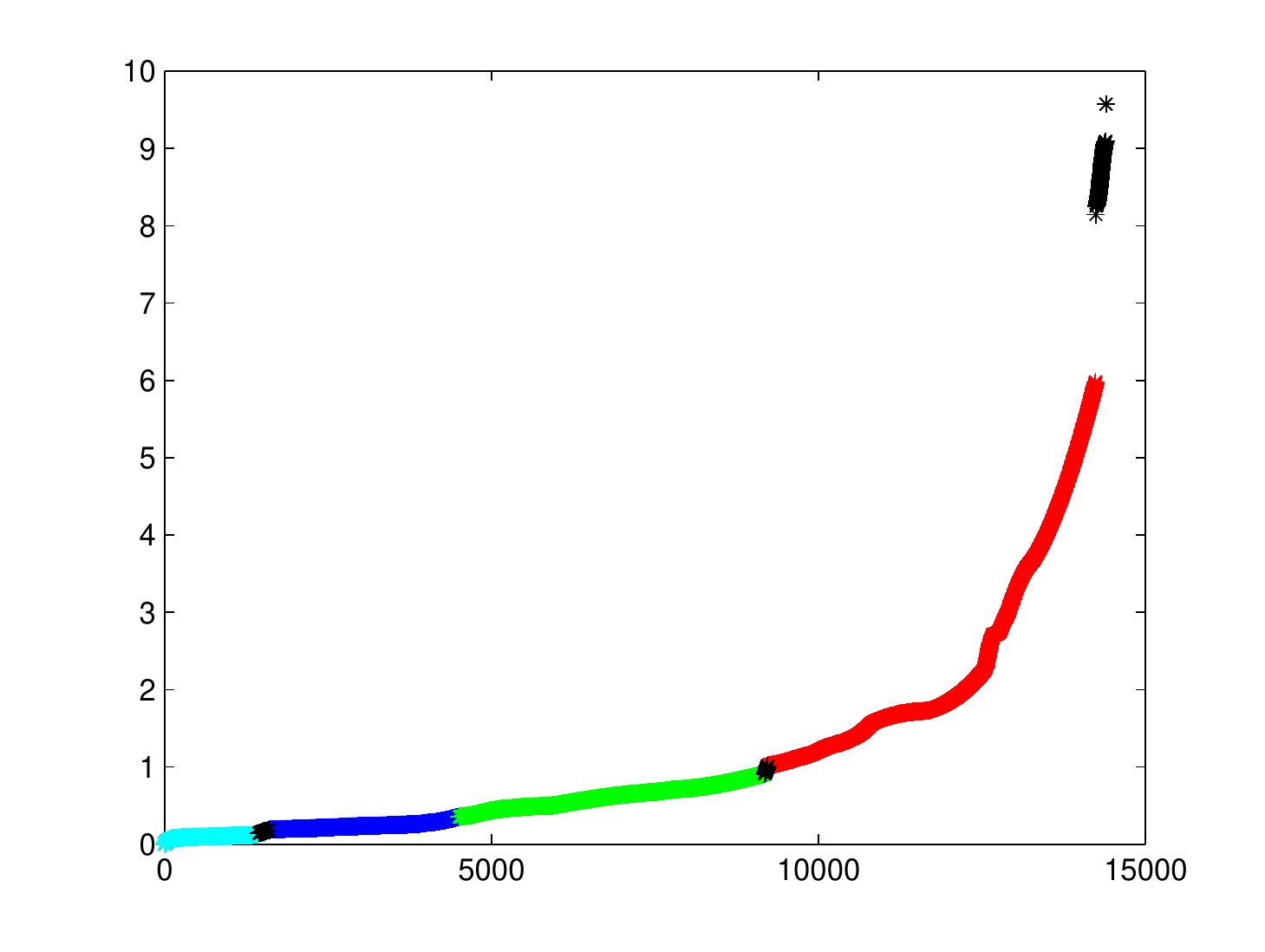}
\caption{Eigenvalues of $K_N$ for $n=40$ ({\color{black}$\ast$}) together with the eigenvalues of $K_N$ satisfying \eqref{rel_aut} ({\color{cyan}$\ast$})({\color{blue}$\ast$})({\color{green}$\ast$})({\color{red}$\ast$})}\label{autovalori_dumbser_4blocchi_n_40_consaltisullax_outliars}
\end{figure}

Therefore, from relations \eqref{rel_aut40}, we expect a number of eigenvalues of $K_N$ which are in none of the blocks or which are in the \lq wrong' block (5016 effective against 4800 expected eigenvalues in the last block). This is confirmed by Figure \ref{autovalori_dumbser_4blocchi_n_40_consaltisullax_outliars} in which we represent in black the whole spectrum of $K_N$ and highlight by means of different colours the eigenvalues belonging to different blocks. On the other hand, such a phenomenon is in line with relations \eqref{rel_aut} and the order of what is missing/exceeding is infinitesimal in the dimension $N$. As an example, in Table \ref{tab:outliers_K_N_m1_M1} we compare the actual number of eigenvalues of $K_N$ contained in the first interval $[m_1,M_1]$ with the expected number $9n^2/9$. In such way, we succeed in counting the outliers of $K_N$ in $[m_1,M_1]$, whose cardinality behaves as $O(\sqrt{9n^2})$.

\begin{table}[htb]
\begin{center}
\begin{tabular}{|c|c|c|c|c|c|c|}
\hline
$n$&eigs in $[m_1,M_1]$& $9n^2/9$& Out. &Out./$\sqrt{9n^2}$\\
\hline\hline
10 & 64&100&36&$1.20$\\
15 & 169&225&56&$1.24$\\
20&324&400&76&$1.26$\\
25&529&625&96&$1.28$\\
30&784 & 900 &  116 & $1.29$ \\
35&1089&1225&136&$  1.29$\\
40& 1444&1600&156&1.30\\
\hline
\hline
\end {tabular}
\caption{Comparison of the effective number of eigenvalues of $K_N$ contained in the first interval $[m_1,M_1]$ with the expected number $9n^2/9$}\label{tab:outliers_K_N_m1_M1}
\end{center}
\end{table}

A further evidence of relation \eqref{eig_distr_dum_piu} can be obtained by comparing block by block the eigenvalues of $K_N$ with the sampling of the eigenvalue functions of $f$, that is comparing Bl$_1$, Bl$_2$, Bl$_3$, Bl$_4$ with Eval$_1$, Eval$_2$, Eval$_3$, Eval$_4$, respectively.
Two possibilities are available. 
\begin{itemize}
\item On the one hand, we can order Eval$_t$ in ascending way and compare it with Bl$_t$.

    As an example, in Figure \ref{confronto_matrice_simbolo_primo_blocco_dum} we compare Bl$_1$ with Eval$_1$ fixed $n=40$.
\begin{figure}[htb]
\centering
\includegraphics[scale=0.5,keepaspectratio]{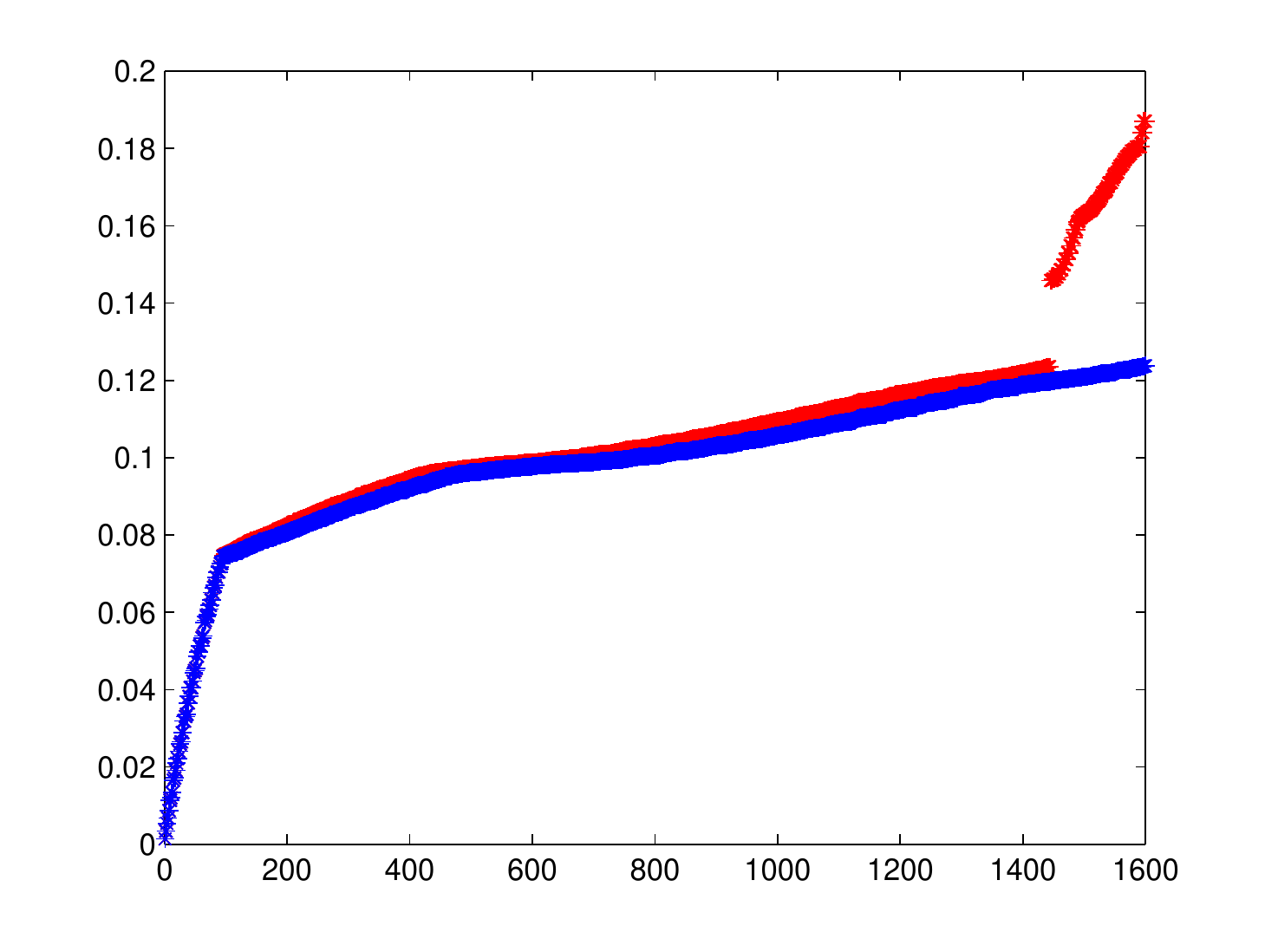}
\caption{Comparison between Bl$_1$ (\textcolor{red}{$\ast$}) and Eval$_1$
 (\textcolor{blue}{$\ast$}), for $n=40$}\label{confronto_matrice_simbolo_primo_blocco_dum}
\end{figure}
Note that a certain number of eigenvalues of $K_N$ seems not to behave as the corresponding sampling of $\lambda_1(f)$. Nevertheless, a direct computation showed that such a number agrees with the one reported in Table \ref{tab:outliers_K_N_m1_M1}. Similar results can be obtained in the comparison between Bl$_2$ with Eval$_2$, Bl$_3$ with Eval$_3$, Bl$_4$ with Eval$_4$.

\item On the other hand, we can compare the elements of Eval$_t$ 
    with the elements of Bl$_t$ by means of the following matching algorithm
    \begin{itemize}
    \item for a fixed $\lambda\in {\rm Bl}_t$ find $\tilde \eta\in{\rm Eval}_t$ such that
\begin{equation*}
\|\lambda-\tilde \eta\|=\min_{\eta\in{\rm Eval}_t}\|\lambda-\eta\|;
\end{equation*}
\item associate $\lambda$ to the couple in $G_n$ corresponding to $\tilde \eta$.
\end{itemize}
Making use of the previous algorithm, in Figure \ref{approx_2d}, we compare the eigenvalues of $K_N$ with $\lambda_l(f)$, $l=1,\ldots,9$ displayed as a mesh on $G_n$, for $n=40$. Once again, the eigenvalues of $K_N$ mimic, up to outliers, the sampling of the eigenvalue functions.

Moreover, looking at Figure \ref{approssimazioni_lambda1_n=40}, we computed the eigenvalues of $K_N$ which do not behave as the corresponding sampling of $\lambda_1(f)$ and, as expected, their order is $O(\sqrt{ 9 n^2})$ (see again Table \ref{tab:outliers_K_N_m1_M1}). As an additional confirmation of such a behaviour, in Table \ref{tab:outliers_9} we show the number of outliers of $K_N$ with respect to the sampling of $\lambda_9(f)$ (see Figure \ref{approssimazioni_lambda9_n=40}). 

\begin{figure}[htb]
\begin{subfigure}[c]{.30\textwidth}
\includegraphics[width=\textwidth]{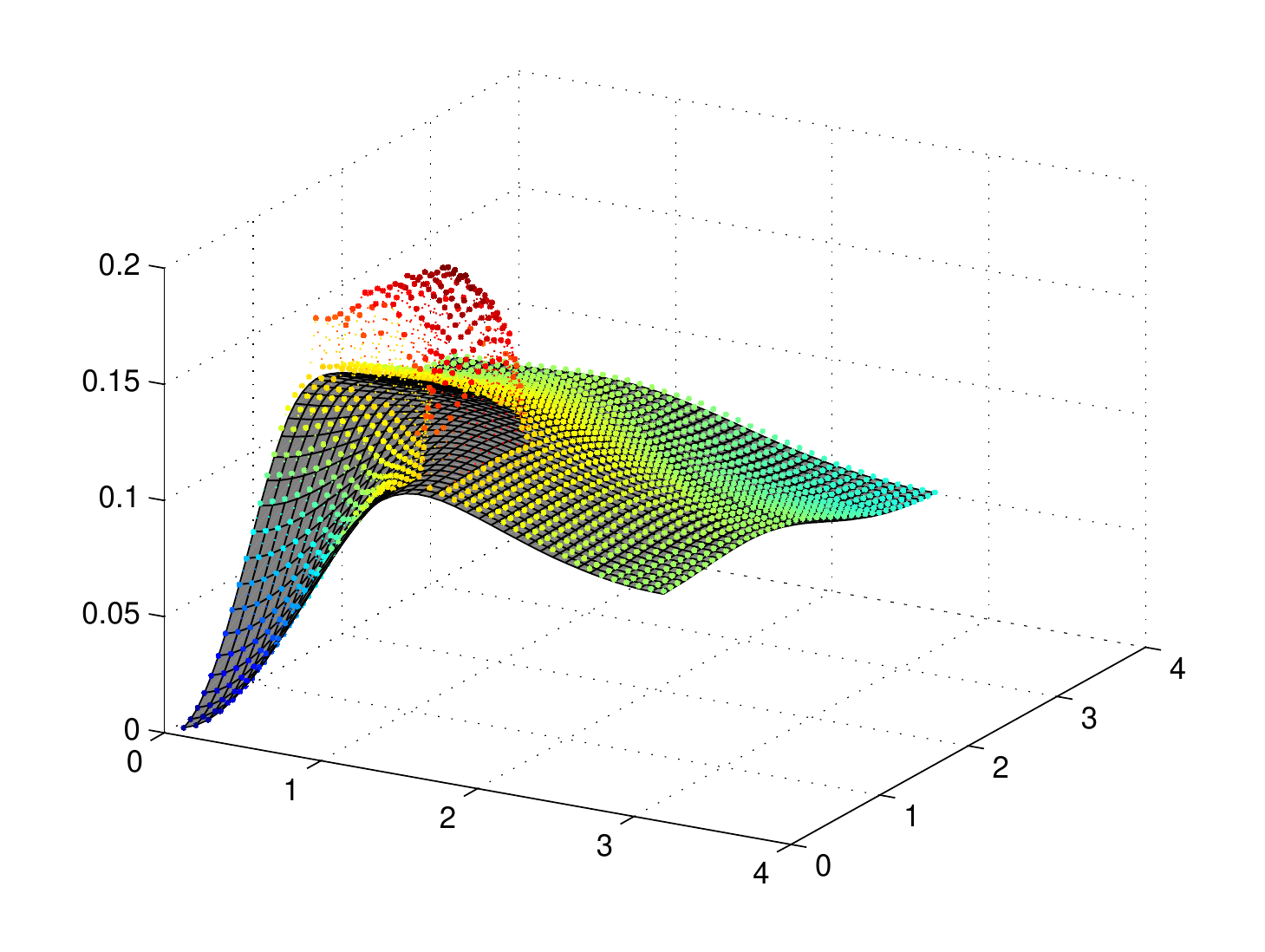}
\subcaption{$\lambda_1(f)$ }\label{approssimazioni_lambda1_n=40}
\end{subfigure}
\begin{subfigure}[c]{.30\textwidth}
\includegraphics[width=\textwidth]{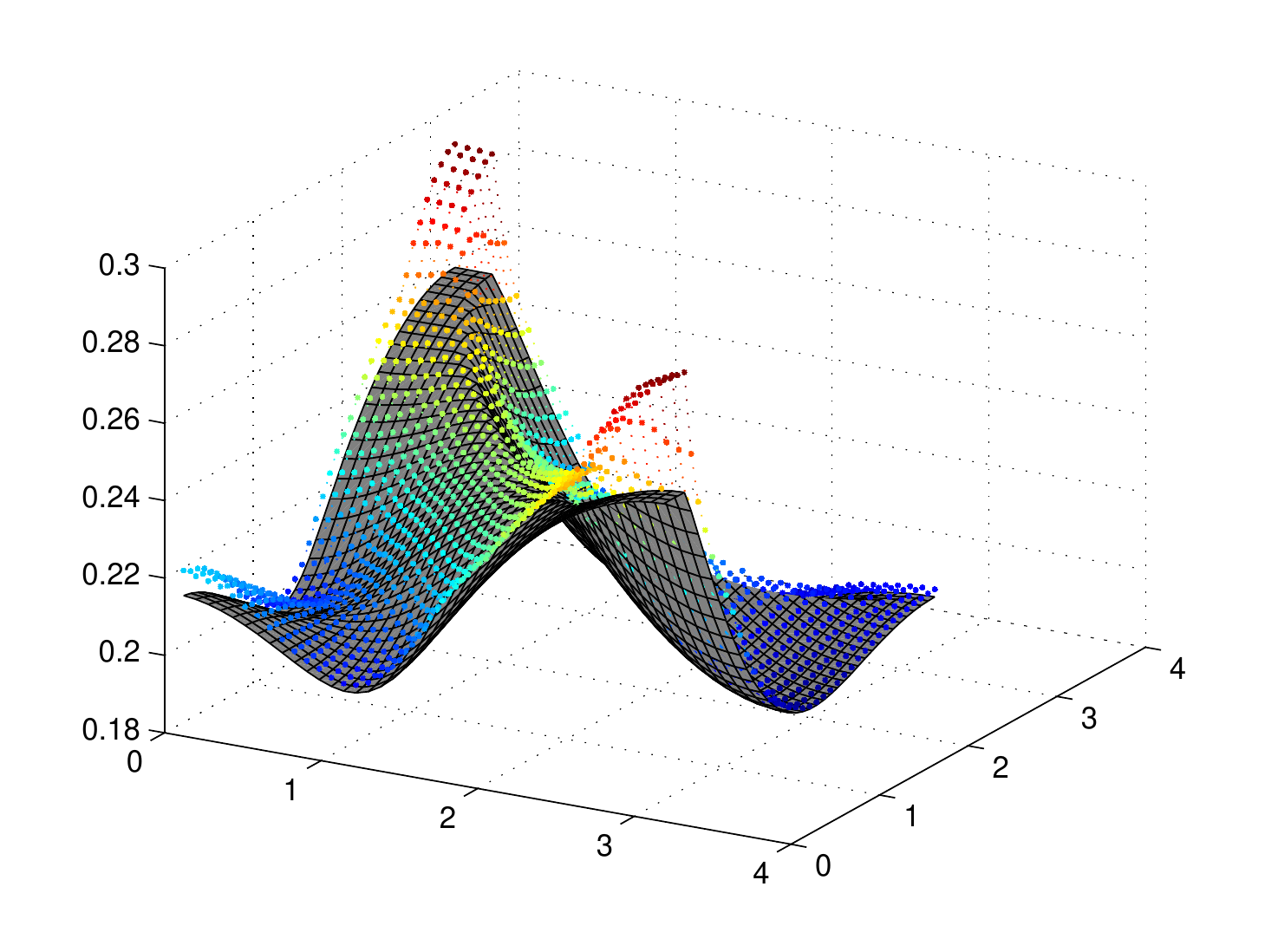}
\subcaption{$\lambda_2(f)$ }\label{approssimazioni_lambda2_n=40}
\end{subfigure}
\begin{subfigure}[c]{.30\textwidth}
\includegraphics[width=\textwidth]{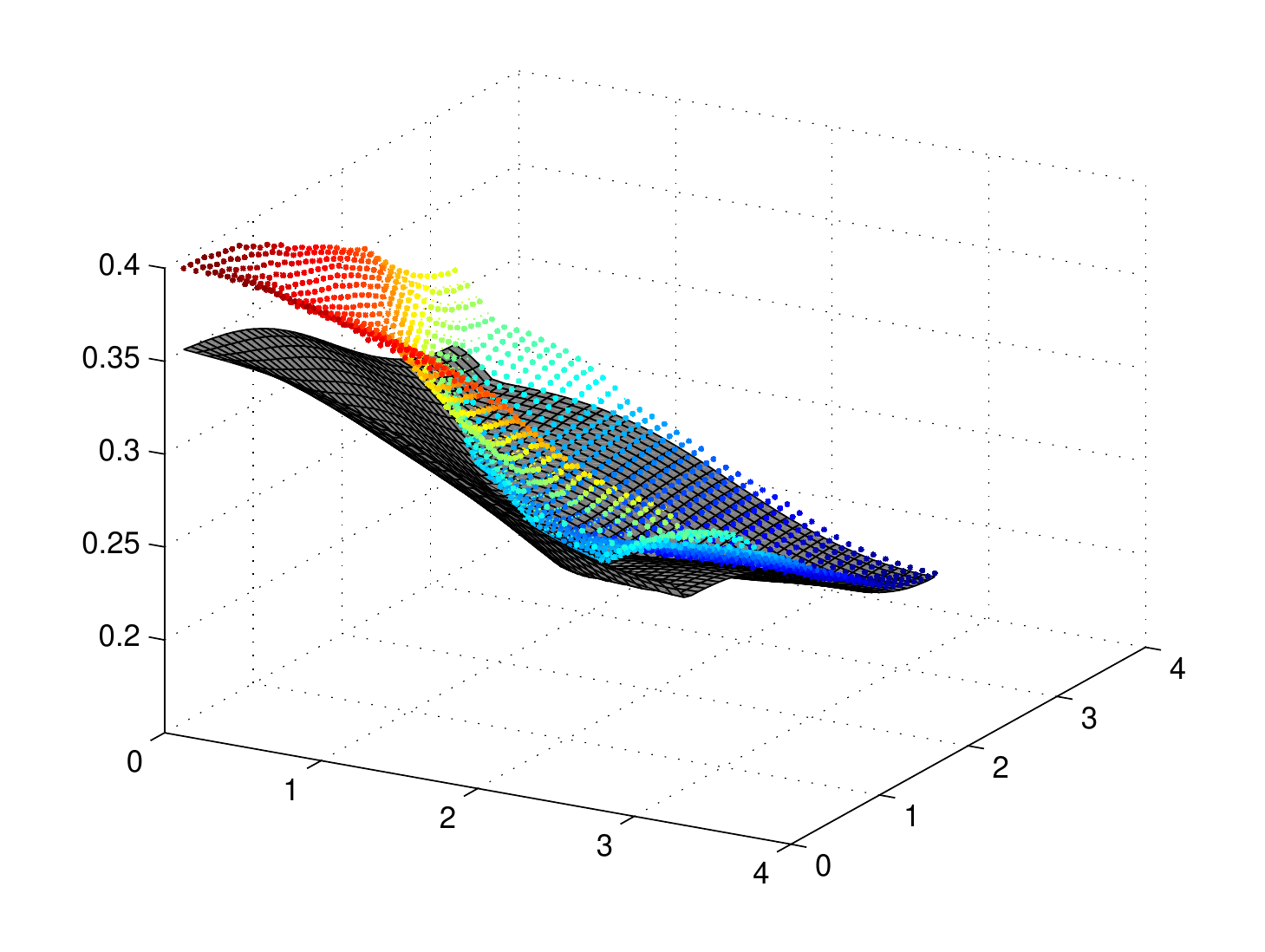}
\subcaption{$\lambda_3(f)$ }\label{approssimazioni_lambda3_n=40}
\end{subfigure}
\\
\begin{subfigure}[c]{.30\textwidth}
\includegraphics[width=\textwidth]{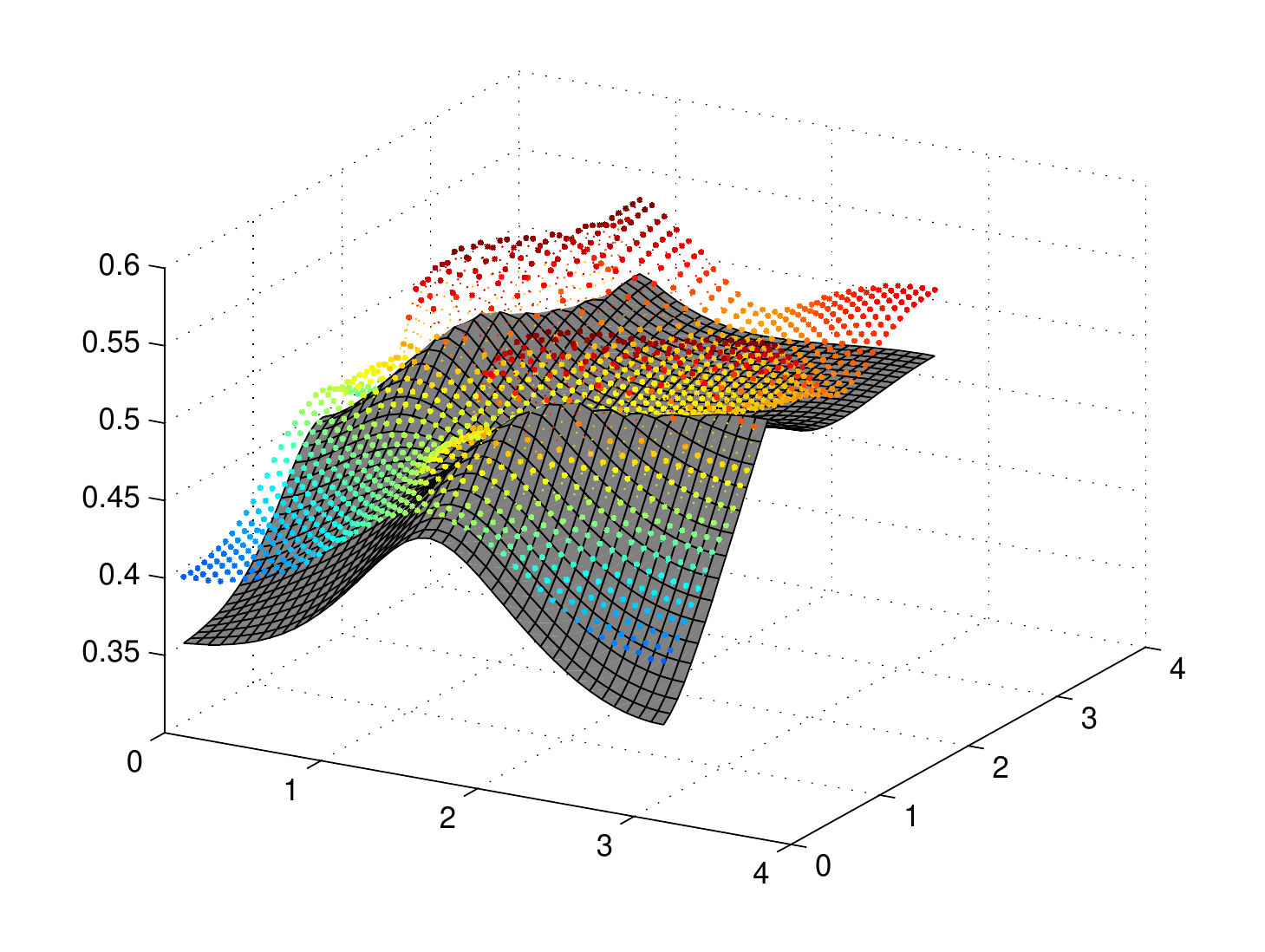}
\subcaption{$\lambda_4(f)$}\label{approssimazioni_lambda4n=40}
\end{subfigure}
\begin{subfigure}[c]{.30\textwidth}
\includegraphics[width=\textwidth]{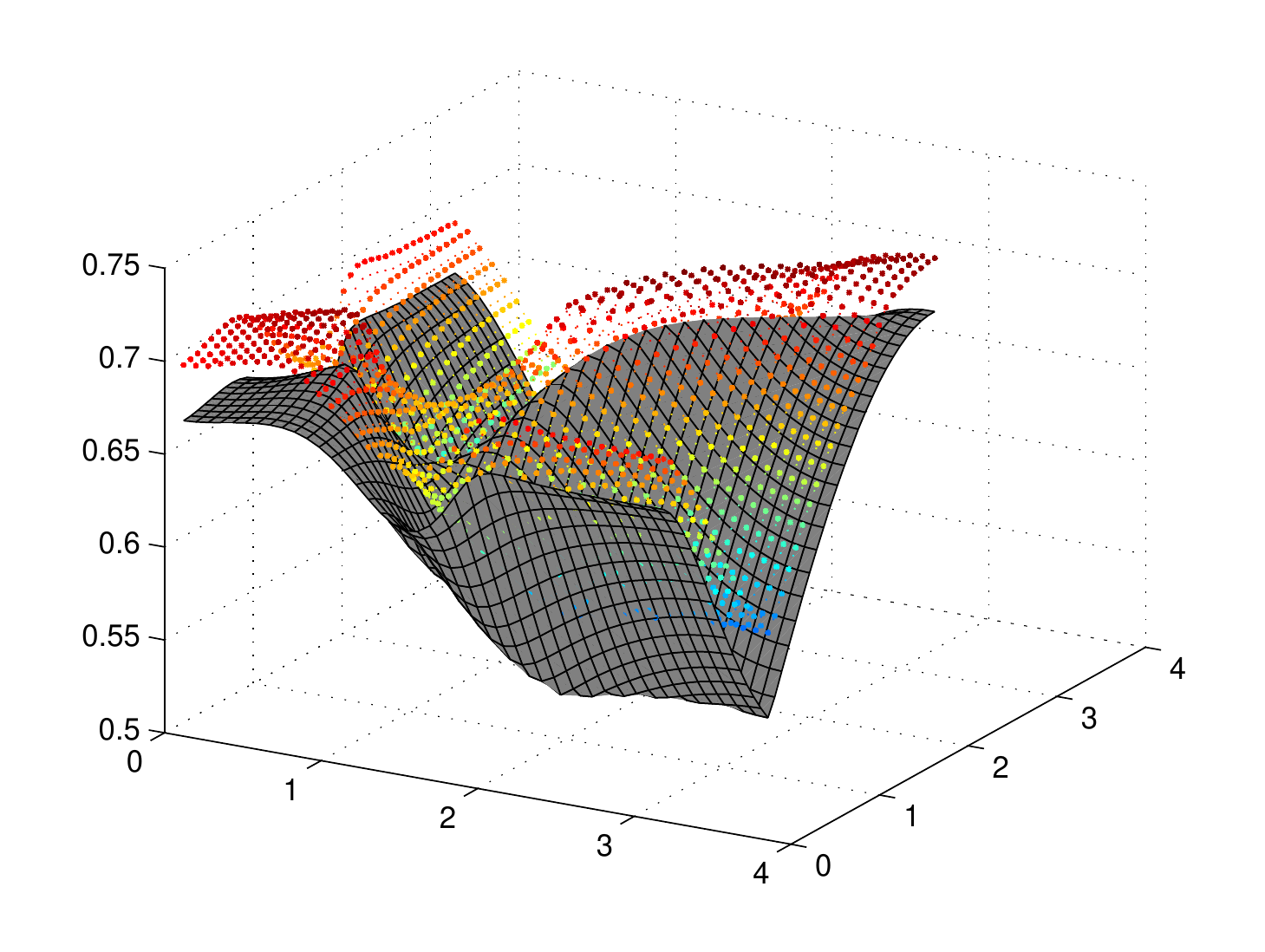}
\subcaption{$\lambda_5(f)$ }\label{approssimazioni_lambda5_n=40}
\end{subfigure}
\begin{subfigure}[c]{.30\textwidth}
\includegraphics[width=\textwidth]{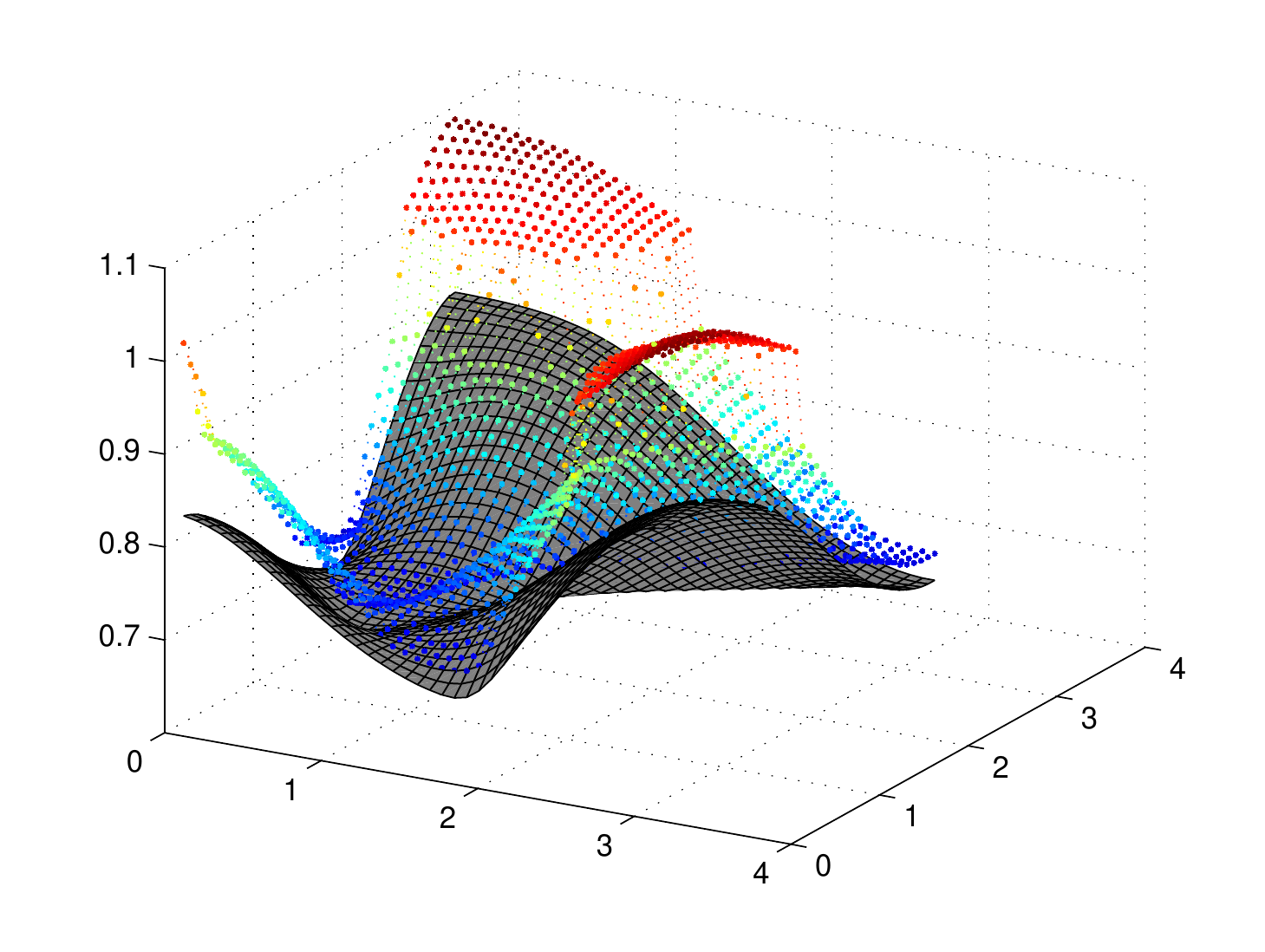}
\subcaption{$\lambda_6(f)$ }\label{approssimazioni_lambda6_n=40}
\end{subfigure}
\\
\begin{subfigure}[c]{.30\textwidth}
\includegraphics[width=\textwidth]{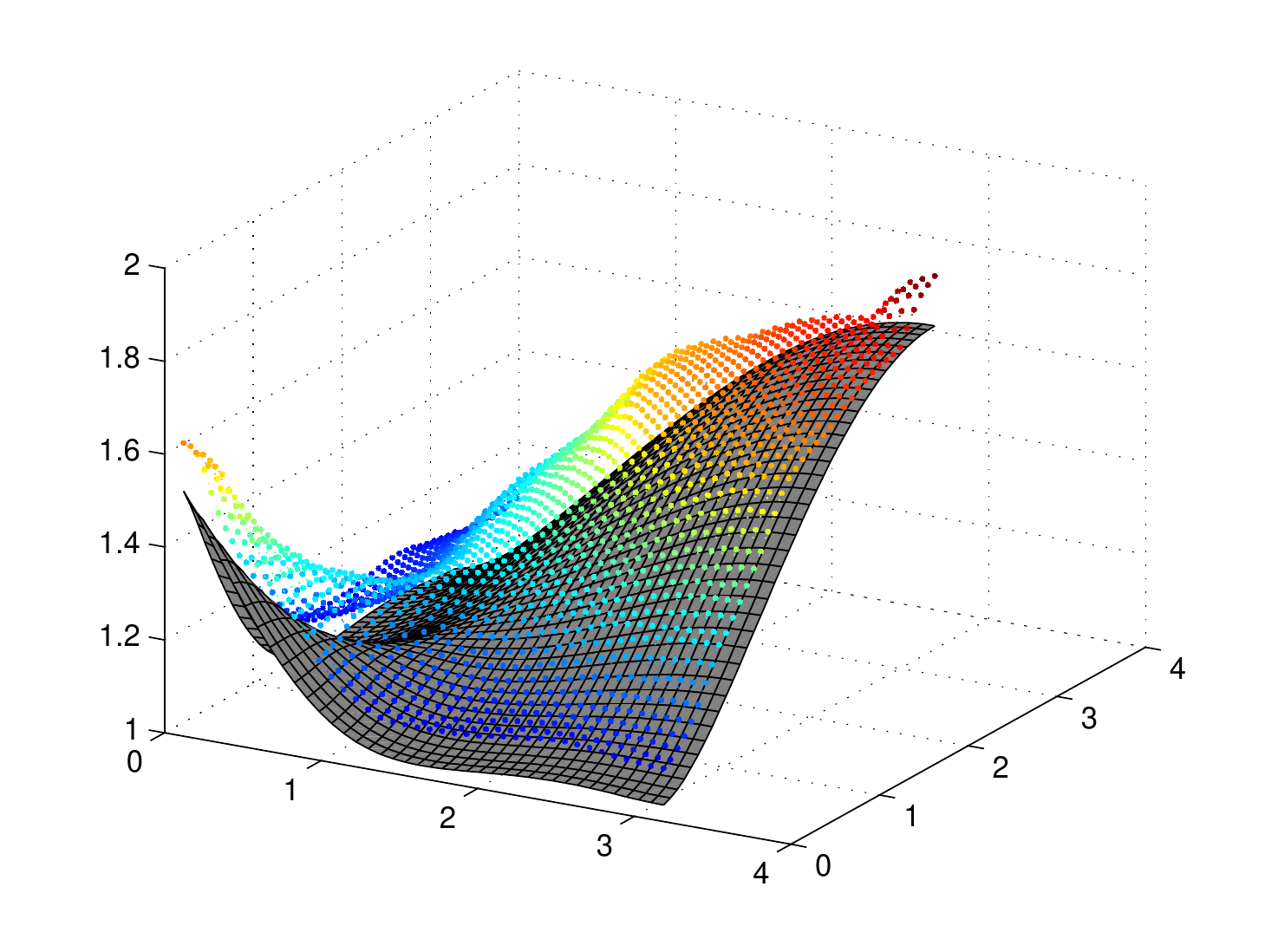}
\subcaption{$\lambda_7(f)$ }\label{approssimazioni_lambda7_n=40}
\end{subfigure}
\begin{subfigure}[c]{.30\textwidth}
\includegraphics[width=\textwidth]{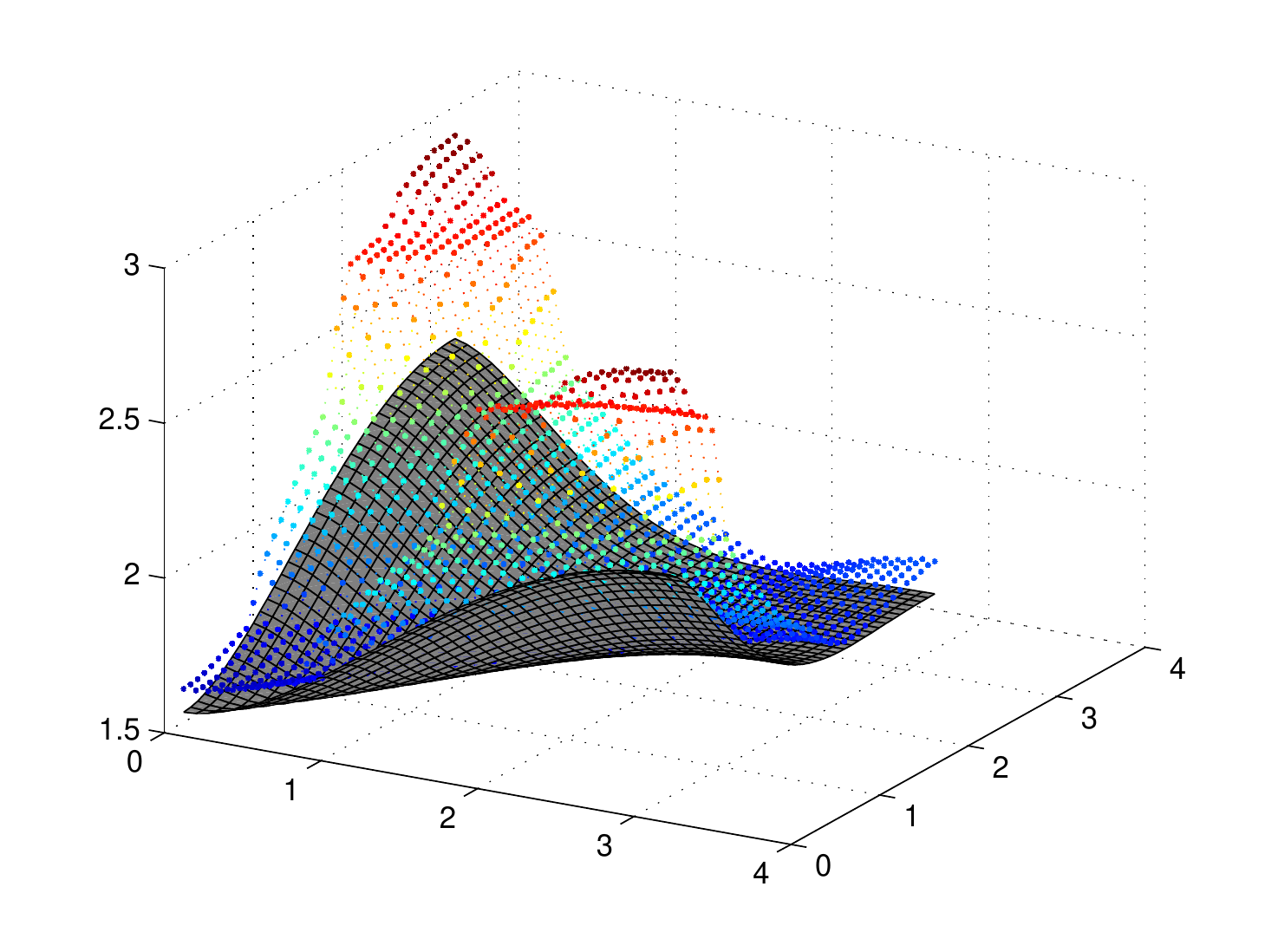}
\subcaption{$\lambda_8(f)$ }\label{approssimazioni_lambda8_n=40}
\end{subfigure}
\begin{subfigure}[c]{.30\textwidth}
\includegraphics[width=\textwidth]{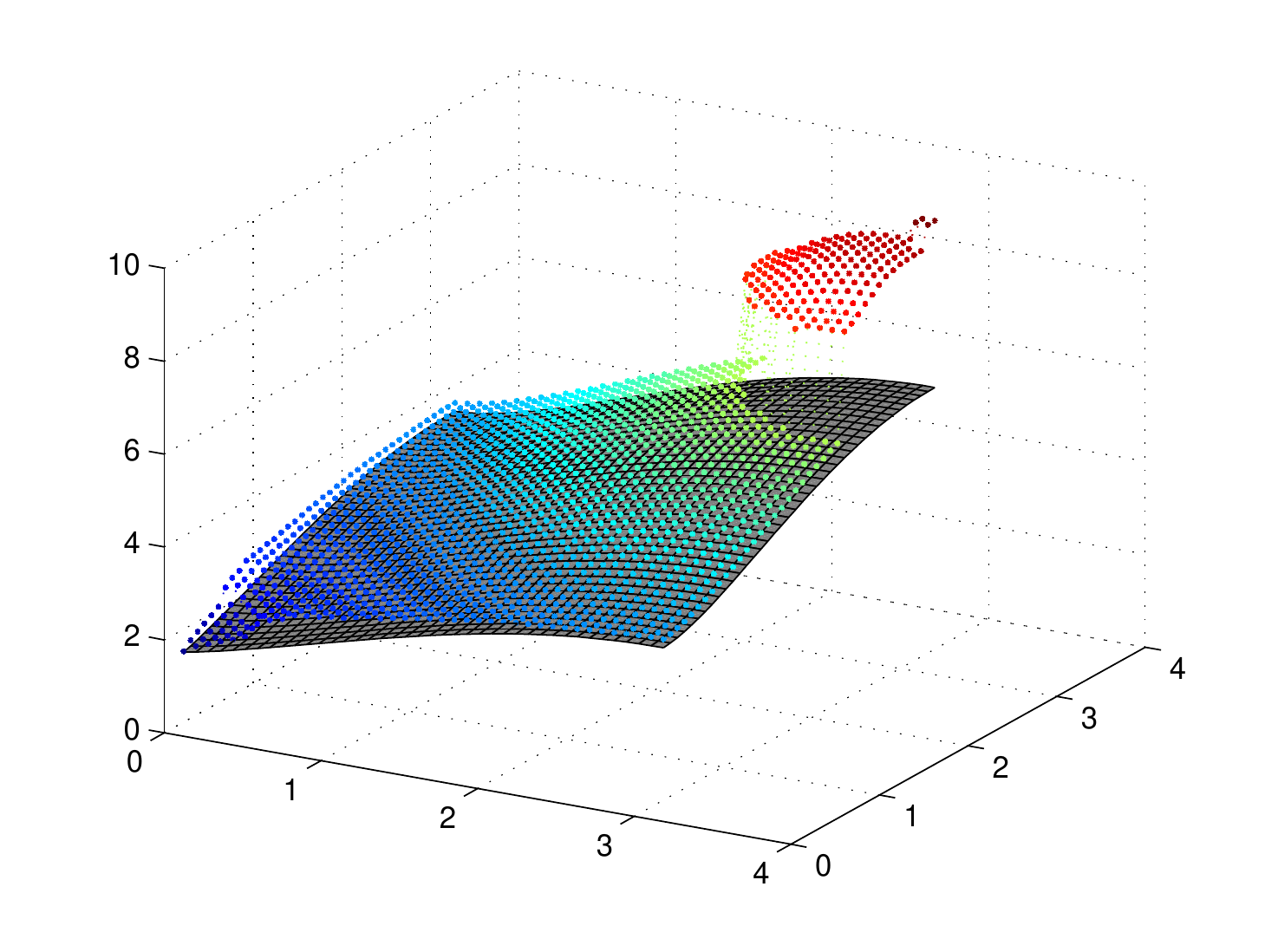}
\subcaption{$\lambda_9(f)$ }\label{approssimazioni_lambda9_n=40}
\end{subfigure}
\caption{Comparison between the eigenvalues of $K_N$ and $\lambda_l(f)$, $l=1,\ldots,9$ displayed as a mesh on $G_n$, when $n=40$}\label{approx_2d}
\end{figure}

\begin{table}[htb]
\begin{center}
\begin{tabular}{|c|c|c|c|c|c|c|}
\hline
$n$& Out. &Out./$\sqrt{9n^2}$\\
\hline\hline
10 &40&1.33\\
15 &60&1.33\\
20 &80&1.33\\
25 &100&1.33\\
30 &120&1.33 \\
35 &140&1.33\\
40 &160&1.33\\
\hline
\hline
\end {tabular}
\caption{Number of eigenvalues of $K_N$ which do not behave as the corresponding sampling of $\lambda_9(f)$.}\label{tab:outliers_9}
\end{center}
\end{table}
\end{itemize}

\subsection{A focus on the eigenvalue functions in a neighborhood of the origin}
In this subsection we study in more detail the behaviour of the eigenvalues $\lambda_l\left(f\right)$, $l=1,\ldots,9$ at $(0,0)$. Such an information is crucial when studying the convergence of a preconditioned Krylov or of a multigrid method.
Since
\begin{equation}\label{lambda_1vsall}
\lambda_1f(\theta_1,\theta_2)< \lambda_lf(\theta_1,\theta_2), \quad l=2,\dots,9, \quad  (\theta_1,\theta_2) \in \mathcal{I}^+_2,
\end{equation}
it is sufficient to study $\lambda_1(f)$ in $(0,0)$. Because of \eqref{lambda_1vsall}, the behaviour of $\lambda_1(f)$ in $(0,0)$ is equivalent to the one of
\begin{equation*}
{\rm det} f(\theta_1,\theta_2)= \prod_{i=1}^9 \lambda_i(f(\theta_1,\theta_2))
\end{equation*}
at the same point, which as a product of nonnegative functions is still a nonnegative function. We numerically checked that
\begin{align*}
{\rm det}\,f(\theta_1, \theta_2)_{|_{(0,0)}}&= 0,\\
\frac{\partial \, {\rm det} \,f(\theta_1, \theta_2)}{\partial{\theta_1}}_{|_{(0,0)}}&=\frac{\partial \, {\rm det} \,f(\theta_1, \theta_2)}{\partial{\theta_2}}_{|_{(0,0)}}=0,
\end{align*}
\begin{align*}
\frac{\partial^2 \, {\rm det} \,f(\theta_1, \theta_2)}{\partial{\theta_2}\partial{\theta_1}}_{|_{(0,0)}}&= \frac{\partial^2 \, {\rm det} \,f(\theta_1, \theta_2)}{\partial{\theta_1}\partial{\theta_2}}_{|_{(0,0)}}=0, 
\end{align*}
\begin{align*}
\frac{\partial^2 \, {\rm det} \,f(\theta_1, \theta_2)}{\partial{\theta_1}^2}_{|_{(0,0)}}&=\frac{\partial^2 \, {\rm det} \,f(\theta_1, \theta_2)}{\partial{\theta_2}^2}_{|_{(0,0)}}= \frac{53}{3912}.\\
\end{align*}
Therefore,
\begin{equation*}
\left( \nabla {\rm det} \, f(\theta_1, \theta_2) \right)_{|_{(0,0)}}=\begin{bmatrix}
\frac{\partial \, {\rm det} \,f(\theta_1, \theta_2)}{\partial{\theta_1}}_{|_{(0,0)}}\\
\frac{\partial \, {\rm det} \,f(\theta_1, \theta_2)}{\partial{\theta_2}}_{|_{(0,0)}}
\end{bmatrix}=
\begin{bmatrix}
0\\
0
\end{bmatrix},
\end{equation*}
and
\begin{equation*}
({H_{{\rm det}\, f}})_{|_{(0,0)}}=\begin{bmatrix}
\frac{\partial^2 \, {\rm det} \,f(\theta_1, \theta_2)}{\partial{\theta_1}^2}_{|_{(0,0)}} &\frac{\partial^2 \, {\rm det} \,f(\theta_1, \theta_2)}{\partial{\theta_1}\partial{\theta_2}}_{|_{(0,0)}} \\
\frac{\partial^2 \, {\rm det} \,f(\theta_1, \theta_2)}{\partial{\theta_2}\partial{\theta_1}}_{|_{(0,0)}} & \frac{\partial^2 \, {\rm det} \,f(\theta_1, \theta_2)}{\partial{\theta_2}^2}_{|_{(0,0)}}\\
\end{bmatrix}=
\begin{bmatrix}
\frac{53}{3912} &0 \\
0 & \frac{53}{3912}\\
\end{bmatrix},
\end{equation*}
that is the Hessian matrix $({H_{{\rm det}\, f}})_{|_{(0,0)}}$ is positive definite. As a consequence,
\begin{eqnarray*}
{\rm det} \,  f(\theta_1, \theta_2) = {\rm det}\,f(\theta_1, \theta_2)_{|_{(0,0)}} + (\nabla {\rm det} \, f(\theta_1,\theta_2))^T_{|_{(0,0)}}
\begin{bmatrix}
\theta_1\\
\theta_2
\end{bmatrix}+ 
\nonumber \\ 
+ \frac{1}{2}\begin{bmatrix}
\theta_1\\
\theta_2
\end{bmatrix}^T
({H_{{\rm det}\, f}})_{|_{(0,0)}}
\begin{bmatrix}
\theta_1\\
\theta_2
\end{bmatrix}+ o \left(\|\theta\|_2^2 \right), 
= \frac{53}{3912}(\theta_1^2+\theta_2^2)+o\left(\|\theta\|_2^2 \right),  
\end{eqnarray*}
where $\|\theta\|_2^2=\theta_1^2+\theta_2^2$.
%
%

\begin{figure}[htb]
\centering
\includegraphics[scale=0.3,keepaspectratio]{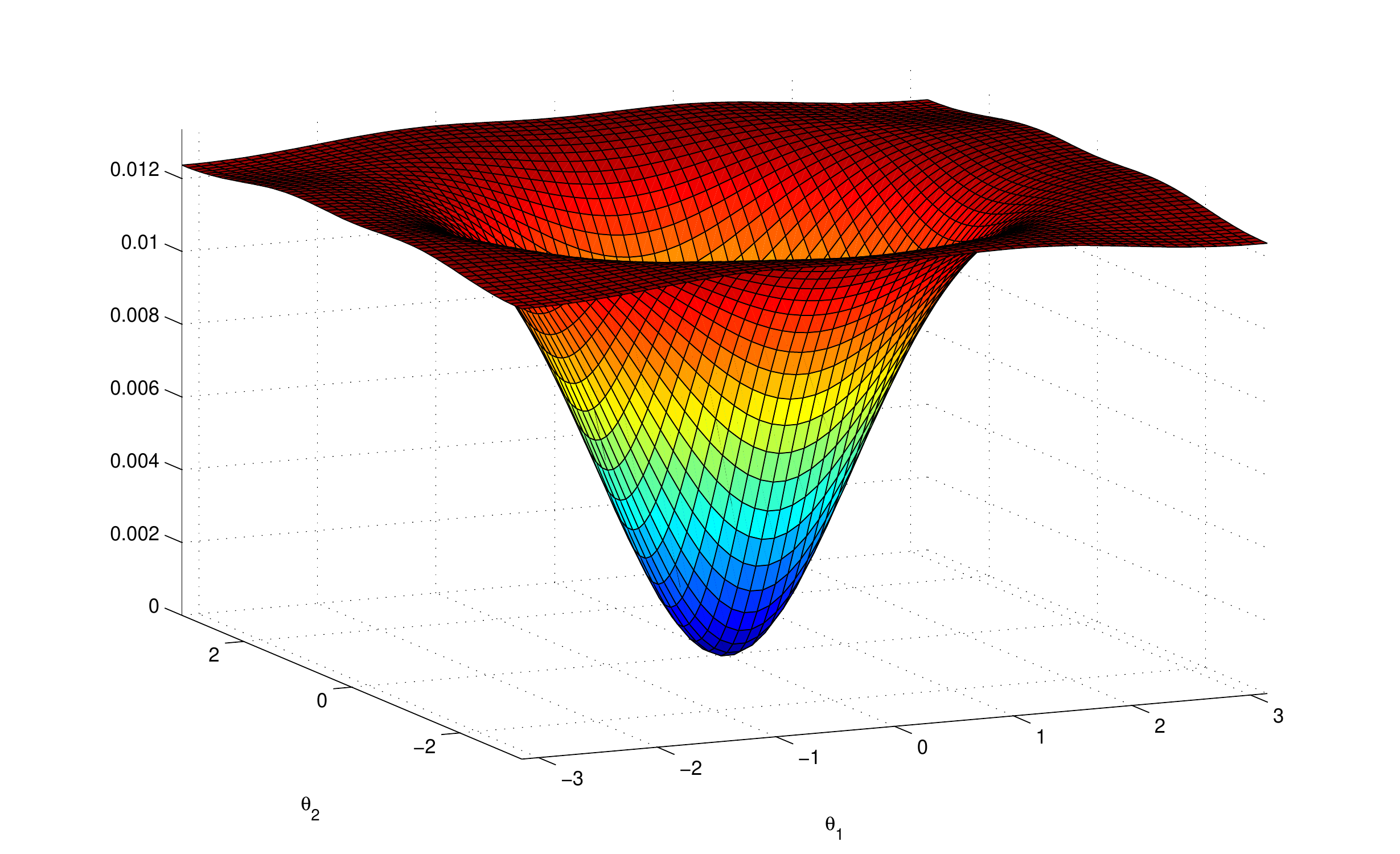}
\caption{${\rm det}\,f(\theta_1,\theta_2)$, $(\theta_1,\theta_2)\in\mathcal{I}_2$}\label{determinante_f_teta1_teta2}
\end{figure}

Hence, in a neighborhood of $(0,0)$ ${\rm det} \,  f(\theta_1, \theta_2)$ behaves as a quadratic form and
\begin{equation*}
\lim_{\|\theta\|_2\rightarrow0} \frac{{\rm det}\,f(\theta_1,\theta_2)}{\|\theta\|_2^2}=\frac{53}{3912},
\end{equation*}
which means that ${\rm det}\,f(\theta_1, \theta_2)$ and then $\lambda_1(f)$ have a zero of order $2$ in $(0,0)$, as confirmed by Figure \ref{determinante_f_teta1_teta2} and Figure \ref{valutazioni_lambda1_n=40}, respectively.

Finally, in the light of the third item of Theorem \ref{loc-extr-s}, we conclude that the minimal eigenvalue of $T_{\bf n}(f)$ goes to zero as $(\hat n)^{-1}$. 

\subsection{Spectral analysis of $K_N$ via low rank perturbations}
In this subsection we study the extremal behaviour of the matrix $K_N$, by making a careful analysis of the low rank matrix $E_\textbf{n}$, defined in Section \ref{sec:spectralsymbol}. In particular, we show that $E_\textbf{n}$ affects the number of outliers of $K_N$ but does not influence the behaviour of minimum eigenvalue of $K_N$ at $(0,0)$ with respect to that of $T_{\bf n}(f)$.\\
As shown in Section \ref{sec:spectralsymbol}, the matrix $K_N$ is the sum of two Hermitian matrices, $T_\textbf{n}(f)$ and $E_\textbf{n}$. The structural and the spectral feature of $T_\textbf{n}(f)$ have been already discussed in Section \ref{sec:spectralsymbol}, while $E_\textbf{n}$ is a block diagonal matrix with $9n \times 9n$ block diagonal blocks.
In particular there are just 3 types of non zero blocks in the matrix $E_{\textbf{n}}$.

\begin{enumerate}
\item $E^{(l)}_{\textbf{n}}$, that is in the top left corner,
\item $E^{(r)}_{\textbf{n}}$, that is in the bottom right corner,
\item $E^{(c)}_{\textbf{n}}$, that is repeated $n-2$ times in the centre of the matrix.
\end{enumerate}
We will prove that $E^{(l)}_{\textbf{n}}$, $E^{(r)}_{\textbf{n}}$ are positive definite, while $E^{(c)}_{\textbf{n}}$ is nonnegative definite. This will let us to conclude that $E_{\textbf{n}}$ is a nonnegative definite matrix.\\
Let us start by observing that $E^{(l)}_{\textbf{n}}$, $E^{(r)}_{\textbf{n}}$ and $E^{(c)}_{\textbf{n}}$ are block diagonal themselves with $9 \times 9$ diagonal blocks. In detail, 
$E^{(l)}_{\textbf{n}}$ and $E^{(r)}_{\textbf{n}}$ are composed by $n$ blocks of fixed dimension $9\times 9$, $e^{(l)}_{i}$ and $e^{(r)}_{i}$, $i=1, \dots, n$, respectively, ordered in ascending way from the top left to the bottom right. Moreover we have

\begin{equation}\label{e_l}
e^{(l)}_{i}= e_{i+1}^{(l)}, \quad i=2, \dots, n-2
\end{equation}

\begin{equation}\label{e_r}
e^{(r)}_{i}= e_{i+1}^{(r)}, \quad i=2, \dots, n-2
\end{equation}

and 
\begin{equation*}
e^{(r)}_{1}= \mathcal{J} e_n^{(l)} \mathcal{J},
\end{equation*}

\begin{equation*}
e^{(r)}_{n}= \mathcal{J} e_1^{(l)} \mathcal{J}
\end{equation*}

\begin{equation*}
e^{(r)}_{i}= \mathcal{J} e_i^{(l)} \mathcal{J}, \quad i=2, \dots, n-1
\end{equation*}

where $\mathcal{J}$ is the $9 \times 9$ flip-matrix 
\begin{equation*}
\mathcal{J}=\begin{bmatrix}
&&&&1&\\
&&&1&\\
&&\iddots&&\\
&1 &&&\\
\end{bmatrix}.
\end{equation*}
Note that $\mathcal{J}=\mathcal{J}^{-1}$, then 

\begin{equation}\label{sim1}
e^{\left(r\right)}_{1} \sim e^{\left(l\right)}_{n},
\end{equation}
\begin{equation}
e^{\left(r\right)}_{n} \sim e^{\left(l\right)}_{1},
\end{equation}
\begin{equation}\label{sim3}
e^{\left(r\right)}_{i} \sim e^{\left(l\right)}_{i}, \qquad i=2,\, \dots,\, n-1.
\end{equation}
A direct computation shows that $e^{\left(l\right)}_{1} $, $e^{\left(l\right)}_{2} $, $e^{\left(l\right)}_{n} $ are positive definite, therefore according to relations (\ref{e_l})-(\ref{e_r}) and (\ref{sim1})-(\ref{sim3}) we can conclude that $E^{(l)}_{\textbf{n}}$, $E^{(r)}_{\textbf{n}}$ are positive definite.

The matrix $E^{(c)}_{\textbf{n}}$ has only $2$ nonzero $9\times 9$  blocks, $e^{(c)}_{1}$, $e^{(c)}_{n} $ in the top left and bottom right corner respectively, such that

\begin{equation}\label{sim_c}
e^{(c)}_{n}= \mathcal{J} e_1^{(c)} \mathcal{J},
\end{equation}
while \[e^{(c)}_{i}=O_9, \quad i=2,\dots,n-1,\]

where $O_9$ is the $9\times 9$ zero matrix.\\
Because of equation (\ref{sim_c}) it holds that 
\begin{equation*}
e^{\left(c\right)}_{1} \sim e^{\left(c\right)}_{n},
\end{equation*}
then, checking directly that $e^{\left(c\right)}_1$ is positive definite, we have proved that $E^{(c)}_{\textbf{n}}$ is nonnegative definite.\\

Summarizing, since  $E^{(l)}_{\textbf{n}}$,  $E^{(r)}_{\textbf{n}}$ are positive definite, while  $E^{(c)}_{\textbf{n}}$ is nonnegative definite, we can conclude that  $E_{\textbf{n}}$ is a nonnegative definite matrix.\\
Let
\begin{equation*}
\lambda_1 (T_{\textbf{n}}(f)) \le \lambda_2 (T_{\textbf{n}}(f)) \le \dots \le \lambda_N (T_{\textbf{n}}(f))
\end{equation*}
be the eigenvalues of $T_{\textbf{n}}(f)$. Since $E_\textbf{n}$ is nonnegative definite, the Interlacing Theorem \cite{interlacing} , applied to the matrices $K_N$, $T_\textbf{n}$ and $E_\textbf{n}$, leads to the relation
\begin{equation}
\label{interlac}
\lambda_j(T_{\textbf{n}}(f))\le \lambda_j(K_N) \le \lambda_{\gamma+j}(T_{\textbf{n}}(f))
\end{equation}
for $1\le j \le N- \gamma$, where $\gamma$ is the rank of $E_{\textbf{n}}(f)$.\\
This relation is useful for the study of the conditioning of the matrix $K_N$.\\
As shown in the last subsection
\[ \lambda_{1}(T_{\textbf{n}}(f)) \overset{ \textbf{n}\rightarrow \infty}{\sim} (\hat n)^{-1},\]
and in addition, from Section \ref{sec:spectralsymbol},
$\{K_N\}_N\sim_\lambda(f,\mathcal{I}_2)$ and $\lambda_1 (f(0,0))=0$, with $f$ nonnegative definite. Hence the minimum eigenvalue of $K_N$, $ \lambda_{1}(K_N) $, has to go to zero. \\
The relation (\ref{interlac}) provides a lower bound for the convergence speed of $ \lambda_{1}(K_N) $ to zero, in fact, choosing in (\ref{interlac}) $j=1$,
\begin{equation}
\lambda_1(T_{\textbf{n}}(f))\le \lambda_1(K_N),
\end{equation}
and this implies that $ \lambda_{1}(K_N) $ does not go to zero faster than $\lambda_{1}(T_{\textbf{n}}(f))$.\\
This means that the system (\ref{system}) has the coefficient matrix $K_N$ with a better conditioning, with respect to that of the matrix $T_{\bf n}(f)$, which is quadratic with the inverse of the mesh size.
In Subsection \ref{sub:Numericaltests}, Table \ref{tab:outliers_9}, we have seen that the ratio between the number of outliers of $K_N$ with respect to the sampling of $\lambda_9(f)$ and $\sqrt{9n^2}$ is constantly equal to $\frac{4}{3}$, so the number of outliers of $K_N$ is $\frac{4}{3} \sqrt{9n^2}=4n$.\\
Due to the fact that the matrix $E_{\bf n}$ is a block diagonal matrix with precisely $2n+ 2(n-2)=4n-4$ of its $9 \times 9$ blocks positive definite, we have that $E_\textbf{n}$ has exactly 
\[9(2n)+9(2(n-2))=36n-36\]
linearly independent rows and then $\gamma$ grows exactly as $36n-36$ (see Table \ref{tab:RangKT}).\\
 This value is greater than the number of outliers, but asymptotically has the same order and the latter is in line with the theoretical forecasts induced by the Interlacing Theorem.


\begin{table}[htb]
\begin{center}
\begin{tabular}{|c|c|c|c|c|c|c|}
\hline
$n$&$Rank(E_{\bf n}(f))$\\
\hline\hline
10 & 324\\
15 & 504\\
20&684\\
25&864\\
30&1044\\
35&1224\\
40&1404\\
\hline
\hline
\end {tabular}
\caption{$Rank (E_{\bf n}(f))$ with increasing n}\label{tab:RangKT}
\end{center}
\end{table}
\subsection{Further variations}

 The numerical tests in Subsection \ref{sub:Numericaltests} are done using Dirichlet pressure boundary conditions everywhere and a standard nodal approach of conforming continuous finite elements, in order to develop the basis functions (the Lagrange interpolation polynomials passing through the given set of nodes), which are needed to compute the values in $K_N$. \\
Two simple but important changes can be considered, but their detailed analysis will be the subject of future research:
\begin{itemize}
\item using periodic boundary conditions;
\item considering another standard basis of Lagrange interpolation polynomials, passing through the Gauss-Legendre quadrature points.
\end{itemize} 
The first is motivated by the fact that several important numerical tests use this kind of boundary condition, the second one by the fact that this important kind of polynomial basis constitute an orthogonal basis. In this way the mass matrices used in the numerical method become diagonal and hence require less memory and computational effort (see e.g. \cite{Fambri2016}). 
Here we give some details on the first item.

Indeed, if we use periodic boundary conditions, then we obtain a sequence of linear systems analogous to (\ref{system}) of the form
\begin{equation}
C_Nx=b, \quad C_N\in\mathbb{R}^{N\times N}, \quad x,b\in\mathbb{R}^{N}.
\end{equation}
The symmetric matrix $ C_N\equiv C_{\bf n}(f)$ is the circulant matrix generated by the symbol $f:\mathcal{I}_2\rightarrow \mathcal{M}_{s}$, $s=(p+1)^2$, described in Section \ref{sec:spectralsymbol}\\

\begin{equation*}
f(\theta_1 , \theta_2 ) =  \hat f_{(0,0)}+  \hat f_{(-1,0)} e^{-\mathbf i \theta_1}+  \hat f_{(0,-1)} e^{-\mathbf i \theta_2}+  \hat f_{(1,0)} e^{\mathbf i \theta_1}+  \hat f_{(0,1)} e^{\mathbf i \theta_2},
\end{equation*}
Because $f$ is a trigonometric polynomial, taking into account Theorem \ref{circ-s} and Remark \ref{fourier-vs-funzione}, for $\bf n$ sufficiently large we have
\begin{equation}
\label{schur-s-bis}
C_{\bf n}(f)=(F_{\bf n}\otimes I_s) D_{\bf n}(f) (F_{\bf n}\otimes I_s)^*, 
\end{equation}
with $D_{\bf n}(f)$ as in (\ref{eig-circ-s}) and $I_s$ the $s \times s$ identity matrix.
\\
In (\ref{schur-s-bis}), as stated in Theorem \ref{circ-s}, the matrix $F_{\bf n}\otimes I_s$ is unitary and $ D_{\bf n}(f) $ is a block diagonal matrix with Hermitian blocks, $f\left(\theta_{\bf r}^{({\bf n})}\right)$, so we have
\begin{equation}
\Lambda(C_{\bf n}(f))=\left\{\lambda_l\left(f\left(\theta_{\bf r}^{({\bf n})}\right)\right)\, :\, {\bf r}= {\bf 0},\dots,{\bf n}-{\bf e}\, ; l=1\dots,s \right\},
\end{equation}
where,
for a fixed $\theta_{\bf r}^{({\bf n})}$,
 $\lambda_l\left(f\left(\theta_{\bf r}^{({\bf n})}\right)\right)$ $l=1\dots,s$ are the eigenvalues of $f\left(\theta_{\bf r}^{({\bf n})}\right)$.\\
Fixed $\bx{n}=(n_1,n_2)$, with $n_1=n_2=n$, and $p=2$ the eigenvalues of $C_N$, with $N=9n^2$, are a sampling of the eigenvalue functions $\lambda_1(f),\ldots,\lambda_9(f)$ on a equispaced grid on $\left[0,2\pi\right]^2$, 
\begin{equation*}
J_n=\left\{\left(\frac{2\pi j}{n},\frac{2\pi k}{n}\right), \quad j,k=0, \dots, n-1\right\}.
\end{equation*}
Regarding the case of a possible change of the basis functions used for representing our numerical solution, we just mention that the 
new coefficient matrix is of the form $\tilde{K}_N ={T}_{{\bf n}}(\tilde{f}) +\tilde{E}_{{\bf n}}$, with the same dimensions and structure seen in (\ref{K_N}) but with different coefficients. The symbol $\tilde{f}$ is again a trigonometric polynomial of the form described before and we obtain, with the same argument, $\{\tilde{K}_N\}_N \sim_\lambda(\tilde{f},\mathcal{I}_2)$.
However, the analytical behavior of $\tilde f$ has to be studied in detail and this will be considered in a future work.



\section{Numerical experiments}\label{sec:alg-num}

In this section we numerically verify the spectral properties derived in Section \ref{MainSec:SpectralAnalysis} on several applications of the staggered DG method \cite{Fambri2016} for the incompressible Navier-Stokes equations $\eqref{eq:intro1}$-$\eqref{eq:intro2}$. In particular we evaluate the computational effort needed for solving the main linear system for the calculation 
of the discrete pressure using successive refinements of a regular grid with $n:=n_1=n_2=\ldots =n_k$ on a square computational domain $\Omega$. From the analysis given in Section \ref{MainSec:SpectralAnalysis} we expect a condition number $\kappa=\kappa(N)\approx c N^{\frac{2}{k}}$ (the analysis has been done for $k=2$ but it is easily extendible to any $k>2$) where $k$ represents the space dimension, $N= n^k (p+1)^k$ is the matrix dimension and $c$ is a positive real constant. Due to the use of the CG method and the spectral distribution/conditioning results, the expected number of iterations for a reaching a precision $\epsilon$ can be expressed as 
\begin{eqnarray}
	Iter(n) \approx \frac{1}{2}\sqrt{c} \log\left(\frac{2 ||r_0||}{\epsilon}\right)(p+1) n \qquad k=2,3.
\label{eq:NE_1}
\end{eqnarray}
where $r_0=\p^{\nt+\delta \tau}-\p^{\nt+\delta \tau}_0$ is the initial residual between the numerical solution $\p$ at the new time step $\nt+\delta \tau$ and the initial guess for the CG method that is indicated with $\p^{\nt+\delta \tau}_0$. In particular we will use a trivial initial guess $\p^{\nt+\delta \tau}_0=b^\nt$ or a better one that is based on the solution at the previous time $\nt$, i.e. $\p^{\nt+\delta \tau}_0=\p^\nt$. In the following we will indicate with the term '\emph{IG}' this second choice for the initial guess. Furthermore, $\epsilon$ is set to $10^{-8}$ for all the simulations.
\subsection{Taylor green vortex}
\begin{figure}[htb]
\begin{subfigure}[c]{.47\textwidth}
\includegraphics[width=\textwidth]{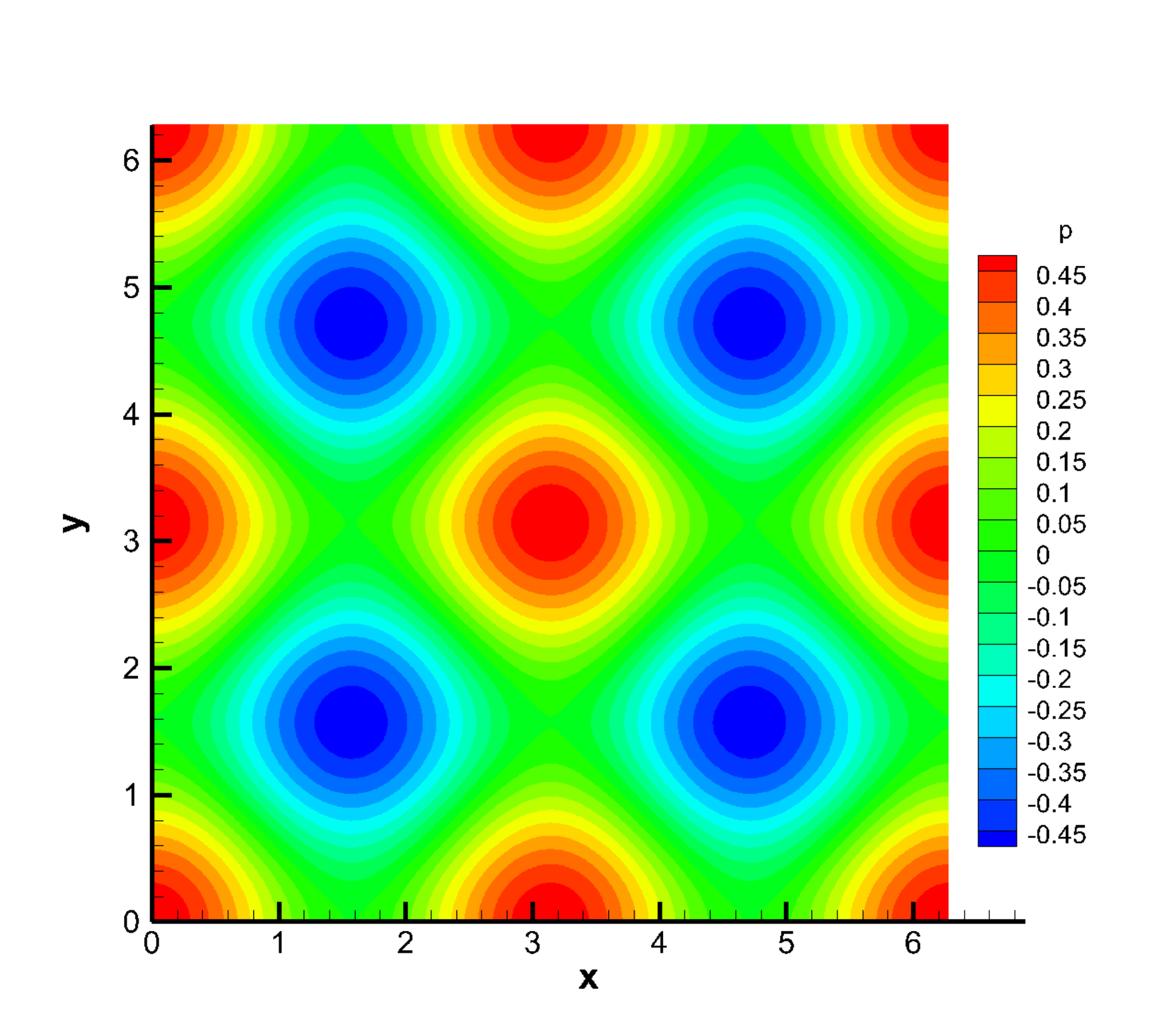}
\subcaption{Contour plot for $k=2$ }
\end{subfigure}
\begin{subfigure}[c]{.47\textwidth}
\includegraphics[width=\textwidth]{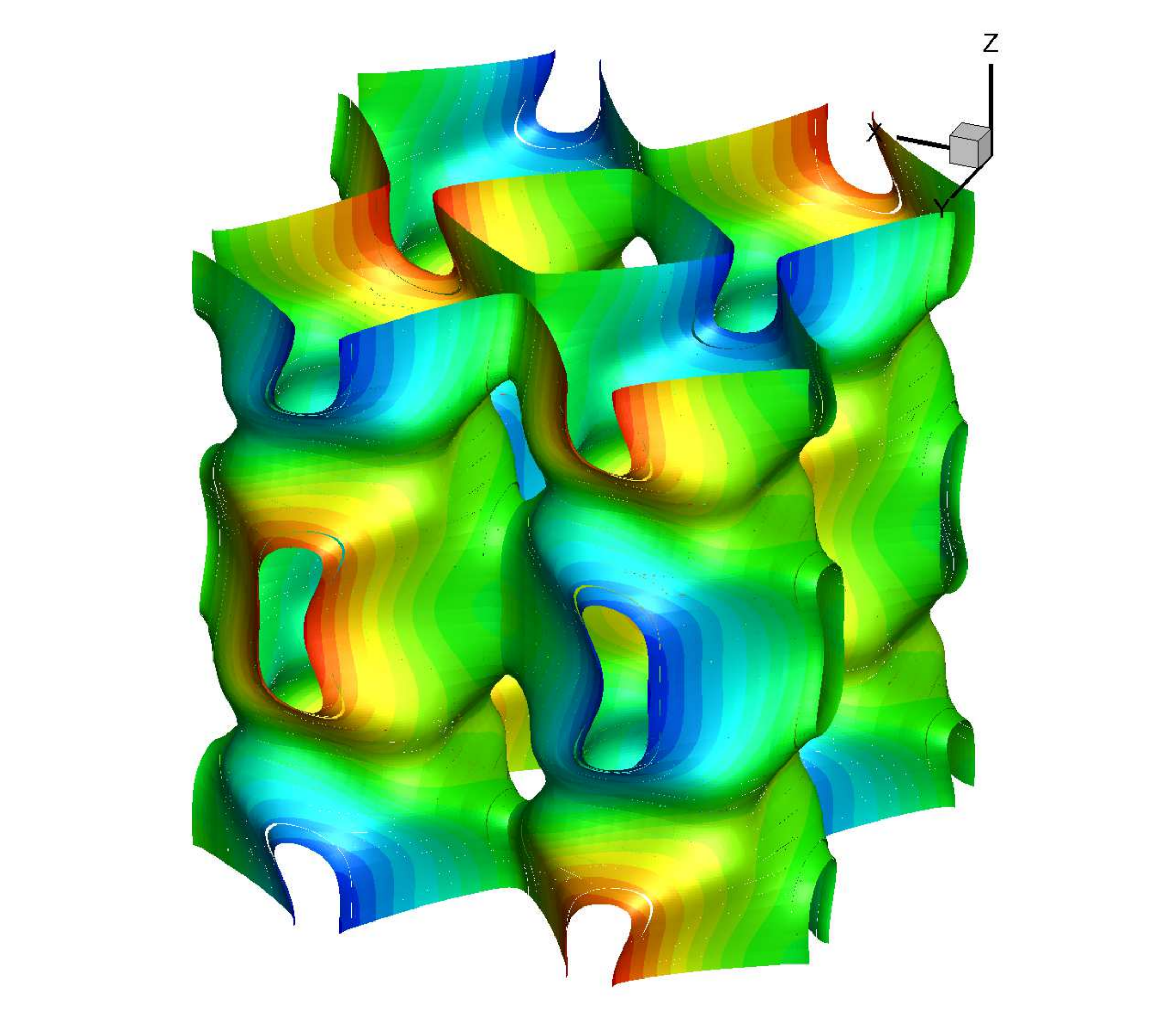}
\subcaption{Isosurface $\p=0$ for $k=3$ }
\end{subfigure}
\caption{Pressure profile at $\tau=\tau_{end}$ }\label{fig:NETGV1}
\end{figure}

\begin{figure}[htb]
\begin{subfigure}[c]{.47\textwidth}
\includegraphics[width=\textwidth]{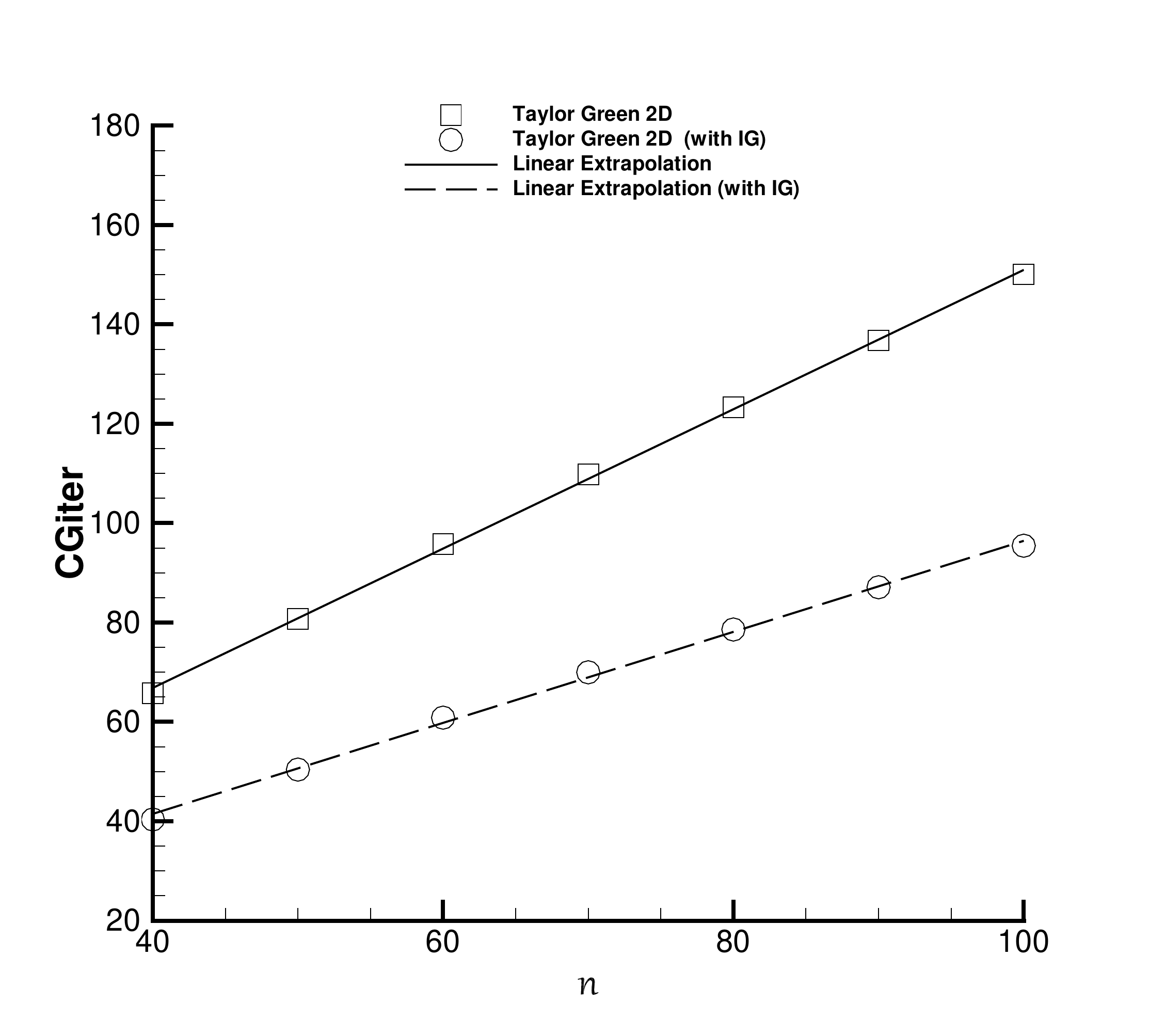}
\subcaption{$k=2$ }
\end{subfigure}
\begin{subfigure}[c]{.47\textwidth}
\includegraphics[width=\textwidth]{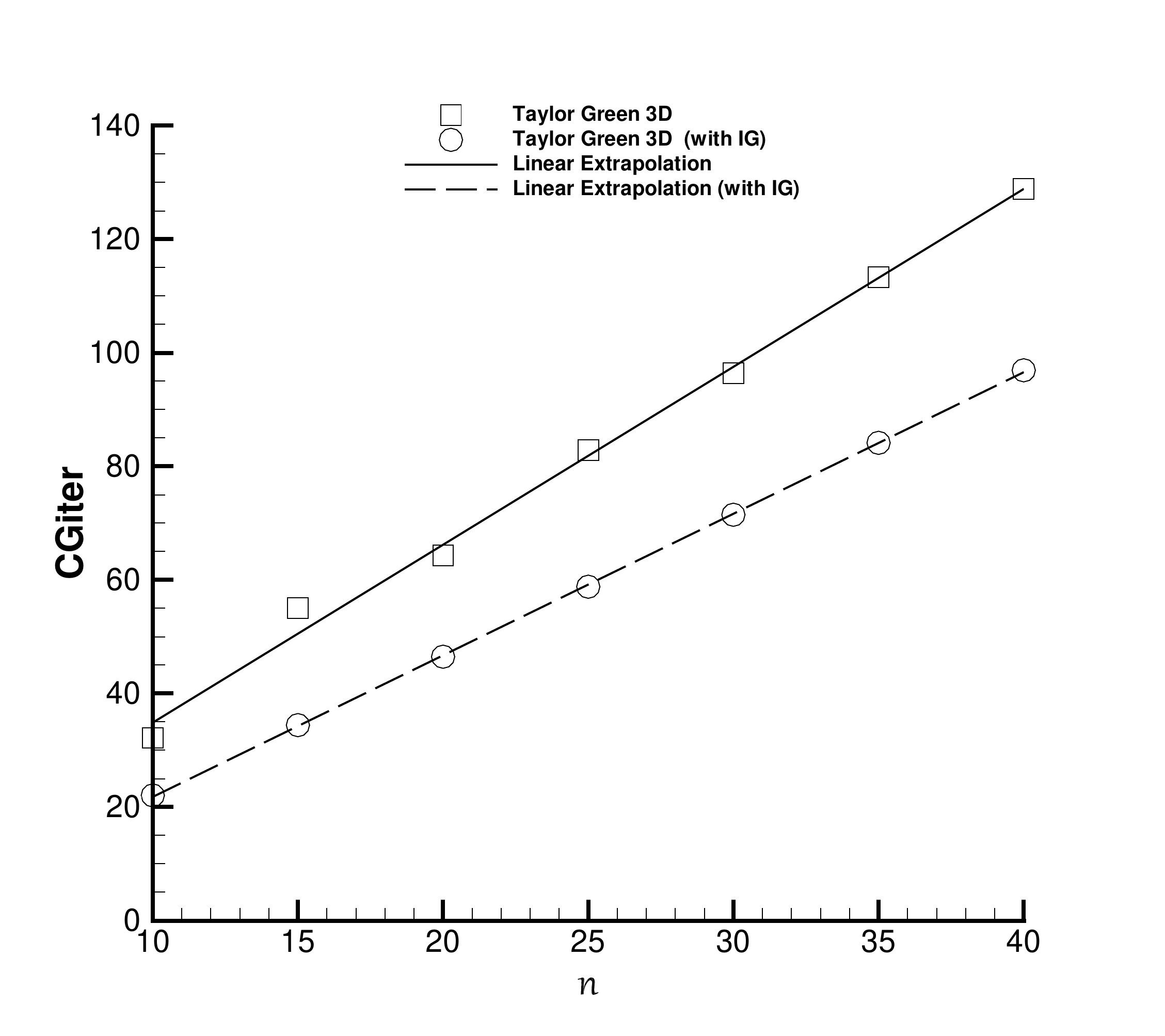}
\subcaption{$k=3$ }
\end{subfigure}
\caption{Obtained average number of iterations as a function of $n$ with and without initial guess compared with the linear extrapolation of the data.}\label{fig:NETGV2}
\end{figure}

First of all we take a classical test problem, the two and three dimensional Taylor Green vortex. The initial condition is given by
\begin{eqnarray}
\left\{
\begin{array}{l}
		u(\vec{x},0)=\sin(x)\cos(y), \\
		v(\vec{x},0)=-\cos(x)\sin(y), \\
		\pan(\vec{x},0)=\frac{1}{4}\left[\cos(2x)+\cos(2y)\right],
\end{array}
\right.
\label{eq:TGV_1_2D}
\end{eqnarray}
for $k=2$ and
\begin{eqnarray}
\left\{
\begin{array}{l}
		u(\vec{x},0)=\sin(x)\cos(y)\cos(z),\\
		v(\vec{x},0)=-\cos(x)\sin(y)\cos(z),\\
		w(\vec{x},0)= 0, \\
		\pan(\vec{x},0)=\frac{1}{16}\left[\cos(2x)+\cos(2y)\right]\left[\cos(2z)+2\right],
\end{array}
\right.
\label{eq:TGV_1_3D}
\end{eqnarray}
for $k=3$. 
The behavior of the solution for $k=3$ was numerically studied by Brachet et al in \cite{Brachet1983} and consists in a fast  generation of small scale structures, whose kinetic energy dissipation was monitored for several Reynolds numbers, see e.g. \cite{Brachet1983,Fambri2016,3STINS}. For $k=2$ and small times there is an analytical representation of the energy dissipation due to friction phenomena and hence this test can be used to check the accuracy of the numerical algorithm, see \cite{Fambri2016}. 
We consider $\Omega=[0,2\pi]^k$; $\delta \tau=5\cdot 10^{-3}$; $\tau_{end}=2$; Reynolds number $Re=800$ and periodic boundary conditions everywhere. The resulting final pressure at $\tau=\tau_{end}$ is shown in Figure $\ref{fig:NETGV1}$ for $k=2$ and $3$. The obtained average number of iterations needed to compute the solution is reported Table $\ref{tab:NETGV}$ and Figure $\ref{fig:NETGV2}$ for the two particular choices of $\p^{\nt+\delta \tau}_0=b^\nt$ and a better initial guess $\p^{\nt+\delta \tau}_0=\p^\nt$. The expected linear behavior for both two and three dimensional case is achieved according to equation \eqref{eq:NE_1}. Note that the  choice of the initial guess $\p^{\nt+\delta \tau}_0=\p^\nt$ becomes particularly good when the solution is steady or quasi-steady, since $\p^{\nt+\delta \tau}-\p^\nt \approx (K_N)^{-1}(b^\nt-b^{\nt-\delta \tau})$ and $b^\nt-b^{\nt-\delta \tau}$ contains essentially the variation of the convective-viscous contribution. Hence, for quasi stationary problems or small perturbations around a steady state, $\p^\nt$ is a good candidate for the initial guess of the CG algorithm. In practice, what we observe is indeed that the needed number of iterations tends to decrease due to a better choice of the initial guess, as suggested in equation \eqref{eq:NE_1}.   
Note, however, that the asymptotic behavior remains the same, i.e. linear in $n$, see Figure $\ref{fig:NETGV2}$ for a graphical representation.  


\begin{table}[htb]
\begin{center}
\begin{tabular}{|c|c|c|c||c|c|c|c|}
\hline
\multicolumn{4}{|c||}{$k=2$} &\multicolumn{4}{c|}{$k=3$} \\
\hline
$n$	&	$N$	& Iter& Iter with IG	&	$n$	&	$N$	& Iter& Iter with IG \\
\hline\hline
$40$ 	& $14400$		&	$65.8$	 &	$40.3$	&$10$ 	& $27000$		&	$32.2$	 &	$22.1$\\
$50$ 	& $22500$		&	$80.7$	 &	$50.3$	&$15$ 	& $91125$		&	$55.1$	 &	$34.6$\\
$60$ 	& $32400$		&	$95.8$	 &	$60.8$	&$20$ 	& $216000$	&	$64.3$	 &	$46.5$\\
$70$ 	& $44100$		&	$109.8$	 &	$69.9$	&$25$ 	& $421875$	&	$82.9$	 &	$58.8$\\
$80$ 	& $57600$		&	$123.3$	 &	$78.5$	&$30$ 	& $729000$	&	$96.4$	 &	$71.5$\\
$90$ 	& $72900$		&	$136.7$	 &	$87.0$	&$35$ 	& $1157625$	&	$113.3$	 &	$84.1$\\
$100$	& $90000$		&	$150.0$	 &	$95.4$	&$40$ 	& $1728000$	&	$128.9$	 &	$96.6$\\
\hline
\hline
\end {tabular}
\caption{Resulting average number of iterations for $\tau\in[0,2]$ with the choice of $\p^{\nt+\delta \tau}_0=b^\nt$ (left) and the use of a proper initial guess $\p^{\nt+\delta \tau}_0=\p^\nt$ (right) for $k=2,3$.}
\label{tab:NETGV}
\end{center}
\end{table}


\clearpage
\subsection{Modified double shear layer}
The previous test manifests at $\tau_{end}=2$ a relatively complex behavior for $k=3$ but a simple one involving sinusoidal functions for $k=2$. In this section we want to test the behavior of the number of iterations in a variant of the classical $2D$ double shear layer originally studied in \cite{Bell1989}. For this test case we consider the same initial condition as the one used in 
\cite{2STINS}. In the original study there is a regular jet region with $\vec{v}=(1,0)$ in a fluid with velocity $\vec{v}=(-1,0)$. The flow is characterized by two shear layers with high velocity gradient in the $y$-direction. This steady state is physically unstable due to the Kelvin-Helmholtz instability and tends to generate also in this case vortical structures close to the shear layers. In order to drive this instability, a small perturbation is introduced in the vertical velocity directly at $\tau=0$. In \cite{Bell1989} the evolution of this instability was performed for periodic boundary conditions everywhere.

\begin{figure}[htb]
\centering
\includegraphics[width=0.23\textwidth]{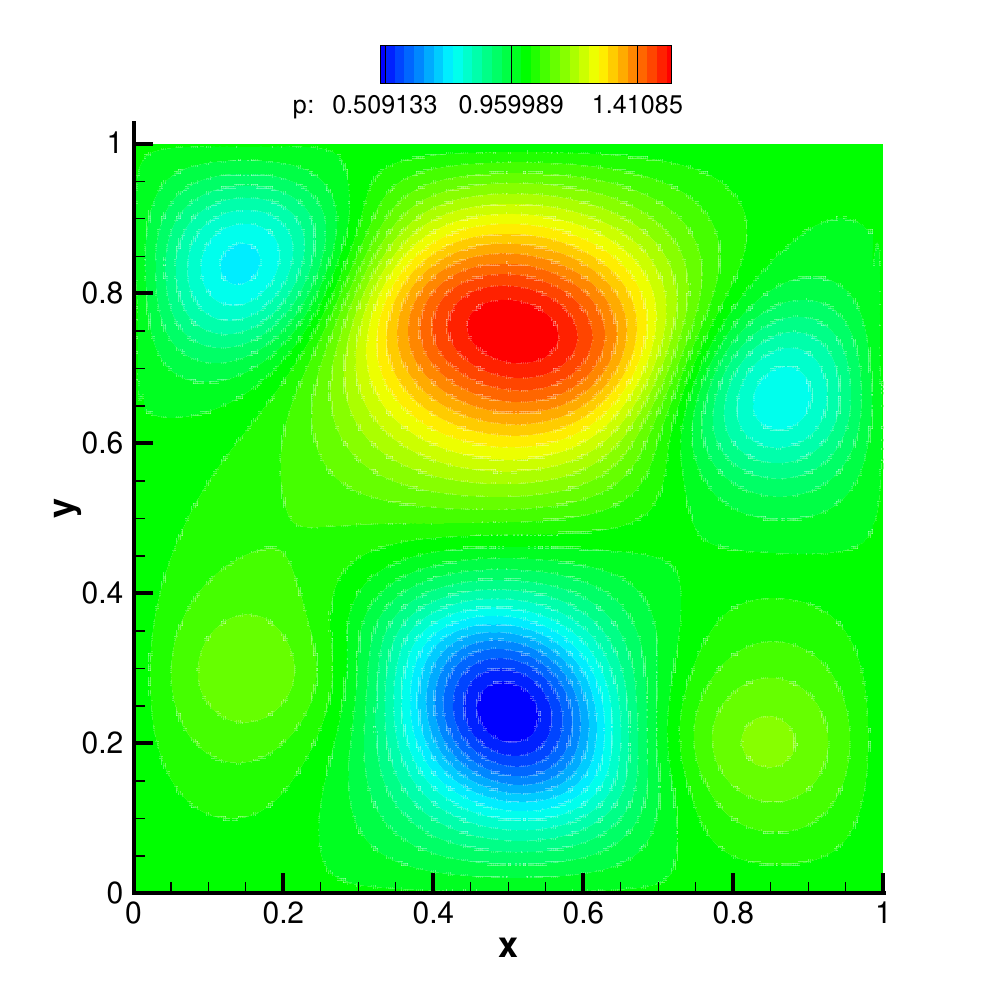}
\includegraphics[width=0.23\textwidth]{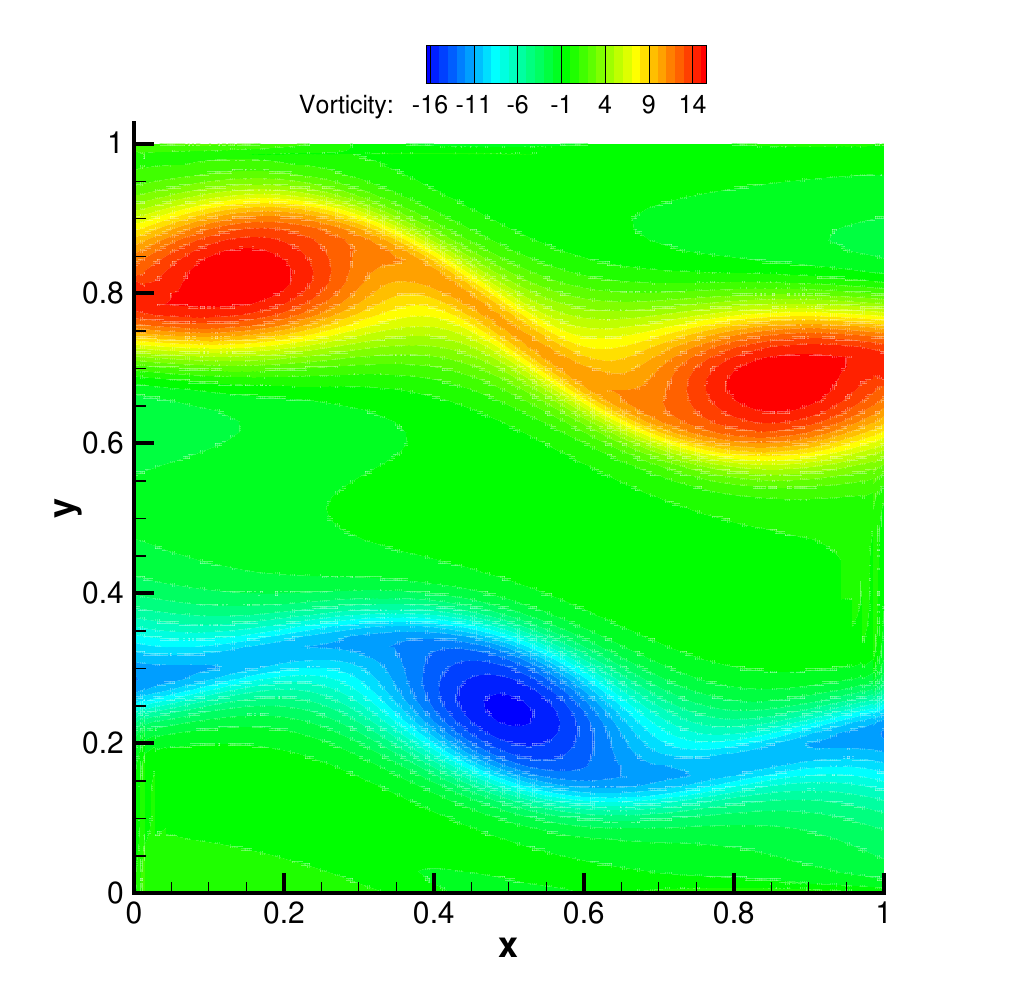}
\includegraphics[width=0.23\textwidth]{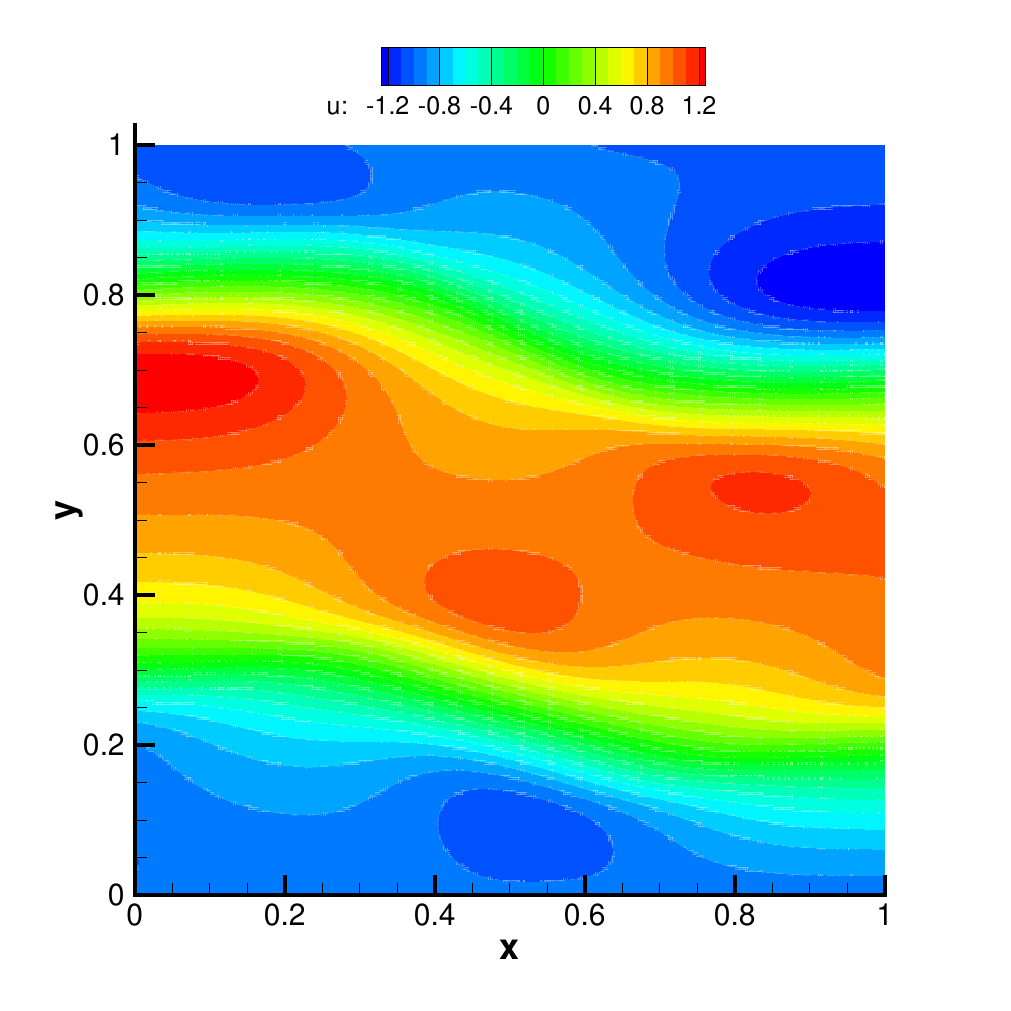}
\includegraphics[width=0.23\textwidth]{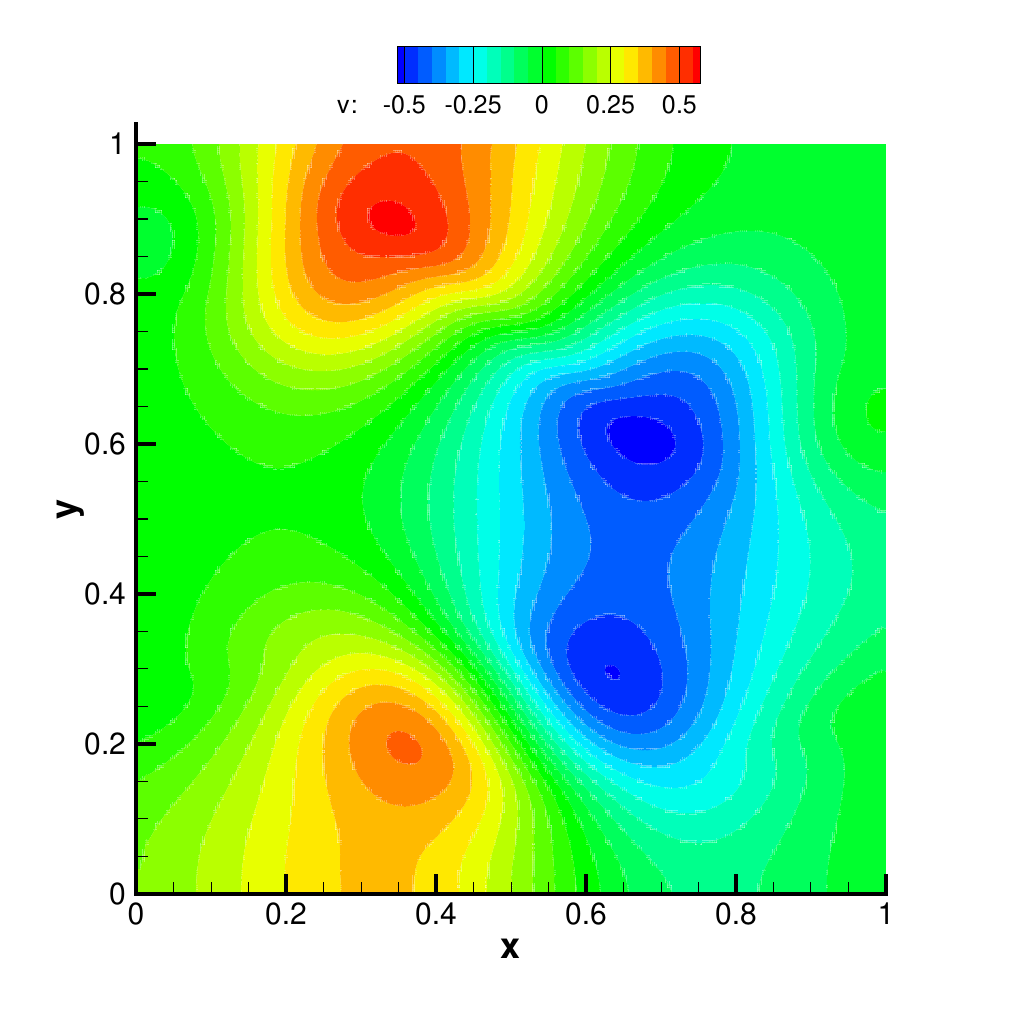}
\caption{Numerical solution for the modified double shear layer at $\tau=1$, from left to right, pressure, vorticity, $u$ and $v$ velocity component, respectively.}
\label{fig:NEDSL1}
\end{figure}

\begin{figure}[htb]
\centering
\includegraphics[width=0.47\textwidth]{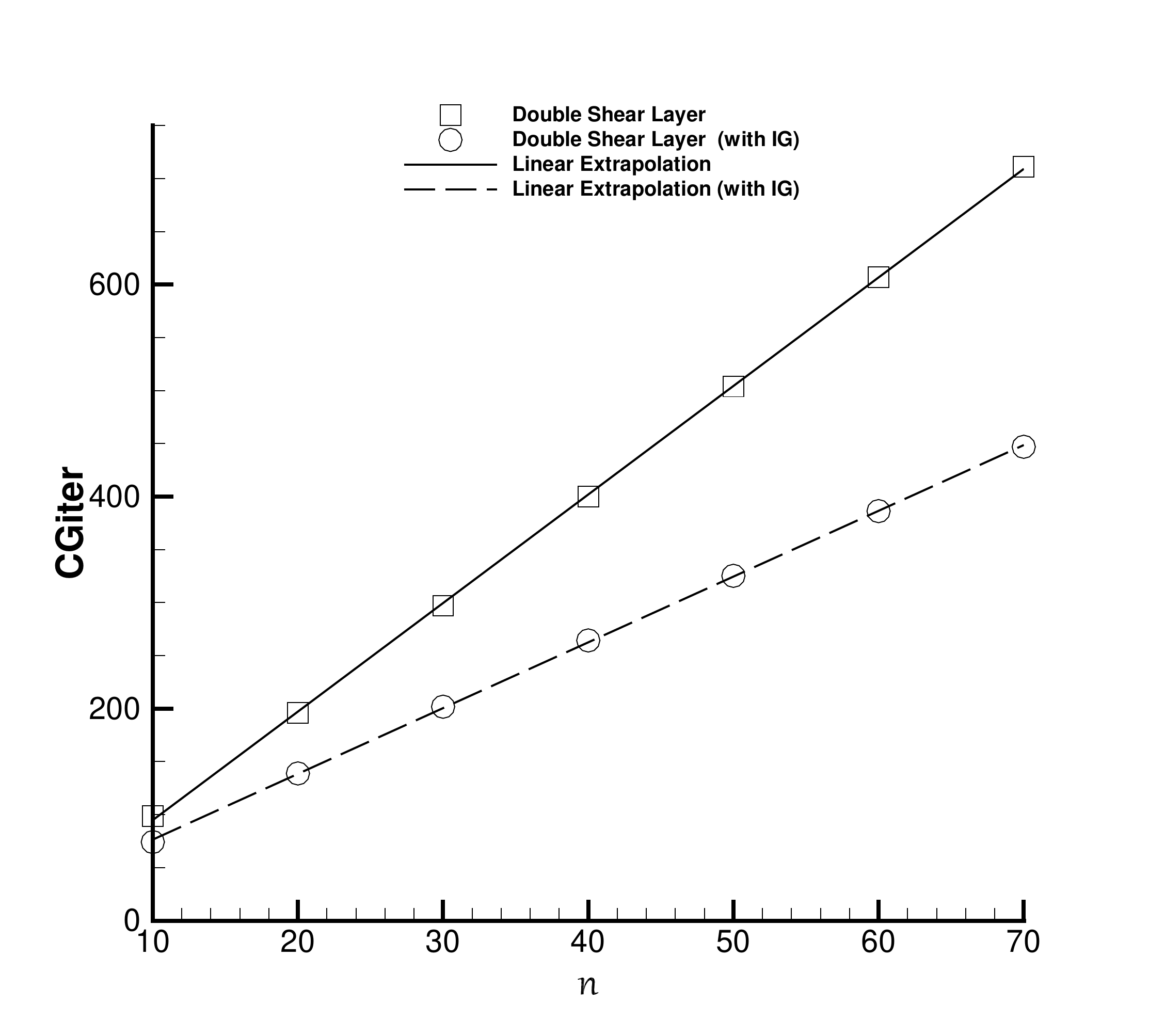}
\caption{Average number of iterations obtained in the modified double shear layer test case.}
\label{fig:NEDSL2}
\end{figure}

For this test we take $p=2$; $\tau_{end}=1$; $Re=800$ but pressure boundary condition everywhere in order to introduce the important perturbation matrix $E_{\textbf{n}}$ discussed in Section $\ref{sec:spectralsymbol}$. In this case we expect a similar but not equal behavior with respect to \cite{Bell1989} due to the use of a different kind of boundary conditions. In any case the resulting pressure field will not maintain a simple sinusoidal one for $k=2$.
The resulting numerical solution at $\tau=\tau_{end}$ for the finest grid is reported in Figure $\ref{fig:NEDSL1}$ while the obtained average number of iterations is shown in Table $\ref{tab:NEDSL}$ and the corresponding plot in Figure $\ref{fig:NEDSL2}$. As expected, also in this case the behavior for the number of iterations is linear with respect to $N^{1/k}$.

\begin{table}[htb]
\begin{center}
\begin{tabular}{|c|c|c|c|}
\hline
$n$	&	$N$	& Iter& Iter with IG \\
\hline\hline
$10$ 	& $900$			&	$98.7$	 &	$74.1$	\\
$20$ 	& $3600$		&	$195.9$	 &	$138.9$	\\
$30$ 	& $8100$		&	$297.1$	 &	$201.8$	\\
$40$ 	& $14400$		&	$400.1$	 &	$264.3$	\\
$50$ 	& $22500$		&	$504.0$	 &	$325.3$	\\
$60$ 	& $32400$		&	$607.1$	 &	$386.3$	\\
$70$	& $44100$		&	$711.3$	 &	$447.0$	\\
\hline
\hline
\end {tabular}
\caption{Resulting average number of iterations for $\tau\in[0,1]$ with and without a proper initial guess.}
\label{tab:NEDSL}
\end{center}
\end{table}

\subsection{Preconditioning}
A simple preconditioner is based on the use of the circulant matrix $C_\textbf{ n}(f)$ that is directly associated to the fully periodic boundary case. In this case we can choose as preconditioner the matrix $C_\textbf{ n}(f)$ with the Strang correction $P_\textbf{ n}(f)=C_\textbf{ n}(f)+e e^\top \frac{1}{N^2}$ where $e=(1,\ldots, 1)$ is the $N-$dimension unitary vector. The inverse of this matrix is still a circulant matrix and so its computation can be done at the cost of $c N\log N$. In this section we want to investigate the impact of this simple preconditioner on the number of iterations in the complete case where the coefficient matrix is $K_N$ (see Subsections \ref{sub:approx}, \ref{sub:spectral_distribution_KN} ).
For this test we take the same framework as in the previous numerical experiment, using $\p^\nt$ as initial guess. The resulting number of iterations is reported in Table \ref{tab:NEDSLpre}

\begin{figure}[htb]
\centering
\includegraphics[width=0.47\textwidth]{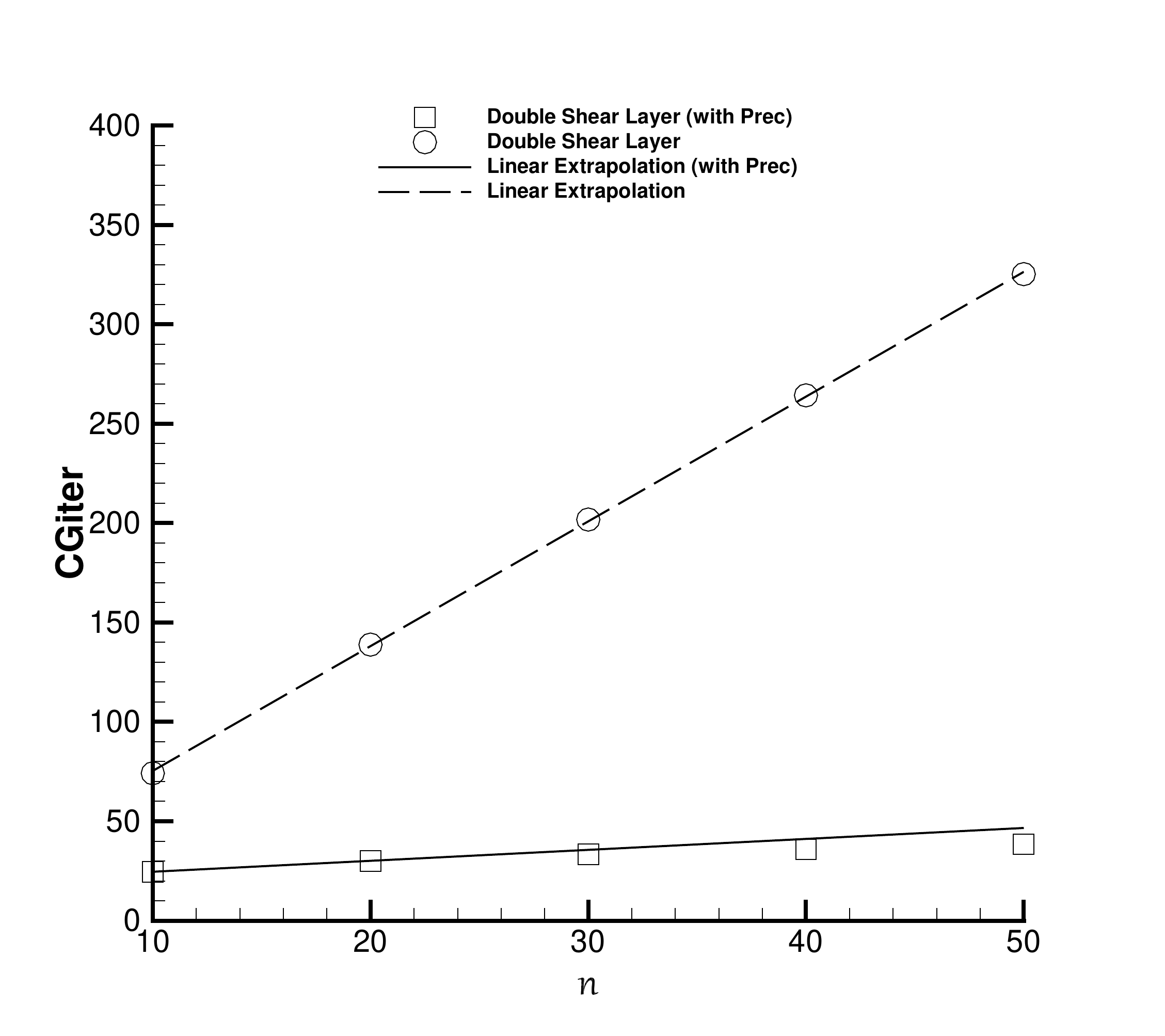}
\caption{Average number of iterations obtained in the modified double shear layer test case with and without preconditioner.}
\label{fig:DSLPre}
\end{figure}

\begin{table}[htb]
\begin{center}
\begin{tabular}{|c|c|c|c|}
\hline
$n$	&	$N$	& Iter (CG method)& Iter (PCG method) \\
\hline\hline
$10$ 	& $900$			&	$74.1$	& $24.6$	\\
$20$ 	& $3600$		&	$138.9$	& $30.1$	\\
$30$ 	& $8100$		&	$201.8$	& $33.3$	\\
$40$ 	& $14400$		&	$264.3$	& $35.8$	\\
$50$ 	& $22500$		&	$325.3$	& $38.3$	\\
\hline
\hline
\end {tabular}
\caption{Resulting average number of iterations for $\tau\in[0,1]$ without and with the preconditioner.}
\label{tab:NEDSLpre}
\end{center}
\end{table} 
The use of this preconditioner drastically reduces the number of iterations as well as the behavior that seems to be sub-linear and almost flat with respect to the case without preconditioner, see Figure $\ref{fig:DSLPre}$.

Let us now take a look at the gain in terms of CPU time obtained by the use of this simple preconditioner. Since $P_n(f)$ is a circulant matrix, we can diagonalize it as $F D F^*$ where $F=F_n \otimes F_n \otimes F_9$ is the three-level Fourier matrix and $D$ is a block diagonal matrix. We can then use the Fast Fourier Transform (FFT) to construct the matrix $D=F^* P_n(f) F$ and then $D^{-1}$ by inverting each single block. Once $D^{-1}$ is known we can easily compute $P_n^{-1}(f)x=F^* D^{-1} F x$ using the three-level FFT algorithm to compute first $x_1=Fx$ at the cost of $N\log N$. Then we have to compute $x_2=D^{-1}x_1$ at a linear cost and finally we obtain  $P_n^{-1}(f)x=F^* x_2$ again at the cost of $N\log N$. A particular test when we can really take advantage of this procedure is the fully periodic case so that the considered test becomes the classical double shear layer test case. 
The resulting total CPU time as well as the total CPU time needed to compute only the linear system is reported in Table $\ref{tab:NEDSLpre2}$ for the fully periodic case (i.e. classical double shear layer). In Table $\ref{tab:NEDSLpre3}$ we report the
obtained results for the case with pressure boundary conditions everyhere (i.e. modified double shear layer).

\begin{table}[htb]
\begin{center}
\begin{tabular}{|c|c|c c|c c|c c|c c|}
\hline
\multicolumn{2}{|c|}{} & \multicolumn{4}{c|}{No Preconditioner} & \multicolumn{4}{c|}{With Preconditioner} \\
\hline
$n$	&	$N_{step}$& $T_{TOT}$ & \small{$\frac{T_{TOT}}{N_{step}}$} & $T_{LS}$ & \small{$\frac{T_{LS}}{N_{step}}$} & $T_{TOT}$ & \small{$\frac{T_{TOT}}{N_{step}}$} & $T_{LS}$ & \small{$\frac{T_{LS}}{N_{step}}$}\\
\hline\hline
$32$ 	&	$709$		&	$195.51$	&	\small{0.28}& $126.92$	& \small{0.18}&	$77.56$	&	\small{0.11}& $6.42$	& \small{0.01}	\\
$64$ 	&	$1411$	&	$2298.6$	&	\small{1.63}& $1793.5$	& \small{1.27}&	$609.06$	&	\small{0.43}& $48.76$	& \small{0.03}	\\
$128$ &	$2829$	&	$31284.$  &	\small{11.06}& $27218.$&\small{9.62}&	$4434.81$&	\small{1.57}& $367.12$&\small{0.13}	\\
\hline
\hline
\end {tabular}
\caption{Number of time steps $N_{step}$, Total and relative (small numbers) CPU time for the solution of the main linear system for the pressure ($T_{LS}$) and the entire CPU time ($T_{TOT}$) for fully periodic boundary conditions. 
Note that in this test $p=2$, $k=2$ and $N=(p+1)^k \, n^k$.}
\label{tab:NEDSLpre2}
\end{center}
\end{table} 

\begin{table}[htb]
\begin{center}
\begin{tabular}{|c|c|c c|c c|c c|c c|}
\hline
\multicolumn{2}{|c|}{} & \multicolumn{4}{c|}{No Preconditioner} & \multicolumn{4}{c|}{With Preconditioner} \\
\hline
$n$		& $N_{step}$& $T_{TOT}$ & \small{$\frac{T_{TOT}}{N_{step}}$} & $T_{LS}$ & \small{$\frac{T_{LS}}{N_{step}}$} & $T_{TOT}$ & \small{$\frac{T_{TOT}}{N_{step}}$} & $T_{LS}$ & \small{$\frac{T_{LS}}{N_{step}}$}\\
\hline\hline
$32$ 	&	$696$		&	$371.36$	&	\small{0.53} & $296.50$	& \small{0.43} &	$219.20$ & \small{0.31}& $142.20$	& \small{0.20}	\\
$64$ 	&	$1410$	&	$4868.4$	&	\small{3.45} & $4280.4$	& \small{3.04} &	$2034.6$ & \small{1.44}& $1419.9$	& \small{1.01}	\\
$128$ &	$2853$	&	$72713.$  &	\small{25.49}& $67693.$ & \small{23.73}&	$20509.$ & \small{7.19}& $15320.$ & \small{5.37}	\\
\hline
\hline
\end {tabular}
\caption{Number of time steps $N_{step}$, total and relative (small numbers) CPU time for the solution of the main linear system for the pressure ($T_{LS}$) and the entire CPU time ($T_{TOT}$) for pressure boundary conditions everywhere. Note that in this test 
$p=2$, $k=2$ and $N=(p+1)^k \, n^k$.}
\label{tab:NEDSLpre3}
\end{center}
\end{table} 
As expected, since the symbol fully represents the periodic case, the gain on $T_{LS}$ obtained by introducing the preconditioner is impressive. In fact, {the computational cost is} \emph{essentially} {the cost of a fully explicit formula for large} $N$. Furthermore, {the number of iterations of the CG method is 1 independently on the time step and the mesh refinement}. In the worst case where we introduce pressure boundary conditions everywhere, we observe a gain factor $T_{LS}^{no pre}/T_{LS}^{pre}$ of $2.0$, $3.0$, $4.4$ for $n=32,64,128$, respectively. Hence, the advantage of using this simple preconditioner suggested by our spectral analysis is verified both for periodic and non periodic case.

\section{Conclusions}\label{sec:final}


We have considered the incompressible Navier-Stokes equations approximated by a novel family of high order semi-implicit DG methods  on staggered grids and we have studied in detail the resulting (structured) matrices.  The theory of Toeplitz matrices generated by a function (in the most general block, multi-level form) and the more recent theory of Generalized Locally Toeplitz matrix-sequences have been the key tools for analyzing the spectral properties of the considered large matrices. We have obtained a quite complete picture of the spectral properties of the underlying linear systems that result after the discretization of the PDE. This information  has been employed for giving a forecast of the convergence history of the CG method. Several numerical tests are provided and critically illustrated in order to show the validity and the potential of our analysis.

The use of these results will be the ground for further researche in the direction of new more advanced techniques (involving preconditioning, multigrid, multi-iterative solvers \cite{Smulti}), by taking into account variable coefficients, unstructured 
meshes in geometrically complex domains, different basis functions and various boundary conditions: we will develop these research lines in future work.

	
\section*{Acknowledgments}
\noindent The work has been partially supported by INDAM-GNCS.
\newline
Maurizio Tavelli, Michael Dumbser and Francesco Fambri have been financed 
by the European Research Council (ERC) under the European Union’s Seventh 
Framework Programme (FP7/2007-2013) with the research project STiMulUs, 
ERC Grant agreement no. 278267

\bibliographystyle{siamplain}
\bibliography{SIDG}

\appendix

\section{Matrix symbol for $k=2$ and $p=2$} 
\label{app.symbol.k2p2}

Recall that for the two-dimensional case ($k=2$) the matrix symbol $f$ is given according to
\eqref{eqn.symbol.f} by 
\begin{equation*}
f(\theta_1 , \theta_2 ) =  \hat f_{(0,0)} +  \hat f_{(-1,0)} e^{-\mathbf i \theta_1}+  \hat f_{(0,-1)} e^{-\mathbf i \theta_2}+  \hat f_{(1,0)} e^{\mathbf i \theta_1}+  \hat f_{(0,1)} e^{\mathbf i \theta_2}. 
\end{equation*}
For the special case $p=2$, the matrices appearing in the above expression (see \cite{Fambri2016} for details concerning their definition) read 
\begin{equation*}
\hat f_{(0,0)} =\left( \begin{array}{ccccccccc} 
\frac{  127}{360    }&\frac{   41}{480   }&\frac{   -43}{320   }&\frac{    41}{480   }&\frac{    -1}{360  }&\frac{     -2}{45}   & \frac{-43}{320}     &  \frac{-2}{45     }&\frac{   13}{288 } \\ 
  &         &    &  &  &    &  & &\\ 
      \frac{41}{480    }&\frac{  103}{90    }&\frac{    41}{480   }&\frac{    -1}{360   }&\frac{     5}{24   }&\frac{     -1}{360}  &\frac{  -2}{45  }    &\frac{  -113}{240  }&\frac{     -2}{45} \\
        &         &    &  &  &    &  & &\\ 
\frac{     -43}{320    }&\frac{   41}{480   }&\frac{   127}{360   }& \frac{-2}{45    }&\frac{    -1}{360  }&\frac{     41}{480}  & \frac{ 13}{288    }&\frac{   -2}{45   }&\frac{    -43}{320}\\ 
  &         &    &  &  &    &  & &\\  
   \frac{ 41}{480    }&\frac{   -1}{360   }&\frac{    -2}{45    }&\frac{   103}{90    }&\frac{     5}{24   }&\frac{   -113}{240}  &\frac{  41}{480    }&\frac{   -1}{360   }&\frac{    -2}{45} \\   
     &         &    &  &  &    &  & &\\ 
    \frac{  -1}{360    }&\frac{    5}{24    }&\frac{    -1}{360   }&\frac{     5}{24    }&\frac{   158}{45   }&\frac{      5}{24 }  &\frac{  -1}{360     }&\frac{   5}{24    }&\frac{    -1}{360}\\  
      &         &    &  &  &    &  & &\\ 
 \frac{     -2}{45     }&\frac{   -1}{360   }&\frac{    41}{480   }&\frac{  -113}{240   }&\frac{     5}{24   }&\frac{    103}{90 }  &\frac{ -2}{45      }&\frac{  -1}{360    }&\frac{   41}{480} \\  
   &         &    &  &  &    &  & &\\  
 \frac{    -43}{320    }&\frac{   -2}{45    }&\frac{    13}{288   }&\frac{    41}{480   }&\frac{    -1}{360  }&\frac{     -2}{45 }  & \frac{127}{360    }&\frac{   41}{480   }&\frac{   -43}{320} \\ 
   &         &    &  &  &    &  & &\\  
    \frac{  -2}{45 }     & \frac{ -113}{240   }&\frac{    -2}{45    }&\frac{    -1}{360   }&\frac{     5}{24   }&\frac{     -1}{360 } &\frac{ 41}{480  }    &  \frac{103}{90     }&\frac{   41}{480} \\  
      &         &    &  &  &    &  & &\\  
 \frac{     13}{288    }&\frac{   -2}{45    }&\frac{   -43}{320   }&\frac{    -2}{45    }&\frac{    -1}{360  }&\frac{     41}{480 } & \frac{ -43}{320    }&\frac{   41}{480  }&\frac{    127}{360 }
\end{array} \right);
\end{equation*} 
\begin{equation*}
\hat f_{(-1,0)} =\left( \begin{array}{ccccccccc} 
   \frac{ 5}{288    }&\frac{    5}{576 }&\frac{      -5}{1152  }&\frac{    23}{720   }&\frac{    23}{1440  }&\frac{   -23}{2880   }&\frac{  -11}{1440 }&\frac{-11}{2880 }&\frac{     11}{5760} \\ 
     &         &    &  &  &    &  & &\\  
\frac{       5}{576  }&\frac{      5}{72 }&\frac{        5}{576 }&\frac{      23}{1440}&\frac{      23}{180 }&\frac{      23}{1440 }&\frac{ -11}{2880  }&\frac{-11}{360   }&\frac{   -11}{2880} \\
  &         &    &  &  &    &  & &\\ 
\frac{      -5}{1152  }&\frac{     5}{576 }&\frac{       5}{288 }&\frac{     -23}{2880 }&\frac{     23}{1440 }&\frac{     23}{720  }&\frac{  11}{5760  }&\frac{-11}{2880  }&\frac{   -11}{1440} \\  
  &         &    &  &  &    &  & &\\ 
\frac{     -17}{144   }&\frac{   -17}{288 }&\frac{      17}{576 }&\frac{     -47}{360  }&\frac{    -47}{720   }&\frac{    47}{1440 }&\frac{  23}{720   }&\frac{  23}{1440 }&\frac{    -23}{2880} \\
  &         &    &  &  &    &  & &\\   
\frac{     -17}{288  }&\frac{    -17}{36   }&\frac{    -17}{288 }&\frac{     -47}{720  }&\frac{    -47}{90    }&\frac{   -47}{720   }&\frac{ 23}{1440  }&\frac{ 23}{180   }&\frac{    23}{1440} \\ 
  &         &    &  &  &    &  & &\\  
\frac{      17}{576   }&\frac{  -17}{288   }&\frac{   -17}{144  }&\frac{     47}{1440  }&\frac{   -47}{720   }&\frac{   -47}{360   }&\frac{ -23}{2880   }&\frac{23}{1440  }&\frac{    23}{720} \\  
  &         &    &  &  &    &  & &\\  
\frac{      -7}{288    }&\frac{   -7}{576  }&\frac{      7}{1152}&\frac{     -17}{144  }&\frac{    -17}{288  }&\frac{     17}{576  }&\frac{   5}{288   }&\frac{5}{576    }&\frac{   -5}{1152} \\  
  &         &    &  &  &    &  & &\\ 
    \frac{  -7}{576     }&\frac{  -7}{72   }&\frac{     -7}{576 }&\frac{     -17}{288  }&\frac{    -17}{36    }&\frac{   -17}{288   }&\frac{ 5}{576    }&\frac{5}{72     }&\frac{    5}{576}  \\  
      &         &    &  &  &    &  & &\\ 
    \frac{   7}{1152     }&\frac{ -7}{576  }&\frac{     -7}{288 }&\frac{      17}{576  }&\frac{    -17}{288   }&\frac{   -17}{144   }&\frac{ -5}{1152  }&\frac{ 5}{576   }&\frac{     5}{288} \\ 
\end{array} \right);
\end{equation*}
 \begin{equation*}
    \hat f_{(0,-1)}  =\left( \begin{array}{ccccccccc} 
      \frac{ 5}{288  }&\frac{  23}{720   }&\frac{   -11}{1440  }&\frac{     5}{576  }&\frac{     23}{1440  }&\frac{   -11}{2880  }&\frac{    -5}{1152}&\frac{-23}{2880  }&\frac{    11}{5760} \\  
        &         &    &  &  &    &  & &\\ 
\frac{     -17}{144  }&\frac{ -47}{360   }&\frac{    23}{720   }&
\frac{    -17}{288  }&\frac{    -47}{720   }&\frac{    23}{1440  }&
\frac{    17}{576 }& \frac{ 47}{1440  }&\frac{   -23}{2880}  \\ 
  &         &    &  &  &    &  & &\\ 
   \frac{    -7}{288  }&\frac{ -17}{144    }&\frac{    5}{288   }&\frac{    -7}{576  }&\frac{    -17}{288   }&\frac{     5}{576   }&\frac{     7}{1152}&\frac{ 17}{576   }&\frac{    -5}{1152} \\  
     &         &    &  &  &    &  & &\\ 
   \frac{     5}{576  }&\frac{  23}{1440  }&\frac{   -11}{2880  }&\frac{     5}{72   }&\frac{     23}{180   }&\frac{   -11}{360   }&\frac{     5}{576 }&\frac{ 23}{1440  }&\frac{   -11}{2880} \\  
     &         &    &  &  &    &  & &\\ 
  \frac{    -17}{288  }&\frac{ -47}{720    }&\frac{   23}{1440  }&\frac{   -17}{36   }&\frac{    -47}{90    }&\frac{    23}{180   }&\frac{   -17}{288 }&\frac{ -47}{720   }&\frac{    23}{1440} \\  
    &         &    &  &  &    &  & &\\ 
  \frac{     -7}{576  }&\frac{ -17}{288    }&\frac{    5}{576   }&\frac{    -7}{72   }&\frac{    -17}{36    }&\frac{     5}{72    }&\frac{    -7}{576 }&\frac{ -17}{288   }&\frac{     5}{576}\\  
    &         &    &  &  &    &  & &\\ 
   \frac{    -5}{1152 }&\frac{ -23}{2880   }&\frac{   11}{5760  }&\frac{     5}{576  }&\frac{     23}{1440  }&\frac{   -11}{2880  }&\frac{     5}{288 }&\frac{  23}{720   }&\frac{   -11}{1440} \\  
     &         &    &  &  &    &  & &\\ 
   \frac{    17}{576  }&\frac{  47}{1440   }&\frac{  -23}{2880  }&\frac{   -17}{288  }&\frac{    -47}{720   }&\frac{    23}{1440  }&\frac{   -17}{144 }&\frac{  -47}{360  }&\frac{     23}{720} \\ 
     &         &    &  &  &    &  & &\\   
    \frac{    7}{1152 }&\frac{  17}{576    }&\frac{   -5}{1152  }&\frac{    -7}{576  }&\frac{    -17}{288   }&\frac{     5}{576   }&\frac{    -7}{288 }&\frac{ -17}{144   }&\frac{     5}{288 }
       \end{array} \right);
\end{equation*}
\begin{equation*}
   \hat f_{(1,0)} = \hat f_{(-1,0)}^T;
\end{equation*} 
\begin{equation*}
   \hat f_{(0,1)} =   \hat f_{(0,-1)}^T.
\end{equation*}

\end{document}